\documentclass[12pt]{amsart}
\usepackage{amscd,amsfonts,amsmath,amsxtra,amssymb}
\usepackage[all]{xy}
\usepackage{graphics}
\usepackage{hyperref}
\usepackage[french]{babel}
\usepackage[T1]{fontenc}

\newtheorem{sub}{}[section]
\newtheorem{subsub}{}[sub]
\font\tte=cmbsy10

\def\kov#1{\overline{#1}}

\def\coker{\mathop{\rm coker}\nolimits}
\def\hom{\mathop{\rm Hom}\nolimits}
\def\ext{\mathop{\rm Ext}\nolimits}
\def\Aut{\mathop{\rm Aut}\nolimits}
\def\End{\mathop{\rm End}\nolimits}
\def\imm{\mathop{\rm Im}\nolimits}
\def\lra{\longrightarrow}
\def\som{\mathop{\hbox{$\displaystyle\bigoplus$}}\limits}
\def\sigg{\mathop{\hbox{$\displaystyle\sum$}}\limits}

\def\timex{\setbox250=\hbox{$\times$}\hskip 5pt\Box\hskip -0.75em
{\raise 1.5pt\vbox{\box250}}\hskip 5pt}
\def\paragra{{\tte \char120}}
\def\para{\paragra~\hskip -2pt}
\def\subsubsubsection#1#2{\noindent{\bf #1} \ \ {\it #2}}
\def\hfl#1#2{\smash{\mathop{\hbox to 12mm{\rightarrowfill}}
\limits^{\scriptstyle#1}_{\scriptstyle#2}}}
\def\vfl#1#2{\llap{$\scriptstyle #1$}\left\downarrow
\vbox to 6mm{}\right.\rlap{$\scriptstyle #2$}}

\def\m#1{{\hbox{$#1$}}}
\def\ot{\otimes}
\def\sepsec{\vskip 1.4cm}
\def\Nligne{\hfil\break}
\newcommand{\N}{{\mathbb N}}
\newcommand{\A}{{\mathbb A}}
\newcommand{\Z}{{\mathbb Z}}

\newcommand{\C}{{\mathbb C}}
\newcommand{\Q}{{\mathbb Q}}
\renewcommand{\P}{{\mathbb P}}
\newcommand{\D}{{\mathbb D}}

\renewcommand{\L}{{\mathbb L}}

\newcommand{\ke}{{\mathcal E}}
\newcommand{\kf}{{\mathcal F}}
\newcommand{\kg}{{\mathcal G}}
\newcommand{\kh}{{\mathcal H}}
\newcommand{\ki}{{\mathcal I}}

\newcommand{\ko}{{\mathcal O}}

\newcommand{\ks}{{\mathcal S}}
\newcommand{\kt}{{\mathcal T}}
\newcommand{\ku}{{\mathcal U}}
\newcommand{\kv}{{\mathcal V}}

\newcommand{\kx}{{\mathcal X}}
\newcommand{\ky}{{\mathcal Y}}

\def\kog{\leavevmode\raise.3ex\hbox{$\scriptscriptstyle\langle\!\langle$}}
\def\fg{\leavevmode\raise.3ex\hbox{$\scriptscriptstyle\,\rangle\!\rangle$}}

\def\hilb{{\bf Hilb}}
\begin{document}
\newtheorem{xprop}{Proposition}[section]
\newtheorem{xlemm}[xprop]{Lemme}
\newtheorem{xtheo}[xprop]{Th\'eor\`eme}
\newtheorem{xcoro}[xprop]{Corollaire}
\newtheorem{quest}{Question}
\newtheorem{defin}{D\'efinition}

\def\refname{R\'ef\'erences}
\def\contentsname{Sommaire}
\def\proofname{D\'emonstration}
\def\abstractname{R\'esum\'e}

\author[{J.M. Dr\'{e}zet}]{Jean--Marc Dr\'{e}zet}
\address{
Institut de Math\'ematiques de Jussieu,
Case 247,
4 place Jussieu,
F-75252 Paris, France}
\email{drezet@math.jussieu.fr}
\urladdr{http://www.math.jussieu.fr/$\sim$drezet}

\title[{\tiny Vari\'et\'es de modules alternatives}]
{Vari\'et\'es de modules alternatives}
\maketitle
\tableofcontents

\section{Introduction}

Soit $X$ une vari\'et\'e alg\'ebrique complexe projective lisse et
irr\'eductible, de dimension $d$, munie d'un fibr\'e tr\`es ample \m{\ko_X(1)}.
On d\'efinit pr\'ecis\'ement au \para 2 ce qu'est une {\em vari\'et\'e de
modules fins d\'efinie globalement}
 de faisceaux coh\'erents sur $X$. C'est la donn\'ee d'une vari\'et\'e
int\`egre $M$, d'un faisceau coh\'erent $\ke$ sur \ \m{M\times X}, plat sur
$M$, tel que pour tous \m{x,y\in M} distincts, les faisceaux \m{\ke_x},
\m{\ke_y} ne soient pas isomorphes, que $\ke$ soit une d\'eformation compl\`ete
de \m{\ke_x}, et que $\ke$ poss\`ede une propri\'et\'e universelle locale
\'evidente. Dans ce cas on peut voir $M$ comme un ensemble de classes
d'isomorphisme de faisceaux muni d'une structure alg\'ebrique naturelle.
Les exemples les plus connus sont certaines vari\'et\'es de modules de
faisceaux stables. Mais toutes ne sont pas dot\'ees d'un {\em faisceau
universel}, et ne sont donc pas des vari\'et\'es de modules fins (cf.
\cite{dr3}, \cite{yo}). Il existe beaucoup d'autres vari\'et\'es de modules
fins que celles qui sont constitu\'ees de faisceaux stables.

Si on s'int\'eresse \`a des faisceaux non {\em simples} il est n\'ecessaire
d'avoir une d\'efinition moins stricte. On introduit la notion de {\em
vari\'et\'e de modules fins d\'efinie localement}, o\`u le faisceau $\ke$ est
remplac\'e par une famille \m{(\ke_i)} de faisceaux, \m{\ke_i} \'etant d\'efini
sur \ \m{U_i\times X}, \m{U_i} \'etant un ouvert de $M$. Pour tout \
\m{x\in U_i\cap U_j} \ on doit avoir \ \m{\ke_{ix}\simeq\ke_{jx}}. Si les
faisceaux \m{\ke_{ix}} ne sont pas simples, on ne peut pas en g\'en\'eral
{\em recoller} les \m{\ke_i} pour obtenir un faisceau global $\ke$.

Il existe ensuite une s\'erie de d\'efinitions de moins en moins strictes
de vari\'et\'es de modules, jusqu'\`a la notion d'{\em espace de modules
grossiers}, o\`u il n'est plus question de faisceau universel.

Le pr\'esent article est consacr\'e \`a l'\'etude des vari\'et\'es de  modules
fins (d\'efinies localement ou non). On peut se fixer le rang $r$ et les
classes de Chern \ \m{c_i\in H^i(X,\Z)} \ des faisceaux \`a consid\'erer.
On peut distinguer trois principaux cas :

\medskip

\noindent (i) Il n'existe pas de faisceau semi-stable de rang $r$ et de
classes de Chern \m{c_i}.

\medskip

\noindent (ii) On s'int\'eresse aux vari\'et\'es de modules fins contenant un
ouvert non vide de faisceaux stables.

\medskip

\noindent (iii) Il existe des faisceaux semi-stables de rang $r$ et de
classes de Chern \m{c_i}, mais on s'int\'eresse aux vari\'et\'es de modules
fins n'en contenant aucun.

\medskip

On s'int\'eresse surtout ici aux deux premiers cas. Le troisi\`eme a d\'ej\`a
\'et\'e explor\'e par S.A. Str\o mme (cf. \cite{st}), mais dans ce cas il
faudrait sans doute introduire des vari\'et\'es de modules non r\'eduites.

On prouve d'abord des r\'esultats g\'en\'eraux sur les vari\'et\'es de modules
fins. On donne ensuite trois types d'exemples.

Dans le premier on s'int\'eresse
aux {\em faisceaux prioritaires} sur le plan projectif. Ce sont des faisceaux
un peu plus g\'en\'eraux que les faisceaux semi-stables. On construit des
vari\'et\'es de modules fins de faisceaux prioritaires de rang et classes de
Chern donn\'es lorsqu'il n'existe aucun faisceau semi-stable ayant les m\^emes
invariants. Ceci permet de trouver un exemple de vari\'et\'e de modules fins
projective d\'efinie localement mais non globalement.

Dans le second on traite les faisceaux simples de rang 1 (sur \m{X=\P_2} \
ou \ \m{X=\P_1\times\P_1}). On donne ici un exemple de vari\'et\'e de
modules fins non projective mais maximale.

Dans le troisi\`eme on construit des vari\'et\'es de modules fins de faisceaux
coh\'erents \`a partir de vari\'et\'es de modules de morphismes de faisceaux.
On en d\'eduit d'abord l'existence de fibr\'es simples sur le plan projectif
qui sont des d\'eformations de fibr\'es stables, mais qui ne
peuvent \^etre inclus dans aucune vari\'et\'e de modules fins (dans ces
exemples la vari\'et\'e de modules des faisceaux stables correspondante est
une vari\'et\'e de modules fins). On donne
ensuite un exemple de vari\'et\'e de modules fins constitu\'ee de faisceaux
simples sans torsion et qui est projective. Cette vari\'et\'e est diff\'erente
de la vari\'et\'e de modules de faisceaux stables correspondante, et poss\`ede
un ouvert non vide commun avec elle.

\vskip 1.2cm

\noindent{\bf Plan des chapitres suivants : }

\medskip

Le chapitre 2 de cet article est consacr\'e \`a la th\'eorie g\'en\'erale des
vari\'et\'es de modules fins de faisceaux coh\'erents sur $X$. On montre que
tous les faisceaux d'une vari\'et\'e de modules fins ont un espace
d'endomorphismes de la m\^eme dimension (proposition 2.2), et que les
vari\'et\'es de modules fins d\'efinies localement de faisceaux simples sont en
fait d\'efinies globalement (proposition 2.8).

On d\'efinit ensuite au \para 2.7 des {\em vari\'et\'es de modules fins
d'extensions}. Soient \m{(M,\ke)}, \m{(N,\kf)} des vari\'et\'es de modules
fins. On suppose que pour tous \m{x\in M}, \m{y\in N}, on a
$$\ext^1(\ke_x,\kf_y)=\hom(\ke_x,\kf_y)=\ext^2(\ke_x,\ke_x)=
\ext^2(\kf_y,\kf_y)= \ \ \ \ \ \ \ \ \ $$
$$\ \ \ \ \ \ \ \ \ \ \ \ext^2(\kf_y,\ke_x)=\ext^2(\ke_x,\kf_y)=\lbrace
0\rbrace.$$
On consid\`ere les extensions
$$0\lra\ke_x\lra E\lra\kf_y\lra 0$$
telles que l'application lin\'eaire induite
$$\hom(\ke_x,\ke_x)\oplus\hom(\kf_y,\kf_y)\lra\ext^1(\kf_y,\ke_x)$$
soit surjective. On montre qu'il existe une vari\'et\'e de modules fins
d\'efinie localement pour de telles extensions, qui est un ouvert de \
\m{M\times N}. C'est un moyen d'obtenir des vari\'et\'es de modules fins de
faisceaux non simples.

\medskip

Le chapitre 3 traite des {\em faisceaux prioritaires} sur le plan projectif.
On dit qu'un faisceau coh\'erent $E$ sur \m{\P_2} est prioritaire s'il est
sans torsion et si on a \ \m{\ext^2(E,E(-1))=\lbrace 0\rbrace}. Ces faisceaux
ont \'et\'e introduits par A. Hirschowitz et Y. Laszlo dans \cite{hi_la}.
On s'int\'eresse aux rangs $r$ et classes de Chern \m{c_1}, \m{c_2} tels
qu'il existe des faisceaux prioritaires, mais pas de faisceau semi-stable de
rang $r$ et de classes de Chern \m{c_1}, \m{c_2}. On donne dans ce cas une
description des {\em faisceaux prioritaires g\'en\'eriques}. Suivant les
valeurs de $r$, \m{c_1} et \m{c_2} deux cas peuvent se \hbox{produire :}

\medskip

\noindent 1 - Il existe un unique faisceau prioritaire g\'en\'erique, de la
forme
$$(E\ot\C^{m})\oplus(F\ot\C^{n})\oplus(G\ot\C^{p}),$$
$E$, $F$ et $G$
\'etant des {\em fibr\'es exceptionnels} sur \m{\P_2}. Il existe dans ce
cas une unique vari\'et\'e de modules fins de faisceaux prioritaires, r\'eduite
\`a un seul point.

\medskip

\noindent 2 - Il existe une vari\'et\'e de modules de faisceaux stables $M$,
un entier $p$ et un fibr\'e exceptionnel $F$ tels que le faisceau prioritaire
g\'en\'erique soit de la forme
$$\ke \ = \ (F\ot\C^{p})\oplus E$$
avec $E$ dans $M$. En utilisant le th\'eor\`eme 2.7 on construit alors une
vari\'et\'e de modules fins d\'efinie localement de faisceaux prioritaires qui
sont des extensions
$$0\lra F\ot\C^{p}\lra\ke\lra E\lra 0$$
avec $E$ dans $M$. Dans un certain nombre de cas cette vari\'et\'e de modules
fins est isomorphe \`a $M$. On montre au \para 3.4.2.4 qu'il existe des cas
o\`u la vari\'et\'e de modules fins n'est pas d\'efinie globalement.

\bigskip

Dans le chapitre 4 on s'int\'eresse aux faisceaux simples de rang 1. On montre
(th\'eor\`eme 4.2) que si \ \m{X=\P_2} \ ou \ \m{X=\P_1\times\P_1},
toute vari\'et\'e de modules fins de faisceaux simples de rang 1 contient un
ouvert dense constitu\'e de faisceaux d'id\'eaux de sous-sch\'emas finis.
On donne ensuite (th\'eor\`eme 4.7) des exemples de vari\'et\'es de modules
fins de faisceaux simples de rang 1 qui sont maximales, mais non projectives.

\medskip

Dans le chapitre 5 on s'int\'eresse aux vari\'et\'es de modules fins de
faisceaux sur \m{\P_2} construites \`a partir de vari\'et\'es de modules
de morphismes (cf. \cite{dr2}, \cite{ka1}, \cite{ka2}, \cite{dr_tr},
\cite{dr6}, \cite{dr7}). On peut faire la m\^eme \'etude sur n'importe quelle
vari\'et\'e dot\'ee d'une bonne th\'eorie des fibr\'es exceptionnels (comme
par exemple les surfaces de Del Pezzo ou les espaces projectifs de dimension
sup\'erieure).

Commen\c cons par le type le plus simple de morphisme. Soient \m{(E,F)} une
paire exceptionnelle sur \m{\P_2}, $m$, $n$ des entiers positifs premiers
entre eux et tels que
$$n.rg(F)-m.rg(E) \ > \ 0.$$
On consid\`ere les morphismes
$$E\ot\C^{m}\lra F\ot\C^{n}.$$
Soit $W$ l'espace vectoriel des tels morphismes. Sur \m{\P(W)} op\`ere le
groupe alg\'ebrique r\'eductif \ \m{SL(m)\times SL(n)}. On a donc une notion
de {\em morphisme semi-stable} et une vari\'et\'e de modules de morphismes
$N$ qui est lisse et projective. Soient $r$ le rang et \m{c_1}, \m{c_2} les
classes de Chern des conoyaux de tels morphismes injectifs. Si \m{m/n} est
suffisamment petit on montre dans \cite{dr2} que les morphismes stables sont
injectifs, et que leurs conoyaux sont les faisceaux stables de rang $r$ et
de classes de Chern \m{c_1}, \m{c_2}. En g\'en\'eral, soit \m{N_0} l'ouvert de
$N$ constitu\'e des morphismes injectifs. On construit une vari\'et\'e de
modules fins isomorphe \`a \m{N_0} constitu\'ee des conoyaux des morphismes
stables injectifs. On en d\'eduit au \para 5.3.2 une vari\'et\'e de modules
fins projective comportant un ouvert constitu\'e de faisceaux stables, mais
diff\'erente de la vari\'et\'e de modules de ces faisceaux stables. Cette
vari\'et\'e de modules fins contient aussi des faisceaux ayant de la torsion.
R\'eciproquement, on montre aussi que si le conoyau d'un morphisme injectif $f$
fait partie d'une vari\'et\'e de modules fins, alors $f$ est stable. On en
d\'eduit au \para 5.3.1 des exemples de fibr\'es simples sur \m{\P_2} qui
sont des d\'eformations de faisceaux stables mais qui ne peuvent pas \^etre
inclus dans une vari\'et\'e de modules fins.

Soit \m{(E,G,F)} une suite exceptionnelle sur \m{\P_2} (aussi appel\'ee
{\em triade}). On peut aussi consid\'erer des morphismes du type
$$E\ot\C^{m}\lra(F\ot\C^{n})\oplus(G\ot\C^{p}).$$
On construit dans \cite{dr_tr}, \cite{dr6} et \cite{dr7} des vari\'et\'es de
modules de tels morphismes. On a ici une action du groupe alg\'ebrique non
r\'eductif \ \m{GL(m)\times\Aut((F\ot\C^{n})\oplus(G\ot\C^{p}))}. On d\'efinit
une notion de stabilit\'e pour les morphismes pr\'ec\'edents \`a partir
d'une paire \m{(\lambda,\mu)} de nombres rationnels positifs tels que \
\m{n\lambda+p\mu=1} (appel\'ee
{\em polarisation}). D'apr\`es \cite{dr5} certaines vari\'et\'es de
modules {\em extr\'emales} de faisceaux semi-stables sur \m{\P_2} sont
isomorphes \`a de telles vari\'et\'es de modules de morphismes. On consid\`ere
au \para 5.3.3 des morphismes
$$\ko(-3)\ot\C^{2}\lra\ko(-2)\oplus(\ko(-1)\ot\C^{7}).$$
Si
$$(\lambda,\mu) \ = \ (\frac{3}{10}+\epsilon,\frac{1}{10}-\frac{\epsilon}{7})$$
(avec $\epsilon$ suffisamment petit) les morphismes stables sont injectifs et
leurs conoyaux sont les faisceaux stables de rang 6 et de classes de Chern -3
et 8. Si
$$\lambda \ > \ \frac{1}{2}$$
les morphismes stables sont encore injectifs mais il peut arriver que le
conoyau ne soit pas stable. On obtient dans ce cas une vari\'et\'e de modules
fins (isomorphe \`a la vari\'et\'e de modules de morphismes) constitu\'ee de
faisceaux sans torsion, qui est projective, et qui contient des faisceaux
stables et des faisceaux instables.

\vskip 1.2cm

\noindent{\bf Questions}

\medskip

\noindent 1 - Toute vari\'et\'e de modules fins peut-elle \^etre contenue
dans une vari\'et\'e de modules fins maximale ?
 Existe-t'il une suite infinie de vari\'et\'es de modules fins
$$M_0\subset M_1\subset M_2\ldots,$$
les inclusions \'etant strictes ?

\medskip

\noindent 2 - Existe t'il un nombre fini de vari\'et\'es de modules fins
maximales contenant un faisceau \hbox{donn\'e ?}

\medskip

\noindent 3 - On trouve au \para 5.3.1 des faisceaux simples qui sont des
d\'eformations de faisceaux stables mais qui ne peuvent pas \^etre inclus
dans une vari\'et\'e de modules fins. Comment caract\'eriser les faisceaux
(simples ou non) qui peuvent faire partie de vari\'et\'es de modules fins ?

\medskip

\noindent 4 - On trouve au \para 5.3.3 des vari\'et\'es de modules fins
projectives (et donc maximales) de faisceaux simples sans torsion contenant \`a
la fois des faisceaux stables et des faisceaux instables. Comment
caract\'eriser les vari\'et\'es de modules fins de faisceaux stables ?

\medskip

\noindent 5 - Parmi les faisceaux pouvant \^etre inclus dans une vari\'et\'e
de modules fins on consid\`ere la relation suivante : \m{E\sim F} s'il existe
un germe de courbe lisse \m{(C,x_0)} et des faisceaux coh\'erents $\ke$, $\kf$
sur \m{X\times C}, plats sur $C$, tels que \ \m{\ke_{x_0}\simeq E},
\m{\kf_{x_0}\simeq F}, et \ \m{\ke_x\simeq\kf_x} \ si \ \m{x\not =x_0}. Quelles
sont les classes d'\'equivalence d'un faisceau coh\'erent donn\'e (pour la
relation d'\'equivalence d\'efinie par la relation pr\'ec\'edente) ?
(cette relation d'\'equivalence est \'etudi\'ee au \para 2.8).

\vskip 1.2cm

\noindent{\it Notations et rappels : }

\medskip


Les {\it points} d'une vari\'et\'e alg\'ebrique
sont ses points ferm\'es.

Soit $P$ un polyn\^ome \`a coefficients rationnels, $E$ un faisceau coh\'erent
sur $X$. On note \m{\hilb^P(E)} le sch\'ema de Hilbert-Grothendieck des
quotients de $E$ de polyn\^ome de Hilbert $P$ (cf \cite{gr}, \cite{lp2}).

Si \m{\ke}, \m{\kf} sont des faisceaux coh\'erents sur une vari\'et\'e
alg\'ebrique, \m{\hom(\ke,\kf)} d\'esigne l'espace vectoriel des morphismes de
$\ke$ dans $\kf$, et \m{\underline{\rm Hom}(\ke,\kf)} le faisceau des
morphismes de $\ke$ dans $\kf$.

On note \ \m{j:X\lra\P_n} \ le plongement d\'efini par \m{\ko_X(1)}, et
$Q$ le fibr\'e quotient canonique sur \m{\P_n}.

Rappelons que le th\'eor\`eme de Riemann-Roch s'\'ecrit pour un faisceau
coh\'erent $E$ de rang \ \m{r>0} \ sur une surface projective lisse
irr\'eductible $X$ s'\'ecrit
$$\chi(E) \ = \ rg(E).(P(\mu(E))-\Delta(E)),$$
\m{\chi(E)} d\'esignant la caract\'erisitique d'Euler-Poincar\'e de $E$, avec
$$P(\mu) \ = \ \frac{\mu(\mu-\omega_X)}{2}+\chi(\ko_X)$$
pour tout \ \m{\mu\in A^1(X)\ot\Q}, \m{\omega_X} d\'esignant le fibr\'e
canonique sur $X$ et
$$\mu(E) \ = \ \frac{c_1(E)}{r}, \ \ \ \ \Delta(E) \ = \ \frac{1}{r}(c_2(E)-
\frac{r-1}{2r}c_1(E)^2).$$
Si $E$, $F$ sont des faisceaux coh\'erents sur $X$, on pose
$$\chi(E,F) \ = \ \sigg_{0\leq i\leq 2}(-1)^i\dim(\ext^i(E,F)).$$
On a, si \ \m{rg(E)>0} \ et \ \m{rg(F)>0},
$$\chi(E,F) \ = \ rg(E).rg(F).(P(\mu(F)-\mu(E))-\Delta(E)-\Delta(F)).$$
Si $E$ et $F$ ne sont plus n\'ecessairement de rang positif, on a
$$\chi(E,F) \ = \ -c_1(E)c_1(F) -rg(E)c_2(F)-rg(F)c_2(E)+rg(E)rg(F)
\chi(\ko_X)$$
$$ \ \ \ \ \ \ \ \ \ + \frac{1}{2}(rg(F).\omega_Xc_1(E)-rg(E).\omega_Xc_1(F)+
rg(E)c_1(F)^2+rg(F)c_1(E)^2).$$
On a en g\'en\'eral, pour tout entier $i$, un isomorphisme canonique

\ \m{\ext^i(E,F) \ \simeq \ \ext^{2-i}(F,E\ot K_X)^*}
\ \ (dualit\'e de Serre, cf. \cite{dr_lp}, prop. (1.2)).

\bigskip

Soient $F$ un faisceau coh\'erent sur une vari\'et\'e alg\'ebrique et
\[ 0=F_0\subset F_1\subset\cdots\subset F_{n-1}\subset F_n=F \]
une filtration
de $F$. Soit \ \m{E_i=F_i/F_{i-1}}. Il existe alors une suite spectrale
bien connue \m{E_r^{pq}} d'aboutissement \m{\ext^\cdot(F,F)} et telle que
\[ E_1^{pq} \ = \ \som_i\ext^{p+q}(E_i,E_{i-p}). \]
La construction de cette suite spectrale est similaire \`a celle qui est
effectu\'ee dans le \para 1.5 de \cite{dr_lp}.

\vskip 2.5cm

\section{Vari\'et\'es de modules fins}

\begin{sub}{\bf Familles de faisceaux coh\'erents}\end{sub}

Soit $S$ une vari\'et\'e alg\'ebrique.

\medskip

\begin{defin}
On appelle {\em famille de faisceaux sur $X$ param\'etr\'ee par $S$} un
faisceau coh\'erent $\kf$ sur \ \m{S\times X}, plat sur $S$.
\end{defin}

\medskip

\begin{defin}
On appelle {\em polyfamille de faisceaux sur $X$ param\'etr\'ee par $S$} la
donn\'ee d'un recouvrement ouvert \m{(U_i)_{i\in I}} de $S$, et pour tout
\m{i\in I} d'une famille \m{\ke_i} de faisceaux sur $X$ param\'etr\'ee par
\m{U_i}, tels que pour tous \m{i,j\in I} et tout \m{x\in U_i\cap U_j} on ait
\ \m{\ke_{i,x}\simeq\ke_{j,x}}.
\end{defin}

\medskip

\noindent{\em Notations : } Pour toute sous-vari\'et\'e localament ferm\'ee
\m{S'} de $S$, on note
$$\kf_{S'} \ = \ \kf_{\mid S'\times X}.$$
Si \m{S'} est r\'eduite \`a un point $s$ on notera plus simplement \
\m{\kf_{S'} = \kf_s}.
On note \m{p_S} (resp. \m{p_X}) la projection \ \m{S\times X\lra S} \
(resp.\ \m{S\times X\lra X}). Si $\ke$, $\kf$, $\kg$ sont des
faisceaux coh\'erents sur \ \m{S\times X}, $S$ et $X$
respectivement, on notera plus simplement
$$\ke\ot p_S^*(\kf)\ot p_X^*(\kg) \ = \ \ke\ot\kf\ot\kg.$$
Si \ \m{f : T\lra S} \ est un morphisme de vari\'et\'es alg\'ebriques, et
$\kf$ une famille de faisceaux sur $X$ param\'etr\'ee par $S$, on note
$$f^\sharp(\kf) \ = \ (f\times I_{X})^*(\kf).$$
C'est une famille de faisceaux sur $X$ param\'etr\'ee par $T$.

Si $\kx$ est un ensemble non vide de classes d'isomorphisme de faisceaux
coh\'erents sur $X$, on appelle {\em famille de faisceaux de $\kx$
param\'etr\'ee par $S$} une famille $\kf$ de faisceaux coh\'erents sur $X$
param\'etr\'ee par $S$ telle que pour tout point $s$ de $S$ la classe
d'isomorphisme de \m{\kf_s} soit dans $\kx$.

\vskip 1.2cm

\begin{sub}{\bf Ensembles ouverts de faisceaux coh\'erents sur $X$}
\end{sub}

Soient $r$, \m{c_i\in H^i(X,\Z)}, \m{1\leq i\leq d}, avec \ \m{r\geq 0}.
Soit $\kx$ un ensemble non vide de classes d'isomorphisme de faisceaux
coh\'erents sur $X$, de rang $r$ et de classes de Chern \m{c_i}.

\medskip

\begin{defin}
On dit que $\kx$ est un {\em ensemble ouvert} si pour toute vari\'et\'e
alg\'ebrique $S$ et toute famille $\kf$ de faisceaux coh\'erents sur $X$
de rang $r$ et de classes de Chern \m{c_i} param\'etr\'ee par $S$,
l'ensemble des points $s$ de $S$ tels que la classe d'isomorphisme
de \m{\kf_s} soit dans $\kx$ est un ouvert de Zariski de $S$.
\end{defin}

\medskip

\begin{defin}
Supposons que $\kx$ soit un ensemble ouvert. On dit que $\kx$ est
{\em irr\'eductible} si pour toutes polyfamilles de faisceaux de $\kx$
param\'etr\'ees par des vari\'et\'es alg\'ebriques \m{X_1}, \m{X_2}
respectivement, il existe un polyfamille de faisceaux de $\kx$ param\'etr\'ee
par une vari\'et\'e alg\'ebrique irr\'eductible $Y$ contenant tous les
faisceaux des polyfamilles param\'etr\'ees par \m{X_1} et \m{X_2}.
\end{defin}

\vskip 1.2cm

\begin{sub}{\bf Vari\'et\'es de modules fins - D\'efinition}\end{sub}

Soient $r$, \m{c_i\in H^i(X,\Z)}, \m{1\leq i\leq d}, avec \ \m{r\geq 0}. Soit
$\kx$ un ensemble ouvert de faisceaux coh\'erents sur $X$ de rang $r$ et de
classes de Chern \m{c_i}. On pose
$$\Delta \ = \ \frac{1}{r}(c_2-\frac{r-1}{2r}c_1^2).$$

\medskip

\begin{defin}
On appelle {\em vari\'et\'e de modules fins globale}, ou plus simplement
{\em vari\'et\'e de modules fins}
pour $\kx$ la donn\'ee d'une vari\'et\'e alg\'ebrique int\`egre $M$ et
d'une famille $\ke$ de faisceaux de $\kx$ param\'etr\'ee par $M$ telles que :

\noindent (i) Pour tout \'el\'ement $x$ de $\kx$, il existe un unique point
$m$ de $M$ tel que la classe d'isomorphisme de \m{\ke_m} soit $x$.

\noindent (ii) Pour toute famille $\kf$ de faisceaux de $\kx$ param\'etr\'ee
par une vari\'et\'e alg\'ebrique $S$, il existe un morphisme \
\m{f:S\lra M} \ tel que pour tout point $s$ de $S$ il existe un ouvert $U$ de
$S$ contenant $s$ et un isomorphisme \ \m{\kf_U\simeq f^\sharp(\ke)_U}.

\medskip

{\em
On dit aussi dans ce cas que \m{(M,\ke)} est une vari\'et\'e de modules fins
pour $\kx$. On remarquera que d'apr\`es (i), il existe une bijection canonique
entre $\kx$ et l'ensemble des points de $M$. D'autre part le morphisme
$f$ de (ii) est unique : l'image d'un point $s$ de $S$ est le point de $M$
correspondant \`a \m{\kf_s}.}

\medskip

On appelle {\em vari\'et\'e de modules fins d\'efinie localement} pour $\kx$
la donn\'ee d'une vari\'et\'e alg\'ebrique int\`egre $M$, d'un recouvrement
ouvert \m{(U_i)_{i\in I}} de $M$, et pour tout \m{i\in I} d'une famille
\m{\ke_i} de faisceaux de $\kx$ param\'etr\'ee par \m{U_i} tels que

\noindent (i) Pour tout \'el\'ement $x$ de $\kx$, il existe un unique point
$m$ de $M$ tel que pour tout \m{i\in I} tel que \m{x\in U_i}, la classe
d'isomorphisme de \m{\ke_{im}} soit $x$.

\noindent (ii) Pour toute famille $\kf$ de faisceaux de $\kx$ param\'etr\'ee
par une vari\'et\'e alg\'ebrique $S$, il existe un morphisme \
\m{f:S\lra M} \ tel que pour tout point $s$ de $S$ et tout \m{i\in I} tel
qu'il existe un \m{x\in U_i} tel que \m{\ke_{ix}\simeq\kf_s} il existe un
ouvert $U$ de $S$ contenant $s$ tel que \ \m{f(U)\subset U_i} \ et un
isomorphisme \ \m{\kf_U\simeq f^\sharp(\ke_i)_U}.
\end{defin}

\medskip

\noindent La plupart des r\'esultats qui suivent dont la d\'emonstration est
donn\'ee pour les vari\'et\'es de modules fins sont aussi valables pour les
vari\'et\'es de modules fins d\'efinies localement.

\medskip

\noindent{\bf Remarque : } Dans ce qui pr\'ec\`ede, \m{(U_i),(\ke_i)} est donc
une {\em polyfamille} de faisceaux de $\kx$ param\'etr\'ee par $M$ (cf.
d\'efinition 4).

\medskip

\begin{xprop}
Si \m{(M,\ke)} et \m{(M',\ke')} sont deux vari\'et\'es de modules fins pour
$\kx$, il existe un isomorphisme \ \m{\phi : M\lra M'} \ tel que pour tout
point $m$ de $M$ il existe un ouvert $U$ de $M$ contenant $m$ et un
isomorphisme \ \m{\ke_U\simeq\phi^\sharp(\ke'_{\phi(U)})}.

R\'eciproquement, si \m{(M,\ke)} est une vari\'et\'e de modules fin pour $\kx$
et si \m{\ke'} est une famille de faisceaux de $\kx$ param\'etr\'ee par une
vari\'et\'e alg\'ebrique r\'eduite \m{M'} poss\'edant la propri\'et\'e
pr\'ec\'edente, alors \m{(M',\ke')} est une vari\'et\'e de modules fin pour
$\kx$.
\end{xprop}

\begin{proof} Imm\'ediat.

\medskip

On a un r\'esultat analogue pour les vari\'et\'es de modules d\'efinies
localement.

Une vari\'et\'e de modules fins induit naturellement une vari\'et\'e de modules
d\'efinie localement (avec un recouvrement ouvert r\'eduit \`a un seul
\'el\'ement). On donne au \para 3.4.2.4 un exemple de vari\'et\'e de modules
fins d\'efinie localement qui ne peut pas l'\^etre globalement.

Il d\'ecoule de la th\'eorie des d\'eformations des faisceaux (cf.
\cite{si_tr}) que si \m{(M,\ke)} est une vari\'et\'e de modules fins pour $\kx$
et si \m{x\in M}, alors l'espace tangent \m{T_x(M)} s'identifie canoniquement
\`a \m{\ext^1(\ke_x,\ke_x)}. Il en d\'ecoule en particulier que si cette
dimension est ind\'ependante de $x$, alors $M$ est lisse. On a un r\'esultat
analogue pour les vari\'et\'es de modules d\'efinies localement.

\vskip 1.2cm

\begin{sub}{\bf Endomorphismes des faisceaux d'une vari\'et\'e de modules fins}
\end{sub}

\begin{xprop}
Supposons qu'il existe une vari\'et\'e de modules fins \m{(M,\ke)} pour $\kx$.
Alors pour tout \m{x\in M}, \m{\Aut(\ke_x)} agit trivialement sur
\m{\ext^1(\ke_x,\ke_x)}, et \m{\dim(\hom(\ke_x,\ke_x))} est ind\'ependant de 
\m{x}.
\end{xprop}

\begin{proof} Soit $\kf$ une famille de faisceaux de $\kx$ param\'etr\'ee par
une vari\'et\'e alg\'ebrique $Y$. Il existe un entier $m$ tel que pour
tout \m{m'\geq m-n} et tout \m{y\in Y} on ait \m{h^i(\kf_y(m'))=0} \ pour
\m{i>0}. On a alors
un complexe
$$(\Lambda^n Q^*)(-m)\ot G_n=R^{-n}\lra\ldots\lra
(\Lambda^iQ^*)(-m)\ot G_i=R^{-i}\lra\ldots\ \ \ \ \ \ \ \ $$
$$ \ \ \ \ \ \ \ \ \ \ \ \ \ \ \ \ \ \ \ldots\lra \ko(-m)\ot G_0=R^0$$
donn\'e par la suite spectrale de Beilinson sur \m{\P_n}, avec
$$G_i \ = \ H^0(\kf_y(m-i)).$$
Ce complexe d\'epend fonctoriellement de $\kf_y$, est exact en degr\'e $<0$ et
sa cohomologie en degr\'e 0 est isomorphe \`a $\kf_y$. On montre ais\'ement que
le morphisme canonique
$$\End(R^.)\lra\End(\kf_y)$$
est un isomorphisme (cela peut se faire directement ou \`a l'aide d'une suite
spectrale comme dans le lemme 23 de \cite{lp0}). Compte tenu du fait que les
fibr\'es \m{\Lambda^iQ^*} sont simples on en d\'eduit que
\m{\End(\kf_y)} est isomorphe au noyau de l'application canonique
\[
\som_{0\leq i\leq n}L(G_i,G_i)\lra\som_{1\leq i\leq n}L(G_i,G_{i-1})\ot
\hom(\Lambda^iQ^*,\Lambda^{i-1}Q^*).
\]
La r\'esolution pr\'ec\'edente de $\kf_y$ se globalise en une r\'esolution
localement libre sur \ \m{Y\times\P_n} :
$$0\lra (\Lambda^n Q^*)(-m)\ot \kg_n\lra\ldots\lra (\Lambda^iQ^*)(-m)\ot
\kg_i\lra\ldots \ \ \ \ \ \ \ \ \ $$
$$\ \ \ \ \ \ \ \ \ \ \ \lra \ko(-m)\ot \kg_0\lra j_*(\kf)\lra 0$$
(o\`u $j$ d\'esigne ici l'inclusion \ \m{Y\times X\subset Y\times\P_n}).
On obtient alors un morphisme canonique de fibr\'es vectoriels sur $Y$
$$\Psi : \som_{0\leq i\leq n}(\kg_i^*\ot \kg_i)\lra
\som_{1\leq i\leq n}(\kg_i^*\ot \kg_{i-1}\ot\hom(\Lambda^iQ^*,
\Lambda^{i-1}Q^*))$$
tel que pour tout \m{y\in Y}, on ait
$$\hom(\kf_y,\kf_y) \ \simeq \ \ker(\Psi_y).$$
On suppose maintenant que \ \m{Y=M} \ et \ \m{\kf=\ke}.
Soit $k$ le rang de \m{\ker(\Psi)}. On a donc pour tout \m{x\in M},
\ \m{\dim(\hom(\ke_x,\ke_x))\geq k}, l'\'egalit\'e \'etant v\'erifi\'ee sur un
ouvert dense de $M$.

\medskip

\begin{xlemm}
Soient \ \m{\psi : E\lra F} \ un morphisme de fibr\'es vectoriels sur une
vari\'et\'e alg\'ebrique irr\'eductible $Y$, $p$ le rang de \m{\ker(\psi)}.
Soient \m{y\in Y}, et \m{E'_y\subset E_y} l'image dans \m{E_y} de
\m{\ker(\psi)}. Alors on a \ \m{\dim(E'_y)\leq p}.
\end{xlemm}

\begin{proof} On a une suite exacte
$$0\lra E'_y\lra E_y\lra \imm(\psi)_y\lra 0,$$
qui prouve que \m{\dim(E'_y)} est plus petit aux points o\`u \m{\psi_y} n'est
pas de rang maximal. \end{proof}

\medskip

\begin{xlemm}
Soit $\kf$ une famille de faisceaux de $\kx$ param\'etr\'ee par une vari\'et\'e
alg\'ebrique irr\'eductible $Y$ telle que l'image du morphisme induit
\ \m{f : Y\lra M} \
contienne un ouvert de $M$. Soit $k$ la dimension g\'en\'erique de
\m{\hom(\kf_y,\kf_y)} (c'est-\`a-dire que la dimension de cet espace vectoriel
est $k$ pour tout $y$ dans un ouvert non vide de $M$).
Alors pour tout \m{y\in Y} il existe un sous-espace
vectoriel $W$ de \m{\hom(\kf_y,\kf_y)} tel que \ \m{\dim(W)\leq k}, et que
pour tout \ \m{\sigma\in\hom(\kf,\kf)}, on ait \ \m{\sigma_y\in W}.
\end{xlemm}

\begin{proof} D\'ecoule imm\'ediatement du lemme pr\'ec\'edent. \end{proof}

\medskip

\noindent{\it Fin de la d\'emonstation de la proposition 2.2.}
Soit $m$ un entier tel
que pour tout point $x$ de $M$, $\ke_x(m)$ soit engendr\'e par ses sections
globales, et que \ \m{h^i(\ke_x(m)) = \lbrace 0\rbrace} \ pour \m{i\geq 1}.
Soient \ \m{p=h^0(\ke_x(m))}, qui est ind\'ependant de $x$, et $P$ le
polyn\^ome de Hilbert des faisceaux de $\kx$. Soient enfin \
$${\bf Q} \ = \ \hilb^P(\ko_X(-m)\ot\C^{p}),$$
et
$$\pi : \ko_X(-m)\ot\C^{p}\lra\kf$$
le morphisme surjectif universel sur \ \m{{\bf Q}\ot X}. Sur $\bf Q$ on
a une action canonique de \m{GL(p)}. Soit $x$
un point de $M$. Soit $q$ un point de $\bf Q$ tel que \ \m{\kf_q\simeq\ke_x},
et que \m{\pi_q} induise un isomorphisme \ \m{\C^{p}\simeq H^0(\kf_q(m))}. Il
existe par d\'efinition de \m{(M,\ke)} un ouvert $U$ de $\bf Q$ contenant $q$
et un morphisme \ \m{f : U\lra M} \ tel que
$$f^\sharp(\ke) \ \simeq \ \kf_U.$$
Soit $W$ le sous-espace vectoriel de dimension $k$ de \m{\hom(\ke_x,\ke_x)}
d\'efini dans le lemme pr\'ec\'edent, correspondant \`a la famille
\m{f^\sharp(\ke)}.
Rappelons que le groupe \m{\Aut(\kf_q)} s'identifie canoniquement au
stabilisateur de $q$ dans \m{GL(p)}. Soit \ \m{g\in\Aut(\kf_q)}. On a un
isomorphisme canonique \ \m{g^*(\kf)\simeq\kf}. Au point $q$, cet isomorphisme
n'est autre que $g$. Sur \ \m{Z = g^{-1}U\cap U}, on a
$$g^*(f^\sharp(\ke)) \ = \ f^\sharp(\ke).$$
On obtient donc un automorphisme de \m{f^\sharp(\ke)_Z} :
$$\sigma : f^\sharp(\ke)_Z \ \simeq \ \kf_Z \ \simeq \ g^*(\kf)_Z \ \simeq \
g^*(f^\sharp(\ke)).$$
On en d\'eduit que \ \m{g\in W}. Donc \ \m{\dim(\hom(\ke_x,\ke_x))=k}. Pour
d\'emontrer la premi\`ere assertion de la proposition 2.2, on remarque que
le morphisme de d\'eformation infinit\'esimale de Koda\"\i ra-Spencer de
\m{g^*(f^\sharp(\ke))} au point $q$ n'est autre que la compos\'ee de celui
de \m{f^\sharp(\ke)} et de la multiplication par $g$ dans
\m{\ext^1(\ke_x,\ke_x)}. Mais les deux morphismes de d\'eformations sont
\'egaux \`a l'application tangente de $f$ en $q$.
\end{proof}

\medskip

On d\'eduit de la proposition 2.2 la

\medskip

\begin{xprop}
Supposons que $X$ soit une surface K3 et que $\kx$ admette une vari\'et\'e de
modules fins \m{(M,\ke)}. Alors $M$ est lisse et de dimension \
\m{k+r^2(2\Delta-\chi(\ko_X))}.
\end{xprop}

\medskip

Les propositions 2.2 et 2.5 sont aussi vraies pour les vari\'et\'es de modules
fins d\'efinies localement.

Pour qu'il existe une vari\'et\'e de modules fins pour
$\kx$ il est n\'ecessaire d'apr\`es la proposition 2.2 que \m{\dim(\hom(E,E))}
soit le m\^eme pour tous les faisceaux $E$ de $\kx$.
La proposition 2.2 admet une sorte de r\'eciproque :

\bigskip

\begin{xprop}
On suppose que pour tout faisceau $E$ de $\kx$ on a \
\m{\ext^2(E,E)=\lbrace 0\rbrace}, et que \ \m{\dim(\hom(E,E))=p} \ est
ind\'ependant de $E$. Soit $\ke$ une famille de faisceaux de $\kx$
param\'etr\'ee par une vari\'et\'e lisse $M$ telle que pour tout faisceau $E$
de $\kx$ il existe un unique point $x$ de $M$ tel que \ \m{\ke_x\simeq E}, et
que $\ke$ soit un d\'eformation compl\`ete de \m{\ke_x}. Alors \m{(M,\ke)} est
une vari\'et\'e de modules fins pour $\kx$.
\end{xprop}

\medskip

\begin{xlemm}
Soient $M$ une vari\'et\'e alg\'ebrique int\`egre et $\ke$ une famille de
faisceaux de $\kx$ param\'etr\'ee par $M$. On suppose que tout faisceau de
$\kx$ est isomorphe \`a un unique \m{\ke_x}, \m{x\in M}, que
\ \m{\ext^2(\ke_x,\ke_x)=\lbrace 0\rbrace}, que \ \m{\dim(\hom(\ke_x,\ke_x))=p}
\ est ind\'ependant de $x$ et que pour toute famille
$\kf$ de faisceaux de $\kx$ param\'etr\'ee par une vari\'et\'e alg\'ebrique
$Y$, il existe un morphisme \ \m{f:Y\lra M} \ associant \`a \m{y\in Y}
l'unique point
$x$ de $M$ tel que \ \m{\ke_x\simeq\kf_y}. Alors \m{(M,\ke)}
est une vari\'et\'e de modules fins pour $\kx$.
\end{xlemm}

\begin{proof} Il faut montrer que $\ke$ est {\em universelle} localement. On se
ram\`ene ais\'ement au cas o\`u $Y$ est un ouvert du sch\'ema de Grothendieck
$Q$ de la d\'emonstration de la proposition 2.2 et $\kf$ la restriction du
faisceau universel sur $Q$. Notons que $Y$ est lisse. Soient
\m{y\in Y} et
$$\phi:\kf_y\lra\ke_{f(x)}$$
un isomorphisme. Le faisceau
$$G=p_{Y*}(\underline{\hom}(\kf,f^\sharp(\ke)))$$
est localement libre de rang $p$. Soit $G_0$ l'ouvert de $G$ constitu\'e des
isomorphismes. On peut voir $\phi$ comme un point de \m{G_0}. Puisque \m{G_0}
est ouvert il existe une section locale de \m{G_0} d\'efinie en $y$ dont
la valeur en $y$ est $\psi$. Ceci donne un isomorphisme local \
\m{\kf\simeq f^\sharp(E)}.
\end{proof}

\medskip

\noindent{\em D\'emonstration de la proposition 2.6.}
Rappelons qu'on dit que $\ke$ est une {\em d\'eformation compl\`ete} de
\m{\ke_x} si le morphisme de d\'eformation infinit\'esimale de
Koda\"\i ra-Spencer
$$\omega_x : TM_x\lra\ext^1(\ke_x,\ke_x)$$
est surjectif. Soit \m{((S,s_0),\kg)} une d\'eformation verselle de \m{\ke_x}
(cf. \cite{si_tr}), \m{(S,s_0)} \'etant donc un germe de vari\'et\'e
analytique lisse et \m{\kg} un faisceau analytique sur \m{S\times X} tel que \
\m{\kg_{s_0}\simeq\ke_x} (la lissit\'e de $S$ d\'ecoule du fait que \
\m{\ext^2(\ke_x,\ke_x)=\lbrace 0\rbrace}). Soit \m{M_0} le germe de vari\'et\'e
analytique d\'efini par $M$ au point $x$. De $\ke$ on d\'eduit un morphisme \
\m{\phi:M_0\lra S} \ tel que \ \m{\phi^*(\kg)\simeq\ke_{M_0}} dont
l'application tangente est \m{\omega_x}. Les fibres de $\phi$ sont de dimension
0 car sinon tous les faisceaux de la famille $\ke$ ne seraient pas distincts.
Il en d\'ecoule que \m{\omega_x} est en fait un isomorphisme, ainsi que $\phi$.
Donc \m{\ke_x} admet une d\'eformation verselle {\em alg\'ebrique}.

Soit \m{\ke'} une famille de faisceaux de $\kx$ param\'etr\'ee par une
vari\'et\'e alg\'ebrique $Y$. Il suffit d'apr\`es le lemme 2.7 de montrer
que l'application \ \m{f:Y\lra M} \ associant \`a $y$ l'unique point $x$ tel
que \ \m{\ke'_y\simeq\ke_x} \ est r\'eguli\`ere. C'est clair car c'est le
recollement de tous les morphismes universels de $Y$ vers les
d\'eformations verselles des faisceaux \m{\ke'_y}. \end{proof}

\medskip

La proposition 2.6 se g\'en\'eralise aussi sans peine aux
vari\'et\'es de modules fins d\'efinies localement.

\vskip 1.2cm

\begin{sub}{\bf Faisceaux simples}\end{sub}

Soit $\kx$ un ensemble ouvert de faisceaux de rang $r$ et de classes de
Chern \m{c_i} sur $X$, admettant une vari\'et\'e de modules
fins \m{(M,\ke)}. Il d\'ecoule de la proposition 2.2 que s'il existe dans $\kx$
un faisceau {\em simple} (c'est-\`a-dire dont tous les seuls endomorphismes
soient les homoth\'eties), alors tous les faisceaux de $\kx$ sont simples.
C'est le cas en particulier lorsque $\kx$ contient un faisceau stable
relativement \`a \m{\ko_X(1)}.

On suppose maintenant que les faisceaux de $\kx$ sont simples. Il n'y a
pas lieu alors de distinguer les vari\'et\'es de modules fins et les
vari\'et\'es de modules fins d\'efinies localement :

\medskip

\begin{xprop}
Si tous les faisceaux de $\kx$ sont simples, et s'il existe une vari\'et\'e de
modules fins d\'efinie localement pour $\kx$, il existe une vari\'et\'e de
modules fins pour $\kx$.
\end{xprop}

\begin{proof} Supposons qu'il existe une vari\'et\'e de modules fins d\'efinie
localement pour $\kx$, donn\'ee par la vari\'et\'e $M$, le recouvrement ouvert
\m{(U_i)_{i\in I}} de $M$, et les familles \m{\ke_i} de faisceaux de $\kx$.
Il faut construire un {\em faisceau universel} $\ke$ sur \ \m{M\times X}.
Supposons d'abord que les faisceaux de $\kx$ soient localement libres.
Puisque les fibr\'es de $\kx$ sont simples, on peut construire un fibr\'e en
espaces projectifs $\bf P$ sur \ \m{M\times X}\ en recollant les
\m{\P(\ke_i)}. On a pour tout \m{i\in I}, \ \m{{\bf P}_{\mid U_i}\simeq
\P(\ke_i)}. Il en d\'ecoule d'apr\`es \cite{hi_na} qu'il existe un fibr\'e
vectoriel $\ke$ sur \ \m{M\times X} \ tel que \ \m{{\bf P}\simeq\P(\ke)}.
Il est alors facile de voir que \m{(M,\ke)} est une vari\'et\'e de modules
fins pour $\kx$.

Si tous les faisceaux de $\kx$ ne sont pas localement libres, on consid\`ere
le plongement \ \m{j : X\subset\P_n} \
d\'efini par \m{\ko_X(1)}, et on voit les
faisceaux coh\'erents sur $X$ comme des faisceaux coh\'erents sur \m{\P_n}.
On note $Q$ le fibr\'e quotient universel sur \m{\P_n}.
Les faisceaux de $\kx$ constituent en {\em ensemble limit\'e} de faisceaux :
en effet, un nombre fini de \m{U_i} suffit \`a recouvrir $M$. Il existe donc
un entier $m$ tel que pour tout \ \m{m'\geq m} et tout faisceau $E$ de $\kx$
on ait \ \m{H^i(j_*(E)(m'))=\lbrace 0\rbrace} \ pour \ \m{i\geq 1}. On a alors
une r\'esolution {\em canonique} de \m{j_*(E)} :


$$0\lra\ko(-m-1)\ot V_n(E)\ \hfl{\alpha_n(E)}{} \ \ldots\lra(\Lambda^iQ^*)(-m)
\ot V_i(E)\ \hfl{\alpha_i(E)}{} \ \ldots
\ \ \ \ \ \ \ \ \ \ \ $$
$$\ \ \ \ \ \ \ \ \ \ \ \ \ \ \ \ \ \ \lra\ko(-m)\ot V_0(E)\lra j_*(E)\lra 0,$$
avec \ \m{V_i(E)=H^0(E(m-i))} \ pour \m{0\leq i\leq n}. Comme dans le cas
o\`u les faisceaux de $\kx$ sont localement libres, les $V_i(E)$, $E$
parcourant $\kx$, se recollent pour d\'efinir un fibr\'e vectoriel \m{V_n} sur
$M$. Les morphismes \m{\alpha_i(E)} se recollent aussi pour d\'efinir
$$\alpha_i : (\Lambda^iQ^*)\ot V_i\lra (\Lambda^{i-1}Q^*)\ot V_{i-1}.$$
Le faisceau \m{\coker(\alpha_1)} est le faisceau universel $\ke$ sur \
\m{M\times X} \ recherch\'e. \end{proof}

\medskip

On conserve les notations du \para 2.4. Soit \m{{\bf Q}_0} l'ouvert de
\m{\bf Q} constitu\'e des points $q$ tels que la classe d'isomorphisme de
\m{\kf_q} soit dans $\kx$ et \m{\pi_q} induise un isomorphisme \
\m{\C^{p}\simeq H^0(\kf_q(m))}. Alors on peut reconstruire \m{{\bf Q}_0},
\m{\kf_{Q_0}} et \m{\pi} sur \m{{\bf Q}_0} \`a partir de \m{(M,\ke)}. On
consid\`ere pour cela le fibr\'e vectoriel
$$E \ = \ p_{M*}({\rm \underline{Hom}}(\ko_X(-m)\ot\C^{p},\ke))$$
sur \m{M} et l'ouvert \m{U_0} de \m{\P(E)} correspondant aux
isomorphismes. Soit \ \m{\phi : U_0\lra M} \ la projection. On a
un morphisme canonique
$$\pi_0 : \ko_X(-m)\ot\C^{p}\lra\phi^*(\ke)$$
de faisceaux sur \ \m{U_0\times X}. Alors on a un isomorphisme canonique
\ \m{U_0\simeq {\bf Q_0}} \ et \m{\pi_0} s'identifie naturellement \`a la
restriction de $\pi$. On en d\`eduit imm\'ediatement la

\medskip

\begin{xprop}
Le morphisme \ \m{\phi : {\bf Q_0}\lra M} \ est un quotient g\'eom\'etrique de
\m{\bf Q_0} par \m{PGL(p)}.
\end{xprop}

\medskip

On suppose que $X$ est {\em une surface de Del Pezzo} (c'est-\`a-dire que
\m{-K_X} est tr\`es ample). Alors si $E$ est un faisceau coh\'erent simple
sur $X$, on a \ \m{\ext^2(E,E)=\lbrace 0\rbrace}. En effet on a par dualit\'e
de Serre \ \m{\ext^2(E,E)\simeq\hom(E,E\ot K_X)^*}, et si ce dernier n'est
pas nul on peut en composant avec des morphismes \ \m{K_X\lra\ko_X} \ trouver
des endomorphismes de $E$ qui ne sont pas des homoth\'eties. On en d\'eduit la

\medskip

\begin{xprop}
Supposons que $X$ soit une surface de Del Pezzo, que les faisceaux de $\kx$
soient simples et que $\kx$ admette une vari\'et\'e de modules fins
\m{(M,\ke)}. Alors $M$ est une vari\'et\'e lisse de dimension \
\m{1+r^2(2\Delta-1)}.
\end{xprop}

\newpage

\begin{sub}{\bf Vari\'et\'es de modules fins maximales}\end{sub}

Soit \m{(M,\ke)} une vari\'et\'e de modules fins d'un ensemble ouvert $\kx$ de
faisceaux coh\'erents de rang et classes de Chern donn\'es sur $X$.
On dit que $M$ (ou \m{(M,\ke)}, $\kx$) est {\em maximale} s'il
n'existe pas d'ensemble ouvert \m{\kx'} de faisceaux coh\'erents de m\^emes
rang et classes de Chern contenant strictement $\kx$ et admettant une
vari\'et\'e de modules fins \m{(M',\ke')}. J'ignore
si toute vari\'et\'e de modules fins est contenue dans une vari\'et\'e de
modules fins maximale. On a une notion analogue de vari\'et\'e de modules fins
d\'efinie localement maximale.

\vskip 1.2cm

\begin{sub}{\bf Vari\'et\'es de modules fins d'extensions}\end{sub}

\begin{subsub}{Extensions}\end{subsub}

Soient $\kx$ et \m{\kx'} des ensembles ouverts de faisceaux coh\'erents sur
$X$. On suppose qu'il existe des vari\'et\'es de modules fins \m{(M,\ke)} et
\m{(M',\ke')} pour $\kx$ et \m{\kx'} respectivement. Soient $E$ un faisceau de
$\kx$, \m{E'} un faisceau de \m{\kx'}. On consid\`ere une extension
$$(1) \ \ \ \ \ \ \ 0\lra E'\lra G\lra E\lra 0$$
donn\'ee par un \'el\'ement $\epsilon$ de \m{\ext^1(E,E')}. Soit
$$\phi : \hom(E,E)\oplus\hom(E',E')\lra\ext^1(E,E')$$
l'application lin\'eaire d\'eduite de $\epsilon$. On suppose que
$$\ext^1(E',E)=\hom(E',E)=\ext^2(E',E')=\ext^2(E,E)=\ \ \ \ \ \ \ \ \ $$
$$\ \ \ \ \ \ \ \ \ \ \ \ext^2(E,E')=\ext^2(E',E)=\lbrace 0\rbrace.$$

\medskip

\begin{xlemm}
Si $\phi$ est surjective, on a
$$\ext^2(G,G) \ = \ \lbrace 0\rbrace ,$$
un isomorphisme canonique 
$$\ext^1(G,G) \ \simeq \ \ext^1(E,E)\oplus\ext^1(E',E')$$
et une suite exacte canonique
$$0\lra\hom(E,E')\lra\hom(G,G)\lra\ker(\phi)\lra 0.$$
\end{xlemm}

\begin{proof} On consid\`ere la suite spectrale \m{E_r^{p,q}} convergeant vers
\m{\ext^{p+q}(G,G)} d\'efinie par la filtration \ \m{E'\subset G}. Les termes
\m{E_1^{p,q}} \'eventuellement non nuls qui nous int\'eressent pour calculer
\m{\hom(G,G)}, \m{\ext^1(G,G)} et \m{\ext^2(G,G)} sont repr\'esent\'es
ci-dessous :
\[
\begin{array}{ccccccc}\ext^1(E',E)=\lbrace 0\rbrace &  & . &  & . &  & . \\
\hom(E',E)=\lbrace 0\rbrace & & \ext^1(E,E)\oplus\ext^1(E',E') & &
\ext^2(E,E')= \lbrace 0\rbrace & & .\\
0 &  & \hom(E,E)\oplus\hom(E',E') & \hfl{\phi}{} & \ext^1(E,E') & & 0\\
0 &  &  0 & & \hom(E,E') &  & 0 \\ \end{array}
\]
Le lemme en d\'ecoule imm\'ediatement. \end{proof}

\medskip

{\bf Remarque : }
On peut bien s\^ur d\'emontrer le r\'esultat pr\'ec\'edent
sans utiliser de suite spectrale, mais moins rapidement, \`a l'aide de suites
exactes longues et de diagrammes commutatifs.

\vskip 0.8cm

On suppose que
$$\hom(E',E)=\lbrace 0\rbrace$$
pour tous $E$ dans $\kx$ et \m{E'} dans \m{\kx'}.
Soit \m{\ky'} le sous-ensemble de \ \m{\kx\times\kx'} \ constitu\'e des paires
\m{(E,E')} telles que
$$\ext^1(E',E)=\ext^2(E',E')=\ext^2(E,E)=\ext^2(E',E)=\ext^i(E,E')
=\lbrace 0\rbrace$$
et qu'il existe un \'el\'ement $\sigma$ de \m{\ext^1(E,E')} tel que
l'application associ\'ee
$$\phi : \hom(E,E)\oplus\hom(E',E')\lra\ext^1(E,E')$$
soit surjective. L'orbite de cet \'el\'ement sous l'action de \
\m{\Aut(E)\times\Aut(E')} \ est l'unique orbite ouverte de \m{\ext^1(E,E')}.
Il en d\'ecoule que le faisceau $G$ extension de $E$ par \m{E'} d\'efini par
$\sigma$ est uniquement d\'etermin\'e par $E$ et \m{E'}. Soit \m{\ky}
l'ensemble des classes d'isomorphisme de tels $G$, \m{(E,E')} parcourant
\m{\ky'}.

\medskip

\begin{xlemm}
L'application \ \m{\theta : \ky'\lra\ky} \ associant $G$ \`a \m{(E,E')} est
une bijection.
\end{xlemm}

\begin{proof} Cette application est surjective par d\'efinition de $\ky$.
Montrons qu'elle est injective. Soient \m{(E,E')\in\ky'}, \m{(E_0,E'_0)\in\ky'}
auxquels sont associ\'ees des extensions
$$0\lra E'\lra G\lra E\lra 0,$$
$$0\lra E'_0\lra G_0\lra E_0\lra 0,$$
\m{G_0} \'etant isomorphe \`a $G$. Fixons un isomorphisme \ \m{f:G\lra G_0}.
Puisque

\noindent\m{\hom(E',E_0)=\lbrace 0\rbrace}, cet isomorphisme envoie
\m{E'} dans \m{E'_0}, et induit un morphisme \
 \m{E\lra E_0}. En consid\'erant
\m{f^{-1}}, on voit que \ \m{(E,E')=(E_0,E'_0)}. \end{proof}

\bigskip

\noindent{\bf Remarque : } Pour $E$ dans $\kx$ et \m{E'} dans \m{\kx'},
soient \ \m{p=\dim(\hom(E,E))} \ et \

\noindent
\m{p'=\dim(\hom(E',E'))}, qui ne d\'ependent pas de $E$ et \m{E'} d'apr\`es
la proposition 2.2. Soit
$$e \ = \ \chi(E,E')$$
pour $E$ dans $\kx$ et \m{E'} dans \m{\kx'}, qui est aussi ind\'ependant de $E$
et \m{E'}. Soit \ \m{(E,E')\in\ky'}. On a alors
$$e \ = \ \dim(\hom(E,E'))-\dim(\ext^1(E,E')).$$
Soit \ \m{G=\theta(E,E')}. On a alors
$$\dim(\ext^1(E,E')) \ \leq \ p+p'$$
et
$$\dim(\hom(G,G))=p+p'+e.$$

\newpage

\begin{subsub}{Construction des vari\'et\'es de modules
d'extensions}\end{subsub}

On se place sous l'hypoth\`ese pr\'ec\'edente, c'est-\`a-dire que
$$\hom(E',E)=\lbrace 0\rbrace$$
pour tous $E$ dans $\kx$ et \m{E'} dans \m{\kx'}.

\begin{xtheo}
Soit \m{\ky'} le sous-ensemble de \ \m{\kx\times\kx'} \ constitu\'e des paires
\m{(E,E')} telles que
$$\ext^1(E',E)=\ext^2(E',E')=\ext^2(E,E)=\ext^2(E',E)=\ext^i(E,E')
=\lbrace 0\rbrace$$
pour \ \m{i\geq 2},
et qu'il existe un \'el\'ement $\sigma$ de \m{\ext^1(E,E')} tel que
l'application associ\'ee
$$\phi : \hom(E,E)\oplus\hom(E',E')\lra\ext^1(E,E')$$
soit surjective. Soit $\ky$ l'ensemble des classes d'isomorphisme des faisceaux
qui sont des extensions d\'efinies par de tels $\sigma$. Alors
$\ky$ est un ensemble ouvert, et il existe
une vari\'et\'e de modules fins d\'efinie localement \m{(N,(\ke_i))} pour
$\ky$, $N$ \'etant l'ouvert de \ \m{M\times M'} \ constitu\'e des paires
\m{(x,x')} telles que \ \m{(\ke_x,\ke'_{x'})\in\ky'}.
\end{xtheo}

\begin{proof} Soient \m{y\in \ky}, \m{(E,E')} le point correspondant de
\m{\ky'} et \m{(x_0,x'_0)} le point de $N$ tel que
\ \m{\ke_{x_0}\simeq E}, \m{\ke'_{x'_0}\simeq E'}. D'apr\'es la proposition 2.6
(adapt\'ee aux vari\'et\'es de modules fins d\'efinies localement) il suffit
de d\'efinir un ouvert $U$ de $N$ contenant \m{(x_0,x'_0)} tel qu'on puisse
construire une famille de faisceaux \m{\kg} de $\ky$ param\'etr\'ee par $U$
telle que pour tout \m{(z,z')\in U}, la classe d'isomorphisme de
\m{\kg_{(z,z')}} soit \m{\theta^{-1}(\ke_z,\ke'_{z'})} et que $\kg$ soit une
d\'eformation compl\`ete de \m{\kg_{(z,z')}}.

\medskip

\begin{xlemm}
Soit $m$ un entier. Alors il existe des
fibr\'es vectoriels \m{\kg_0,\ldots,\kg_{d}} sur \m{M} et un fibr\'e vectoriel
$\kh$ sur \m{M\times X}, des entiers \ \m{m<m_0<m_1<\ldots <m_{d}} \ et une
suite exacte de faisceaux coh\'erents sur \ \m{M\times X}
$$0\lra\kh\lra\kg_{d}\ot\ko_X(-m_{d})\lra\ldots \ \ \ \ \ \ \ \ \ \ $$
$$\ \ \ \ \ \ \ \ \ \ \ldots\lra\kg_1\ot\ko_X(-m_1)\lra
\kg_0\ot\ko_X(-m_0)\lra\ke\lra 0.$$
\end{xlemm}

\begin{proof} Puisque $\kx$ est limit\'e il existe un entier \ \m{m_0>m} \ tel
que pour tout \m{y\in M} le faisceau \m{\ke_y(m_0)} soit engendr\'e par ses
sections globales (th\'eor\`eme A de Serre). On prend 
$$\kg_0 \ = \ p_{M*}(\ke(m_0)).$$
Le morphisme d'\'evaluation
$$ev : \kg_0\ot\ko_X(-m_0)\lra\ke$$
est alors surjectif. On construit ensuite \m{G_1} en consid\'erant \m{\ker(ev)}
\`a la place de $\ke$, et on obtient une suite exacte
$$\kg_1\ot\ko_X(-m_1)\lra\kg_0\ot\ko_X(-m_0)\lra\ke\lra 0.$$
Si on poursuit ce proc\'ed\'e jusqu'\`a \m{\kg_{d}} le noyau $\kh$ du morphisme
$$\kg_d\ot\ko_X(-m_d)\lra\kg_{d-1}\ot\ko_X(-m_{d-1})$$
est alors localement libre. \end{proof}

\medskip

Puisque $\ke$ est plat sur $M$, pour tout \m{y\in M}, la restriction de la
suite exacte du lemme pr\'ec\'edent \`a  \ \m{\lbrace y\rbrace\times X} \ est
exacte. La suite exacte du lemme 2.13 se d\'ecompose en suites exactes courtes
$$0\lra\kh\lra\kg_{d}\ot\ko_X(-m_{d})\lra\kv_{d-1}\lra 0,$$
$$0\lra\kv_{d-1}\lra\kg_{d-1}\ot\ko_X(-m_{d-2})\lra\kv_{d-1}\lra 0,$$
$$ . $$
$$ . $$
$$ . $$
$$0\lra\kv_1\lra\kg_1\ot\ko_X(-m_1)\lra\kv_0\lra 0,$$
$$0\lra\kv_0\lra\kg_0\ot\ko_X(-m_0)\lra\ke\lra 0,$$
dont les restrictions \`a  \ \m{\lbrace y\rbrace\times X} \ sont aussi
exactes.

\begin{xlemm}
Si $m$ est assez grand, on a pour tous \m{(y,y')\in N} et \m{i>0}
$$\ext^i(\kg_{ly}\ot\ko_X(-m_l),\ke'_{y'})=\lbrace 0\rbrace$$
pour \m{0\leq l\leq d} et
$$\ext^i(\kh_y,\ke'_{y'})=\lbrace 0\rbrace.$$
\end{xlemm}

\begin{proof}
Puisque \m{\kg_{ly}} est trivial, la premi\`ere \'egalit\'e sera v\'erifi\'ee
d\'es que $m$ est assez grand, d'apr\`es le th\'eor\`eme B de Serre. Pour
d\'emontrer la seconde on consid\`ere les suites exactes pr\'ec\'edentes. On
en d\'eduit, pour tout entier \ \m{i>0} \ des isomorphismes
$$\ext^i(\kh_y,\ke'_{y'})\simeq\ext^{i+1}(\kv_{d-1,y},\ke'_{y'})\simeq\ldots
\simeq\ext^{i+d}(\kv_{0y},\ke'_{y'})=\lbrace 0\rbrace.$$
\end{proof}

\medskip

D'apr\'es le lemme pr\'ec\'edent, les faisceaux
$$F_{d+1} \ = \ p_{N*}(\kh^*\ot p_{M'}^*(\ke')), \ \ \ \
F_i \ = \ p_{N*}(\kg_i^*\ot\ko_X(m_i)\ot p_{M'}^*(\ke')), \ \ (0\leq i\leq d)$$
sont localement libres sur $N$. D'apr\`es ce qui pr\'ec\`ede, la suite exacte
du lemme 2.13 induit un complexe
 $$0\lra F_0 \ \hfl{\alpha_0}{} \ F_1 \ \hfl{\alpha_1}{} \
\ldots\lra F_{d+1}\lra 0$$
exact en degr\'es diff\'erents de $0$ et $1$. Pour tout \m{(y,y')\in N} on a
$$\ker(\alpha_0) \ \simeq \ \hom(\ke_y,\ke'_{y'}), \ \ \ \
\ker(\alpha_1)/\imm(\alpha_0) \ \simeq \ \ext^1(\ke_y,\ke'_{y'}).$$
Le faisceau
$$F'_1 \ = \ \ker(\alpha_1)$$
est localement libre. On obtient finalement un morphisme de fibr\'es vectoriels
sur $N$
$$\alpha : F_0\lra F'_1$$
tel que pour tout \m{(y,y')\in N} on a
$$\ker(\alpha) \ \simeq \ \hom(\ke_y,\ke'_{y'}), \ \ \ \
\coker(\alpha) \ \simeq \ \ext^1(\ke_y,\ke'_{y'}).$$

Soit $W$ un ouvert affine de $N$ contenant \m{(x_0,x'_0)}. Alors on a
$$\ext^i(p_M^\sharp(\kh)_W,p_{M'}^\sharp(\ke')_W)=
\ext^i(p_M^*(\kg_l)_W\ot\ko_X(-m_l),p_{M'}^\sharp(\ke')_W)=\lbrace 0\rbrace$$
pour \m{i>0} et \ \m{0\leq l\leq d} . En effet, les images directes
$$R^ip_{N*}(p_M^*(\kh^*)_W\ot p_{M'}^\sharp(\ke')_W), \ \ \ \
R^ip_{N*}(p_M^*(\kg_l^*)_W\ot\ko_X(m_l)\ot p_{M'}^\sharp(\ke')_W)$$
sont nulles, donc
$$\ext^i(p_M^\sharp(\kh)_W,p_{M'}^\sharp(\ke')_W) \ \simeq \
H^i(p_{N*}(p_M^\sharp(\kh)_W^*\ot p_{M'}^\sharp(\ke')_W)),$$
$$\ext^i(p_M^*(\kg_l)_W\ot\ko_X(-m_l),p_{M'}^\sharp(\ke')_W) \ \simeq \
H^i(p_{N*}(p_M^*(\kg_l^*)_W\ot\ko_X(m_l)\ot p_{M'}^\sharp(\ke')_W)),$$
et ces groupes de cohomologie sont nuls car $W$ est affine.

En employant le m\^eme m\'ethode que pr\'ec\'edemment on en d\'eduit
qu'on a des isomorphismes canoniques
$$\hom(p_M^\sharp(\ke)_W, p_{M'}^\sharp(\ke')_W)\simeq
\ker(H^0(\alpha_{0\mid W})),$$
$$\ext^1(p_M^\sharp(\ke)_W, p_{M'}^\sharp(\ke')_W)\simeq
\coker(H^0(\alpha_{0\mid W})).$$
On obtient aussi que pour tout \m{(x,x')\in W}, le morphisme canonique
$$\pi_x : \ext^1(p_M^\sharp(\ke)_W, p_{M'}^\sharp(\ke')_W)
\lra\ext^1(\ke_x,\ke'_{x'})$$
est surjectif.

Soit \m{s\in\ext^1(p_M^\sharp(\ke)_W, p_{M'}^\sharp(\ke')_W)} \ tel
que le morphisme associ\'e
$$\hom(\ke_{x_0},\ke_{x_0})\oplus\hom(\ke'_{x'_0},\ke'_{x'_0})\lra
\ext^1(\ke_{x_0},\ke'_{x'_0})$$
soit surjectif. Ceci est vrai en tout point d'un ouvert affine \m{U\subset W}
contenant \m{(x_0,x'_0)}. On consid\`ere l'extension
$$0\lra p_M^\sharp(\ke)_{U}\lra\kg\lra p_{M'}^\sharp(\ke')_{U}\lra 0$$
d\'efinie par $s$. Alors $\kg$ est plat sur \m{U} (cf \cite{sga1}, expos\'e IV,
prop. 1.1). C'est la famille de faisceaux de $\ky$ qui convient au voisinage
de \m{(x_0,x'_0)}. Ceci ach\`eve la d\'emonstration du th\'eor\`eme 2.13.
\end{proof}

\bigskip

Le th\'eor\`eme 2.13 est aussi vrai pour des vari\'et\'es de modules fins
d\'efinies localement.

\medskip

\noindent{\bf Remarque : } La d\'emonstration du th\'eor\`eme 2.13 montre que
pour toute sous-vari\'et\'e localement ferm\'ee $Z$ de $N$, tout \
\m{s\in H^0(Z,\coker(\alpha))} \ induit une extension
$$0\lra p_M^\sharp(\ke)_Z\lra\kg\lra p_{M'}^\sharp(\ke')_Z\lra 0.$$

\vskip 1cm

\begin{subsub}{Exemple}\end{subsub}

Soit $C$ une courbe alg\'ebrique projective irr\'eductible lisse, de genre $g$.
Soient $r$, $d$, \m{r'}, \m{d'} des entiers, avec \ \m{r\geq 1}, \m{r'\geq 1},
$r$, $d$ premiers entre eux ainsi que \m{r'}, \m{d'}. On note \m{M(r,d)}
(resp. \m{M(r',d')}) la vari\'et\'e de modules des fibr\'es stables de rang
$r$ et de degr\'e $d$ (resp. de rang \m{r'} et de degr\'e \m{d'}) sur $C$,
et $\ke$ (resp. \m{\ke'}) le fibr\'e universel sur \ \m{M(r,d)\times C} (resp.
\m{M(r',d')\times C}). Soient $F$ un fibr\'e vectoriel alg\'ebrique de rang $2$
sur $C$, et
$$X \ = \ \P(F)$$
la surface r\'egl\'ee sur $C$ des droites de $F$. Soient
$$\pi:X\lra C$$
la projection et \m{\ko(1)} le fibr\'e en droites tautologique sur $X$. On
s'int\'eresse \`a des extensions du type
$$0\lra\pi^*(\ke'_{y'})\ot\ko(1)\lra\ku\lra\pi^*(\ke_y)\lra 0.$$
Ceci est bien un cas particulier de ce qu'on a \'etudi\'e pr\'ec\'edemment car

\noindent\m{(M(r',d'),\pi^\sharp(\ke')\ot\ko(1))} et \m{(M(r,d),
\pi^\sharp(\ke))} sont des vari\'et\'es de modules fins sur $X$.

V\'erifions que les hypoth\`eses du th\'eor\`eme 2.13 sont bien v\'erifi\'ees.
On a pour tout $i$
$$R^i\pi_*((\pi^*(\ke'_{y'})\ot\ko(1))^*\ot\pi^*(\ke_y)) \ = \
{\ke'_{y'}}^*\ot\ke_y\ot R^i\pi_*(\ko(-1)) \ = \lbrace 0\rbrace,$$
donc
$$\ext^i(\pi^*(\ke'_{y'})\ot\ko(1),\pi^*(\ke_y)) \ = \ \lbrace 0\rbrace$$
pour tout $i$. On a d'autre part
$$R^i\pi_*(\pi^*(\ke_y)^*\ot\pi^*(\ke'_{y'})\ot\ko(1)) \ \simeq \
\ke_y^*\ot\ke'_{y'}\ot R^i\pi_*(\ko(1)) \ = \lbrace 0\rbrace$$
si \ \m{i>0}. On a donc
$$\ext^i(\pi^*(\ke_y),\pi^*(\ke'_{y'})\ot\ko(1)) \ \simeq \
H^i(\ke_y^*\ot\ke'_{y'}\ot F^*).$$
En particulier c'est nul si \ \m{i\geq 2}.

\medskip

On peut donc envisager la construction de vari\'et\'es de modules d'extensions
du type pr\'ec\'edent. Soient \ \m{y\in M(r,d)} \ et \ \m{y'\in M(r',d')}. Les
fibr\'es \m{\ke_y} et \m{\ke'_{y'}} \'etant simples, le morphisme


$$\phi : \hom(\pi^*(\ke'_{y'})\ot\ko(1),\pi^*(\ke'_{y'})\ot\ko(1))\oplus
\hom(\pi^*(\ke_y),\pi^*(\ke_y))\ \ \ \ \ \ \ \ \ \ \ \ $$
$$\ \ \ \ \ \ \ \ \ \ \ \ \lra\ext^1(\pi^*(\ke_y),\pi^*(\ke'_{y'})\ot\ko(1))$$
est surjectif si et seulement si
$$\dim(\ext^1(\pi^*(\ke_y),\pi^*(\ke'_{y'})\ot\ko(1))) \ = \
h^1(\ke_y^*\ot\ke_{y'}\ot F^*) \ \leq 1,$$
et si l'extension est non triviale en cas d'\'egalit\'e. On doit donc avoir
$$\chi(\ke_y^*\ot\ke_{y'}\ot F^*) \ \geq \ -1,$$
c'est-\`a-dire
$$\frac{d'}{r'}-\frac{d}{r} \ \geq \ g-1+\frac{deg(F)}{2}-\frac{1}{2rr'}.$$
On obtient alors une vari\'et\'e de modules fins d\'efinie localement
\m{(N,(\kf_i))} pout $\ky$, $N$ \'etant l'ouvert de \ \m{M(r',d')\times M(r,d)}
\ constitu\'e des \m{(y',y)} tels que
$$h^1(\ke_y^*\ot\ke_{y'}\ot F^*) \ \leq 1.$$

\vskip 1.2cm

\begin{sub}{\bf Faisceaux admissibles}\end{sub}

\begin{defin}
On dit qu'un faisceau coh\'erent $E$ sur $X$ est {\em admissible} s'il existe
un ensemble ouvert de faisceaux coh\'erents sur $X$ contenant la classe
d'isomorphisme de $E$ et admettant une vari\'et\'e de modules fins d\'efinie
localement.
\end{defin}

\medskip

D'apr\`es la proposition 2.2, si $E$ est admissible, \m{\Aut(E)} agit
trivialement sur
 \m{\ext^1(E,E)}. Mais cette condition n'est pas suffisante pour
que $E$ soit admissible. On donne au \para 5.3.1 des exemples de fibr\'es
simples sur \m{\P_2} qui se d\'eforment en fibr\'es stables mais
ne sont pas admissibles.

\medskip

On d\'efinit maintenant une relation d'\'equivalence, not\'ee $\equiv$, sur les
classes d'isomorphisme de faisceaux admissibles sur $X$. Cette relation
d'\'equivalence est engendr\'ee par la relation $\sim$ suivante : soient $E$,
$F$ des faisceaux admissibles. Alors \ \m{E\sim F} \ si et seulement si il
existe un
germe de courbe lisse \m{(C,x_0)}, des faisceaux coh\'erents $\ke$, $\kf$ sur
\ \m{C\times X}, plats sur $C$, tels que \ \m{\ke_{x_0}\simeq E},
\m{\kf_{x_0}\simeq F} \ et \ \m{\ke_x\simeq\kf_x} \ si \ \m{x\not =x_0}. Par
d\'efinition on a donc \ \m{E\equiv F} \ s'il existe une suite
\m{G_1,\ldots,G_n} de faisceaux admissibles tels que
$$E\sim G_1, \ \ G_1\sim G_2, \ldots, G_{n-1}\sim G_n, \ \ G_n\sim F.$$

On donne au \para 4.2.3 des exemples de faisceaux simples de rang 1
\'equivalents et non isomorphes.

\vskip 2.5cm

\section{Faisceaux prioritaires g\'en\'eriques instables sur le plan projectif}

\begin{sub}{\bf \'Enonc\'e des principaux r\'esultats}\end{sub}

Les faisceaux prioritaires sur $\P_2$ ont \'et\'e
introduits par A. Hirschowitz et Y. Laszlo dans \cite{hi_la}. Rappelons qu'un
faisceau coh\'erent $\ke$ sur $\P_2$ est dit {\em prioritaire} s'il est sans
torsion et si \ \m{\ext^2(\ke,\ke(-1))=0}. Par exemple les faisceaux
semi-stables au sens de Gieseker-Maruyama sont prioritaires. On s'int\'eresse
ici \`a la structure pr\'ecise du faisceau prioritaire g\'en\'erique de rang
$r$ et de classes de Chern \m{c_1}, \m{c_2} lorsqu'il n'existe pas de faisceau
semi-stable de m\^emes rang et classes de Chern. On en d\'eduira des exemples
de vari\'et\'es de modules fins de faisceaux non simples.

D'apr\`es \cite{hi_la}, le {\em champ}
des faisceaux prioritaires est lisse et irr\'eductible. Les conditions
d'existence des faisceaux prioritaires sont les suivantes : posons
$$\mu \ = \ \frac{c_1}{r}, \ \ \ \Delta \ = \ \frac{1}{r}(c_2 -
\frac{r-1}{2r}c_1^2),$$
(si $\ke$ est un faisceau coh\'erent $\ke$ sur $\P_2$ de rang $r$ et de
classes de Chern \m{c_1}, \m{c_2}, on appelle $\mu=\mu(\ke)$ la {\em pente} de
$\ke$ et $\Delta=\Delta(\ke)$ le {\em discriminant} de $\ke$).
Alors, si \ \m{-1\leq\mu\leq 0}, il existe un faisceau prioritaire de
pente $\mu$ et de discriminant $\Delta$ si et seulement si on a
$$\Delta \ \geq \ \frac{\mu(\mu+1)}{2}.$$
Les conditions d'existence des faisceaux semi-stables sur $\P_2$ sont
rappel\'ees ci-dessous. On peut voir qu'il existe beaucoup de triplets
\m{(r,c_1,c_2)} tels qu'il existe un faisceau prioritaire de rang $r$ et
de classes de Chern \m{c_1}, \m{c_2} mais pas de faisceau semi-stable avec les
m\^emes invariants. On note \m{M(r,c_1,c_2)} la vari\'et\'e de modules des
faisceaux semi-stables de rang $r$ et de classes de Chern \m{c_1}, \m{c_2} sur
\m{\P_2}.

Les conditions d'existence des faisceaux semi-stables sur $\P_2$ (cf.
\cite{dr_lp}) s'expriment en fonction des seules variables $\mu$ et $\Delta$.
On montre qu'il existe une unique fonction $\delta(\mu)$ telle qu'on ait \
\m{\dim(M(r,c_1,c_2)) > 0} \ si et seulement si \ \m{\Delta\geq\delta(\mu)}. La
fonction \m{\delta(\mu)} est d\'ecrite \`a l'aide des {\it fibr\'es
exceptionnels}.

On dit qu'un faisceau coh\'erent $\ke$ sur $\P_2$ est {\it exceptionnel} si
$\ke$ est {\it simple} (c'est-\`a-dire si les seuls endomorphismes de $\ke$
sont les homoth\'eties), et si
$$\ext^1(\ke,\ke) \ = \ \ext^2(\ke,\ke) \ = \ \lbrace 0\rbrace.$$
Un tel faisceau est alors localement libre et stable, et la vari\'et\'e de
modules de faisceaux semi-stables correspondante contient l'unique point $\ke$.
Il existe une infinit\'e d\'enombrable de fibr\'es exceptionnels, et un
proc\'ed\'e simple permet de les obtenir tous \`a partir des fibr\'es en
droites (cf. \cite{dr1}). Notons qu'un fibr\'e exceptionnel est
uniquement d\'etermin\'e par sa pente.
Soit $F$ un fibr\'e exceptionnel. On note \m{x_F} la
plus petite solution de l'\'equation
$$X^2-3X+\frac{1}{rg(F)^2} \ = \ 0.$$
Alors on montre que les intervalles \
\m{\rbrack\mu(F)-x_F,\mu(F)+x_F\lbrack} \ 
constituent une partition de l'ensemble des nombres rationnels. On va d\'ecrire
la fonction \m{\delta(\mu)} sur cet intervalle. Posons
$$P(X) = \frac{X^2}{2}+\frac{3}{2}X+1.$$
Sur l'intervalle \ \m{\rbrack\mu(F)-x_F,\mu(F)\rbrack}, on a
$$\delta(\mu) \ = \ P(\mu-\mu(F))-\frac{1}{2}(1-\frac{1}{rg(F)^2}),$$
et sur \  \m{\lbrack\mu(F),\mu(F)+x_F\lbrack}, on a
$$\delta(\mu) \ = \ P(\mu(F)-\mu)-\frac{1}{2}(1-\frac{1}{rg(F)^2}).$$
On obtient les courbes $G(F)$ et $D(F)$ repr\'esent\'ees sur la figure qui
suit. Ce sont des segments de coniques.

On consid\`ere maintenant la courbe \ \m{\Delta=\delta'(\mu)} \ d\'efinie de
la fa\c con suivante : sur l'intervalle \
\m{\rbrack\mu(F)-x_F,\mu(F)+x_F\lbrack}, on a
$$\delta'(\mu) = \delta(\mu) - \frac{1}{rg(F)^2}(1-\frac{1}{x_F}\mid\mu(F)-\mu
\mid).$$
On obtient ainsi les segments de coniques $G'(F)$ et $D'(F)$. Le point
\m{(\mu(F),\delta'(\mu(F)))} est la paire \m{(\mu,\Delta)} correspondant au
fibr\'e exceptionnel $F$. Le point \m{(\mu(F),\delta(\mu(F)))} est le
sym\'etrique de $F$ par rapport \`a la droite \ \m{\Delta=1/2}. Notons que si
$\mu$ est un nombre rationnel diff\'erent de la pente d'un fibr\'e
exceptionnel, le nombre $\delta'(\mu)$ est irrationnel.
Ces courbes,
sur l'intervalle \ \m{\rbrack\mu(F)-x_F,\mu(F)+x_F\lbrack} \ , sont
repr\'esent\'ees ci-dessous :

\newpage


\vskip 1cm

\setlength{\unitlength}{0.008500in}
\begin{picture}(410,565)(200,235)
\thicklines
\multiput(600,800)(0.00000,-7.98561){70}{\line( 0,-1){  3.993}}
\multiput(810,520)(-7.96117,0.00000){52}{\line(-1, 0){  3.981}}
\put(600,760){\line(-2,-3){160}}
\put(600,760){\line( 2,-3){160}}
\put(760,520){\line(-2,-3){160}}
\put(600,280){\line(-2, 3){160}}
\multiput(440,520)(0.00000,-8.00000){33}{\line( 0,-1){  4.000}}
\multiput(760,520)(0.00000,-8.00000){33}{\line( 0,-1){  4.000}}
\put(465,630){\makebox(0,0)[lb]{\smash{$G(F)$}}}
\put(695,630){\makebox(0,0)[lb]{\smash{$D(F)$}}}
\put(460,395){\makebox(0,0)[lb]{\smash{$G'(F)$}}}
\put(685,395){\makebox(0,0)[lb]{\smash{$  D'(F)$}}}
\put(610,275){\makebox(0,0)[lb]{\smash{$F$}}}
\put(610,760){\makebox(0,0)[lb]{\smash{$P$}}}
\put(770,530){\makebox(0,0)[lb]{\smash{           $\Delta=1/2$}}}
\put(600,235){\makebox(0,0)[lb]{\smash{$\mu=\mu(F)$}}}
\put(415,245){\makebox(0,0)[lb]{\smash{$\mu=\mu(F)-x_F$}}}
\put(735,245){\makebox(0,0)[lb]{\smash{$\mu=\mu(F)+x_F$}}}
\end{picture}
\bigskip

\bigskip

\bigskip

Elles forment un "losange" dont les c\^ot\'es sont des segments de paraboles.
Pour tout point $x$ de $\P_2$, soit \m{\ki_x} le faisceau d'id\'eaux du point
$x$. On a
$$\ext^1(\ki_x,\ko)\simeq\C.$$
Soit \m{\kv_x} l'unique faisceau
extension non triviale de \m{\ki_x} par $\ko$ (ce faisceau est localement
libre). On va d\'emontrer le

\bigskip

\begin{xtheo}
 Soient $r$, \m{c_1}, \m{c_2} des entiers,
avec \m{r\geq 1}, \m{-1<\mu\leq 0},
$$\Delta \ \geq \ \frac{\mu(\mu+1)}{2},$$
et tels que la vari\'et\'e \m{M(r,c_1,c_2)} soit vide.

\medskip

\noindent 1 - Si \ \m{\Delta < \delta'(\mu)}, il existe des fibr\'es
exceptionnels $E_0$, $E_1$, $E_2$, des espaces vectoriels de dimension finie
$M_0$, $M_1$, $M_2$, dont un au plus peut \^etre nul, tels que le faisceau
prioritaire g\'en\'erique de rang $r$ et de classes de Chern $c_1$, $c_2$ soit
isomorphe \`a
$$(E_0\ot M_0)\oplus(E_1\ot M_1)\oplus(E_2\ot M_2).$$

\medskip

\noindent 2 - On suppose que \m{c_1\not = 0} ou \m{c_2>1}.
Si \ \m{\Delta > \delta'(\mu)}, soit $F$ l'unique fibr\'e
exceptionnel tel que \ \m{\mu\in \ \rbrack\mu(F)-x_F,\mu(F)+x_F\lbrack}. Alors
si \ \m{\mu\leq\mu(F)}, l'entier
$$p \ = \ r.rg(F)(P(\mu-\mu(F))-\Delta-\Delta(F))$$
est strictement positif, et le faisceau prioritaire g\'en\'erique de rang $r$
et de classes de Chern $c_1$, $c_2$ est isomorphe \`a une somme directe
$$(F\ot \C^{p})\oplus\ke,$$
o\`u $\ke$ est un fibr\'e semi-stable situ\'e sur la courbe $G(F)$. De m\^eme,
si \ \m{\mu\geq\mu(F)}, l'entier
$$p \ = \ r.rg(F)(P(\mu(F)-\mu)-\Delta-\Delta(F))$$
est strictement positif, et le faisceau prioritaire g\'en\'erique de rang $r$
et de classes de Chern $c_1$, $c_2$ est isomorphe \`a une somme directe
$$(F\ot \C^{p})\oplus\ke,$$
o\`u $\ke$ est un fibr\'e semi-stable situ\'e sur la courbe $D(F)$.

\medskip

\noindent 3 - Si \ \m{c_1=0}, \m{c_2=1}, le faisceau prioritaire g\'en\'erique
de rang $r$ et de classes de Chern $c_1$, $c_2$ est isomorphe \`a une somme
directe du type
$$(\ko\ot\C^{r-2})\oplus\kv_x.$$
\end{xtheo}

\bigskip

Le r\'esultat pr\'ec\'edent apporte des pr\'ecisions sur ce
qui est d\'emontr\'e dans \cite{hi_la}, c'est-\`a-dire que s'il n'existe pas
des faisceau semi-stable de rang $r$ et de classes de Chern \m{c_1}, \m{c_2},
alors deux cas peuvent se produire : la filtration de Harder-Narasimhan du
faisceau prioritaire g\'en\'erique de rang $r$ et de classes de Chern \m{c_1},
\m{c_2} comporte deux termes, ou elle en comporte trois. Dans le premier cas,
un des termes est semi-exceptionnel (c'est-\`a-dire de la forme
\m{F\ot\C^{k}}, avec $F$ exceptionnel), et dans le second cas les trois termes
sont semi-exceptionnels.

La d\'emonstration du th\'eor\`eme 3.1, 1- repose sur le r\'esultat suivant
\hbox{(proposition 3.6) :} soient $r$, \m{c_1}, \m{c_2} des entiers tels
que \ \m{r\geq 2},
\m{\Delta=\delta(\mu)}. Alors, si \ \m{\mu(F)-x_F<\mu\leq\mu(F)}, il existe
un fibr\'e $E$ de rang $r$ et de classes de Chern \m{c_1}, \m{c_2} tel que
\ \m{\ext^1(E,F)=\lbrace 0\rbrace}, et si \ \m{\mu(F)\leq\mu<\mu(F)+x_F}, il
existe un fibr\'e $E$ de rang $r$ et de classes de Chern \m{c_1}, \m{c_2}
tel que
\ \m{\ext^1(F,E)=\lbrace 0\rbrace}. On montre cependant au \para 3.3.2 qu'il
peut exister des fibr\'es stables de rang $r$ et de classes de Chern \m{c_1},
\m{c_2} n'ayant pas ces propri\'et\'es.

\bigskip

Il est possible de pr\'eciser le 1- du th\'eor\`eme 3.1. On rappelle dans le
\para 3.2.1 la notion de {\em triade}, qui est un triplet particulier
\m{(E,F,G)} de fibr\'es exceptionnels. On ne consid\`ere ici que des
triades de fibr\'es exceptionnels dont les pentes sont comprises entre $-1$
et $0$. A la triade \m{(E,F,G)} correspond le {\em triangle} \m{\kt_{(E,F,G)}}
du plan (de coordonn\'ees \m{(\mu,\Delta)}), dont les c\^ot\'es sont des
segments de paraboles et les sommets les points correspondant \`a $E$, $F$ et
$G$. Ce triangle est d\'efini par les in\'equations
$$\Delta\leq P(\mu-\mu(G))-\Delta(G), \ \
\Delta\geq P(\mu-\mu(H)+3)-\Delta(H), \ \
\Delta\leq P(\mu-\mu(E)+3)-\Delta(E),$$
$H$ \'etant le fibr\'e exceptionnel noyau du morphisme d'\'evaluation
\ \m{E\ot\hom(E,F)\lra F}.
Soit {\bf T} l'ensemble des triades de fibr\'es exceptionnels dont les pentes
sont comprises entre $-1$ et $0$. Soit $\ks$ l'ensemble des points
\m{(\mu,\Delta)} du plan tels que
$$-1\leq\mu\leq 0, \ \ \frac{\mu(\mu+1)}{2}\leq\Delta\leq\delta'(\mu).$$
On d\'emontrera le

\bigskip

\bigskip

\begin{xtheo}
1 - Soient \m{(E,F,G)}, \m{(E',F',G')} des
\'el\'ements distincts de {\bf T}. Alors les triangles \m{\kt_{(E,F,G)}} et
\m{\kt_{(E',F',G')}} ont une intersection non vide si et seulement si cette
intersection est un sommet commun ou un c\^ot\'e commun. Dans le premier cas,
les fibr\'es exceptionnels correspondants sont identiques, et dans le second
les paires de fibr\'es exceptionnels correspondantes le sont.

\medskip

\noindent 2 - On a \ \ \ \m{\displaystyle
\ks\ = \ \bigcup_{(E,F,G)\in{\bf T}}\kt_{(E,F,G)}}.

\medskip

\noindent 3 - Soient \m{r,c_1,c_2} des entiers, avec \ \m{r\geq 1},
$$\mu=\frac{r}{c_1}, \ \ \Delta=\frac{1}{r}(c_2-\frac{r-1}{2r}c_1^2).$$
On suppose que \ \m{(\mu,\Delta)\in\kt_{(E,F,G)}}. Soit $H$ le noyau du
morphisme d'\'evaluation \

\noindent \m{E\ot\hom(E,F)\lra F}. Alors
$$m \ = \ r.rg(E).(P(\mu-\mu(E)+3)-\Delta(E)),$$
$$n \ = \ r.rg(H).(P(\mu-\mu(H)+3)-\Delta(H)),$$
$$p \ = \ r.rg(G).(P(\mu-\mu(G))-\Delta(G))$$
sont des entiers positifs ou nuls, et le fibr\'e prioritaire g\'en\'erique
de rang $r$ et de classes de Chern $c_1$, $c_2$ est de la forme
$$(E\ot\C^{m})\oplus(F\ot\C^{n})\oplus(G\ot\C^{p}).$$
\end{xtheo}

\bigskip

Le lieu du plan constitu\'e des points \m{(\Delta,\mu)} tels qu'il existe des
fibr\'es prioritaires de pente $\mu$ et de discriminant $\Delta$ mais pas de
fibr\'es semi-stables ayant les m\^emes invariants est donc d\'ecompos\'e en
triangles et en losanges (dont les c\^ot\'es sont des segments de paraboles).
Ils peuvent tous \^etre construits de la fa\c con suivante (si on se limite aux
pentes comprises entre -1 et 0) : on part du triangle \
\m{\kt_{(\ko(-1),Q^*,\ko)}}, et des losanges correspondant \`a \m{\ko(-1)},
\m{Q^*} et $\ko$. Une r\'ecurrence
permet ensuite de construire \`a partir de ces donn\'ees tous les triangles et
tous les losanges. Supposons construit le triangle \m{\kt_{(E,F,G)}}, et les
losanges correspondant \`a $E$, $F$ et $G$. Alors le morphisme canonique
\[ F\ot\hom(F,G)\lra G \ \ \ \ \ {\rm (resp. \ \ } E\lra\hom(E,F)^*\ot F\
{\rm)}\]
est surjectif (resp. injectif), et son noyau $H$ (resp. conoyau $K$)
est un fibr\'e
exceptionnel. On obtient deux nouveaux triangles \m{\kt_{(E,H,F)}},
\m{\kt_{(F,K,G)}}, et les losanges correspondant \`a $H$ et $K$. La situation
est r\'esum\'ee par la figure suivante :

\vskip 2cm

\includegraphics{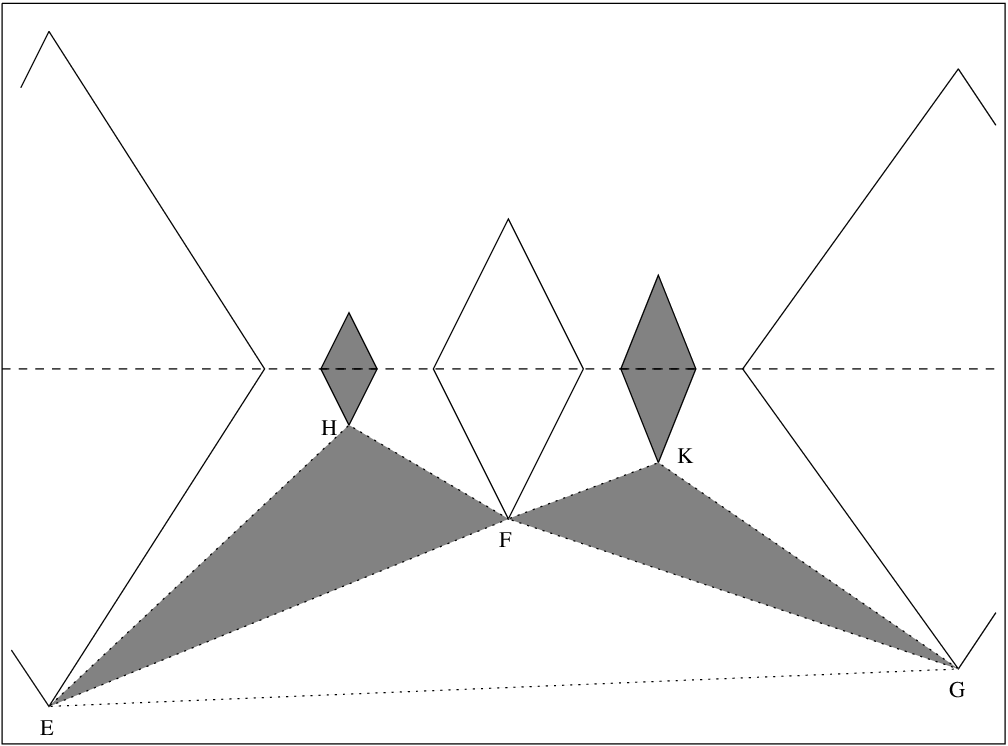}

\vskip 1cm

On construit dans le \para 3.4 des vari\'et\'es de modules fins (d\'efinies
localement ou globalement) constitu\'ees de faisceaux prioritaires de rang
$r$ et de classes de Chern \m{c_1}, \m{c_2}, dans le cas o\`u \
\m{M(r,c_1,c_2)=\emptyset}. Ces vari\'et\'es ne contiennent pas en g\'en\'eral
tous les faisceaux prioritaires de rang $r$ et de classes de Chern \m{c_1},
\m{c_2}, et les faisceaux qui les constituent ne sont pas simples. On donne
en particulier un exemple d'ensemble ouvert de faisceaux prioritaires pour
lequel il existe une vari\'et\'e de modules fins d\'efinie localement mais pas
de vari\'et\'e de modules fins d\'efinie globalement.

\vskip 1.2cm

\begin{sub}{\bf Fibr\'es exceptionnels}\end{sub}

\begin{subsub}{Construction des fibr\'es exceptionnels}\end{subsub}

Les r\'esultats qui suivent ont \'et\'e d\'emontr\'es dans \cite{dr_lp} ou
\cite{dr1}. Un fibr\'e exceptionnel est enti\`erement d\'etermin\'e par sa
pente. Soit $\P$ l'ensemble des pentes de fibr\'es exceptionnels. Si
\m{\alpha\in\P}, on note \m{E_\alpha} le fibr\'e exceptionnel de pente
$\alpha$, et \m{r_\alpha} son rang. On montre que \m{r_\alpha} et
\m{c_1(E_\alpha)} sont premiers entre eux. Soit \ \m{\Delta_\alpha=
\Delta(E_\alpha)}. Alors on a
$$\Delta_\alpha\ = \ \frac{1}{2}(1-\frac{1}{r_\alpha^2}),$$
(ce qui d\'ecoule du fait que \ \m{\chi(E_\alpha,E_\alpha)=1}).

Soit $\D$ l'ensemble des nombres rationnels diadiques, c'est-\`a-dire pouvant
se mettre sous la forme \m{p/2^q}, $p$ et $q$ \'etant des entiers, \m{q\geq 0}.
On a une bijection
$$\epsilon : \D\lra\P.$$
Cette application est enti\`erement d\'etermin\'ee par les propri\'et\'es
suivantes:

\medskip

\noindent - Pour tout entier $k$, on a \ \m{\epsilon(k)=k}.

\noindent - Pour tout entier $k$ et tout \ \m{x\in\D}, on a
\ \m{\epsilon(x+k)=\epsilon(x)+k}.

\noindent - Pour tous entiers $p$, $q$, avec \ \m{q\geq 0}, on a
$$\epsilon(\frac{2p+1}{2^{q+1}}) \ = \
\epsilon(\frac{p}{2^q})\times\epsilon(\frac{p+1}{2^q}),$$
o\`u $\times$ est la loi de composition suivante :
$$\alpha\times\beta\ = \ \frac{\alpha+\beta}{2}+\frac{\Delta_\alpha-\Delta_
\beta}{3+\alpha-\beta}.$$
Cette relation signifie simplement que
$$\chi(E_{\alpha\times\beta},E_\alpha) \ = \
\chi(E_\beta,E_{\alpha\times\beta}) \ = \ 0.$$

\bigskip

La construction des pentes des fibr\'es exceptionnels comprises entre $-1$ et
$0$ se fait donc en partant des pentes $-1$ et $0$, correspondant aux
fibr\'es exceptionnels \m{\ko(-1)} et $\ko$.

On appelle {\em triades} les triplets de fibr\'es exceptionnels de la forme

\noindent\m{(\ko(k),\ko(k+1),\ko(k+2))}, \m{(E_\alpha,E_{\alpha\times\beta},
E_\beta)},\m{(E_{\alpha\times\beta},E_{\beta},E_{\alpha+3})} ou
\m{(E_{\beta-3},E_\alpha,E_{\alpha\times\beta})}, \m{\alpha} et \m{\beta}
\'etant des \'el\'ements de $\P$ de la forme
$$\alpha \ = \ \epsilon(\frac{p}{2^q}), \ \ \ \
\beta \ = \ \epsilon(\frac{p+1}{2^q}).$$
o\`u $p$ et $q$ sont deux entiers avec \ \m{q\geq 0}. Les triades sont
exactement les {\em bases d'h\'elice} de \cite{go_ru}.

On donne maintenant la construction des triades de fibr\'es exceptionnels dont
les pentes sont comprises entre \m{-1} et $0$. Ces triades sont du type
\m{(E_\alpha,E_{\alpha\times\beta},E_\beta)}. La construction
 se fait de la fa\c con suivante, par
r\'ecurrence : on part de la triade \m{(\ko(-1),Q^*,\ko)}, o\`u $Q$ est le
fibr\'e exceptionnel quotient du morphisme canonique \
\m{\ko(-1)\lra\ko\ot H^0(\ko(1))^*}. Supposons la triade \m{(E,F,G)}
construite. Alors on construit les {\it triades adjacentes} \m{(E,H,F)} et
\m{(F,K,G)}. Le fibr\'e $H$ est le noyau du morphisme canonique surjectif
$$F\ot\hom(F,G)\lra G$$
et $K$ est le conoyau du morphisme canonique injectif
$$E\lra F\ot\hom(E,F)^*.$$
De plus, le morphisme canonique
$$E\ot\hom(E,H)\lra H \ \ {\rm \ \ \ (resp. \ } K\lra G\ot\hom(K,G)^*
{\rm \ )}$$
est surjectif (resp. injectif) et son noyau (resp. conoyau) est isomorphe \`a
\m{G(-3)} (resp. \m{E(3)}).

\newpage

\begin{subsub}{Suite spectrale de Beilinson g\'en\'eralis\'ee}\end{subsub}

A toute triade \m{(E,G,F)} et \`a tout faisceau coh\'erent $\ke$ sur $\P_2$
on associe une suite spectrale \m{E^{p,q}_r} de faisceaux coh\'erents sur
$\P_2$, convergeant vers $\ke$ en degr\'e 0 et vers 0 en tout autre degr\'e.
Les termes \m{E^{p,q}_1} \'eventuellement non nuls sont
$$E^{-2,q}_1\simeq H^q(\ke\ot E^*(-3))\ot E, \ \
E^{-1,q}_1\simeq H^q(\ke\ot S^*)\ot G, \ \
E^{0,q}_1\simeq H^q(\ke\ot F^*)\ot F,$$
$S$ d\'esignant le fibr\'e exceptionnel conoyau du morphisme canonique injectif

\noindent\m{G\lra F\ot\hom(G,F)}.

\vskip 0.8cm

\begin{subsub}{S\'erie exceptionnelle associ\'ee \`a un fibr\'e exceptionnel}
\end{subsub}

Soit $F$ un fibr\'e exceptionnel. Les triades comportant $F$ comme terme
de droite sont de la forme \m{(G_n,G_{n+1},F)}, o\`u la suite de fibr\'es
exceptionnels \m{(G_n)} est enti\`erement d\'etermin\'ee par deux de ses
termes cons\'ecutifs, par exemple \m{G_0} et \m{G_1}, par les suites exactes
$$0\lra G_{n-1}\lra (G_n\ot\hom(G_{n-1},G_n)^*)\simeq (G_n\ot\hom(G_n,G_{n+1}))
\lra G_{n+1}\lra 0.$$
On appelle \m{(G_n)} la {\it s\'erie exceptionnelle} \`a gauche associ\'ee
\`a $F$.
Les couples \m{(\mu(G_n),\Delta(G_n))} sont situ\'es sur la conique
d'\'equation
$$\Delta \ = \ P(\mu-\mu(F))-\Delta(F),$$
(ce qui traduit le fait que \ \m{\chi(F,G_n)=0}).

\bigskip

\bigskip


\setlength{\unitlength}{0.240900pt}
\ifx\plotpoint\undefined\newsavebox{\plotpoint}\fi
\begin{picture}(1500,900)(0,0)
\font\gnuplot=cmr10 at 10pt
\gnuplot
\sbox{\plotpoint}{\rule[-0.200pt]{0.400pt}{0.400pt}}%
\put(176,877){\usebox{\plotpoint}}
\multiput(176.58,872.50)(0.493,-1.250){23}{\rule{0.119pt}{1.085pt}}
\multiput(175.17,874.75)(13.000,-29.749){2}{\rule{0.400pt}{0.542pt}}
\multiput(189.58,840.16)(0.492,-1.358){21}{\rule{0.119pt}{1.167pt}}
\multiput(188.17,842.58)(12.000,-29.579){2}{\rule{0.400pt}{0.583pt}}
\multiput(201.58,808.63)(0.493,-1.210){23}{\rule{0.119pt}{1.054pt}}
\multiput(200.17,810.81)(13.000,-28.813){2}{\rule{0.400pt}{0.527pt}}
\multiput(214.58,777.75)(0.493,-1.171){23}{\rule{0.119pt}{1.023pt}}
\multiput(213.17,779.88)(13.000,-27.877){2}{\rule{0.400pt}{0.512pt}}
\multiput(227.58,747.75)(0.493,-1.171){23}{\rule{0.119pt}{1.023pt}}
\multiput(226.17,749.88)(13.000,-27.877){2}{\rule{0.400pt}{0.512pt}}
\multiput(240.58,717.57)(0.492,-1.229){21}{\rule{0.119pt}{1.067pt}}
\multiput(239.17,719.79)(12.000,-26.786){2}{\rule{0.400pt}{0.533pt}}
\multiput(252.58,688.88)(0.493,-1.131){23}{\rule{0.119pt}{0.992pt}}
\multiput(251.17,690.94)(13.000,-26.940){2}{\rule{0.400pt}{0.496pt}}
\multiput(265.58,660.14)(0.493,-1.052){23}{\rule{0.119pt}{0.931pt}}
\put(268,673){$A$}
\multiput(264.17,662.07)(13.000,-25.068){2}{\rule{0.400pt}{0.465pt}}
\multiput(278.58,633.14)(0.493,-1.052){23}{\rule{0.119pt}{0.931pt}}
\multiput(277.17,635.07)(13.000,-25.068){2}{\rule{0.400pt}{0.465pt}}
\multiput(291.58,605.85)(0.492,-1.142){21}{\rule{0.119pt}{1.000pt}}
\multiput(290.17,607.92)(12.000,-24.924){2}{\rule{0.400pt}{0.500pt}}
\multiput(303.58,579.26)(0.493,-1.012){23}{\rule{0.119pt}{0.900pt}}
\multiput(302.17,581.13)(13.000,-24.132){2}{\rule{0.400pt}{0.450pt}}
\multiput(316.58,553.39)(0.493,-0.972){23}{\rule{0.119pt}{0.869pt}}
\multiput(315.17,555.20)(13.000,-23.196){2}{\rule{0.400pt}{0.435pt}}
\multiput(329.58,528.26)(0.492,-1.013){21}{\rule{0.119pt}{0.900pt}}
\multiput(328.17,530.13)(12.000,-22.132){2}{\rule{0.400pt}{0.450pt}}
\multiput(341.58,504.52)(0.493,-0.933){23}{\rule{0.119pt}{0.838pt}}
\multiput(340.17,506.26)(13.000,-22.260){2}{\rule{0.400pt}{0.419pt}}
\multiput(354.58,480.65)(0.493,-0.893){23}{\rule{0.119pt}{0.808pt}}
\multiput(353.17,482.32)(13.000,-21.324){2}{\rule{0.400pt}{0.404pt}}
\multiput(367.58,457.65)(0.493,-0.893){23}{\rule{0.119pt}{0.808pt}}
\multiput(366.17,459.32)(13.000,-21.324){2}{\rule{0.400pt}{0.404pt}}
\put(368,460){\circle*{20}}
\put(383,460){$G_{n-1}$}
\multiput(380.58,434.68)(0.492,-0.884){21}{\rule{0.119pt}{0.800pt}}
\multiput(379.17,436.34)(12.000,-19.340){2}{\rule{0.400pt}{0.400pt}}
\multiput(392.58,413.90)(0.493,-0.814){23}{\rule{0.119pt}{0.746pt}}
\multiput(391.17,415.45)(13.000,-19.451){2}{\rule{0.400pt}{0.373pt}}
\multiput(405.58,392.90)(0.493,-0.814){23}{\rule{0.119pt}{0.746pt}}
\multiput(404.17,394.45)(13.000,-19.451){2}{\rule{0.400pt}{0.373pt}}
\multiput(418.58,372.03)(0.493,-0.774){23}{\rule{0.119pt}{0.715pt}}
\multiput(417.17,373.52)(13.000,-18.515){2}{\rule{0.400pt}{0.358pt}}
\multiput(431.58,351.96)(0.492,-0.798){21}{\rule{0.119pt}{0.733pt}}
\multiput(430.17,353.48)(12.000,-17.478){2}{\rule{0.400pt}{0.367pt}}
\multiput(443.58,333.29)(0.493,-0.695){23}{\rule{0.119pt}{0.654pt}}
\put(445,333){\circle*{20}}
\put(445,333){$G_{n}$}
\multiput(442.17,334.64)(13.000,-16.643){2}{\rule{0.400pt}{0.327pt}}
\multiput(456.58,315.29)(0.493,-0.695){23}{\rule{0.119pt}{0.654pt}}
\multiput(455.17,316.64)(13.000,-16.643){2}{\rule{0.400pt}{0.327pt}}
\multiput(469.58,297.23)(0.492,-0.712){21}{\rule{0.119pt}{0.667pt}}
\multiput(468.17,298.62)(12.000,-15.616){2}{\rule{0.400pt}{0.333pt}}
\multiput(481.58,280.41)(0.493,-0.655){23}{\rule{0.119pt}{0.623pt}}
\multiput(480.17,281.71)(13.000,-15.707){2}{\rule{0.400pt}{0.312pt}}
\multiput(494.58,263.54)(0.493,-0.616){23}{\rule{0.119pt}{0.592pt}}
\multiput(493.17,264.77)(13.000,-14.771){2}{\rule{0.400pt}{0.296pt}}
\multiput(507.58,247.67)(0.493,-0.576){23}{\rule{0.119pt}{0.562pt}}
\multiput(506.17,248.83)(13.000,-13.834){2}{\rule{0.400pt}{0.281pt}}
\multiput(520.58,232.65)(0.492,-0.582){21}{\rule{0.119pt}{0.567pt}}
\multiput(519.17,233.82)(12.000,-12.824){2}{\rule{0.400pt}{0.283pt}}
\multiput(532.58,218.80)(0.493,-0.536){23}{\rule{0.119pt}{0.531pt}}
\multiput(531.17,219.90)(13.000,-12.898){2}{\rule{0.400pt}{0.265pt}}
\multiput(545.58,204.80)(0.493,-0.536){23}{\rule{0.119pt}{0.531pt}}
\multiput(544.17,205.90)(13.000,-12.898){2}{\rule{0.400pt}{0.265pt}}
\multiput(558.00,191.92)(0.539,-0.492){21}{\rule{0.533pt}{0.119pt}}
\multiput(558.00,192.17)(11.893,-12.000){2}{\rule{0.267pt}{0.400pt}}
\put(560,192){\circle*{20}}
\put(575,192){$G_{n+1}$}
\multiput(571.00,179.92)(0.496,-0.492){21}{\rule{0.500pt}{0.119pt}}
\multiput(571.00,180.17)(10.962,-12.000){2}{\rule{0.250pt}{0.400pt}}
\multiput(583.00,167.92)(0.590,-0.492){19}{\rule{0.573pt}{0.118pt}}
\multiput(583.00,168.17)(11.811,-11.000){2}{\rule{0.286pt}{0.400pt}}
\multiput(596.00,156.92)(0.590,-0.492){19}{\rule{0.573pt}{0.118pt}}
\multiput(596.00,157.17)(11.811,-11.000){2}{\rule{0.286pt}{0.400pt}}
\multiput(609.00,145.92)(0.600,-0.491){17}{\rule{0.580pt}{0.118pt}}
\multiput(609.00,146.17)(10.796,-10.000){2}{\rule{0.290pt}{0.400pt}}
\multiput(621.00,135.93)(0.728,-0.489){15}{\rule{0.678pt}{0.118pt}}
\multiput(621.00,136.17)(11.593,-9.000){2}{\rule{0.339pt}{0.400pt}}
\multiput(634.00,126.93)(0.824,-0.488){13}{\rule{0.750pt}{0.117pt}}
\multiput(634.00,127.17)(11.443,-8.000){2}{\rule{0.375pt}{0.400pt}}
\multiput(647.00,118.93)(0.824,-0.488){13}{\rule{0.750pt}{0.117pt}}
\multiput(647.00,119.17)(11.443,-8.000){2}{\rule{0.375pt}{0.400pt}}
\multiput(660.00,110.93)(0.758,-0.488){13}{\rule{0.700pt}{0.117pt}}
\multiput(660.00,111.17)(10.547,-8.000){2}{\rule{0.350pt}{0.400pt}}
\multiput(672.00,102.93)(1.123,-0.482){9}{\rule{0.967pt}{0.116pt}}
\multiput(672.00,103.17)(10.994,-6.000){2}{\rule{0.483pt}{0.400pt}}
\multiput(685.00,96.93)(1.123,-0.482){9}{\rule{0.967pt}{0.116pt}}
\multiput(685.00,97.17)(10.994,-6.000){2}{\rule{0.483pt}{0.400pt}}
\multiput(698.00,90.93)(1.123,-0.482){9}{\rule{0.967pt}{0.116pt}}
\multiput(698.00,91.17)(10.994,-6.000){2}{\rule{0.483pt}{0.400pt}}
\multiput(711.00,84.94)(1.651,-0.468){5}{\rule{1.300pt}{0.113pt}}
\multiput(711.00,85.17)(9.302,-4.000){2}{\rule{0.650pt}{0.400pt}}
\multiput(723.00,80.94)(1.797,-0.468){5}{\rule{1.400pt}{0.113pt}}
\multiput(723.00,81.17)(10.094,-4.000){2}{\rule{0.700pt}{0.400pt}}
\multiput(736.00,76.95)(2.695,-0.447){3}{\rule{1.833pt}{0.108pt}}
\multiput(736.00,77.17)(9.195,-3.000){2}{\rule{0.917pt}{0.400pt}}
\multiput(749.00,73.95)(2.472,-0.447){3}{\rule{1.700pt}{0.108pt}}
\multiput(749.00,74.17)(8.472,-3.000){2}{\rule{0.850pt}{0.400pt}}
\put(761,70.17){\rule{2.700pt}{0.400pt}}
\multiput(761.00,71.17)(7.396,-2.000){2}{\rule{1.350pt}{0.400pt}}
\put(774,68.67){\rule{3.132pt}{0.400pt}}
\multiput(774.00,69.17)(6.500,-1.000){2}{\rule{1.566pt}{0.400pt}}
\put(787,67.67){\rule{3.132pt}{0.400pt}}
\multiput(787.00,68.17)(6.500,-1.000){2}{\rule{1.566pt}{0.400pt}}
\put(812,67.67){\rule{3.132pt}{0.400pt}}
\multiput(812.00,67.17)(6.500,1.000){2}{\rule{1.566pt}{0.400pt}}
\put(825,68.67){\rule{3.132pt}{0.400pt}}
\multiput(825.00,68.17)(6.500,1.000){2}{\rule{1.566pt}{0.400pt}}
\put(838,70.17){\rule{2.700pt}{0.400pt}}
\multiput(838.00,69.17)(7.396,2.000){2}{\rule{1.350pt}{0.400pt}}
\multiput(851.00,72.61)(2.472,0.447){3}{\rule{1.700pt}{0.108pt}}
\multiput(851.00,71.17)(8.472,3.000){2}{\rule{0.850pt}{0.400pt}}
\multiput(863.00,75.61)(2.695,0.447){3}{\rule{1.833pt}{0.108pt}}
\multiput(863.00,74.17)(9.195,3.000){2}{\rule{0.917pt}{0.400pt}}
\multiput(876.00,78.60)(1.797,0.468){5}{\rule{1.400pt}{0.113pt}}
\multiput(876.00,77.17)(10.094,4.000){2}{\rule{0.700pt}{0.400pt}}
\multiput(889.00,82.60)(1.651,0.468){5}{\rule{1.300pt}{0.113pt}}
\multiput(889.00,81.17)(9.302,4.000){2}{\rule{0.650pt}{0.400pt}}
\multiput(901.00,86.59)(1.123,0.482){9}{\rule{0.967pt}{0.116pt}}
\multiput(901.00,85.17)(10.994,6.000){2}{\rule{0.483pt}{0.400pt}}
\multiput(914.00,92.59)(1.123,0.482){9}{\rule{0.967pt}{0.116pt}}
\multiput(914.00,91.17)(10.994,6.000){2}{\rule{0.483pt}{0.400pt}}
\multiput(927.00,98.59)(1.123,0.482){9}{\rule{0.967pt}{0.116pt}}
\multiput(927.00,97.17)(10.994,6.000){2}{\rule{0.483pt}{0.400pt}}
\multiput(940.00,104.59)(0.758,0.488){13}{\rule{0.700pt}{0.117pt}}
\multiput(940.00,103.17)(10.547,8.000){2}{\rule{0.350pt}{0.400pt}}
\multiput(952.00,112.59)(0.824,0.488){13}{\rule{0.750pt}{0.117pt}}
\multiput(952.00,111.17)(11.443,8.000){2}{\rule{0.375pt}{0.400pt}}
\multiput(965.00,120.59)(0.824,0.488){13}{\rule{0.750pt}{0.117pt}}
\multiput(965.00,119.17)(11.443,8.000){2}{\rule{0.375pt}{0.400pt}}
\multiput(978.00,128.59)(0.728,0.489){15}{\rule{0.678pt}{0.118pt}}
\multiput(978.00,127.17)(11.593,9.000){2}{\rule{0.339pt}{0.400pt}}
\multiput(991.00,137.58)(0.600,0.491){17}{\rule{0.580pt}{0.118pt}}
\multiput(991.00,136.17)(10.796,10.000){2}{\rule{0.290pt}{0.400pt}}
\multiput(1003.00,147.58)(0.590,0.492){19}{\rule{0.573pt}{0.118pt}}
\multiput(1003.00,146.17)(11.811,11.000){2}{\rule{0.286pt}{0.400pt}}
\multiput(1016.00,158.58)(0.590,0.492){19}{\rule{0.573pt}{0.118pt}}
\multiput(1016.00,157.17)(11.811,11.000){2}{\rule{0.286pt}{0.400pt}}
\multiput(1029.00,169.58)(0.496,0.492){21}{\rule{0.500pt}{0.119pt}}
\multiput(1029.00,168.17)(10.962,12.000){2}{\rule{0.250pt}{0.400pt}}
\multiput(1041.00,181.58)(0.539,0.492){21}{\rule{0.533pt}{0.119pt}}
\multiput(1041.00,180.17)(11.893,12.000){2}{\rule{0.267pt}{0.400pt}}
\multiput(1054.58,193.00)(0.493,0.536){23}{\rule{0.119pt}{0.531pt}}
\multiput(1053.17,193.00)(13.000,12.898){2}{\rule{0.400pt}{0.265pt}}
\multiput(1067.58,207.00)(0.493,0.536){23}{\rule{0.119pt}{0.531pt}}
\multiput(1066.17,207.00)(13.000,12.898){2}{\rule{0.400pt}{0.265pt}}
\multiput(1080.58,221.00)(0.492,0.582){21}{\rule{0.119pt}{0.567pt}}
\multiput(1079.17,221.00)(12.000,12.824){2}{\rule{0.400pt}{0.283pt}}
\multiput(1092.58,235.00)(0.493,0.576){23}{\rule{0.119pt}{0.562pt}}
\multiput(1091.17,235.00)(13.000,13.834){2}{\rule{0.400pt}{0.281pt}}
\multiput(1105.58,250.00)(0.493,0.616){23}{\rule{0.119pt}{0.592pt}}
\multiput(1104.17,250.00)(13.000,14.771){2}{\rule{0.400pt}{0.296pt}}
\multiput(1118.58,266.00)(0.493,0.655){23}{\rule{0.119pt}{0.623pt}}
\multiput(1117.17,266.00)(13.000,15.707){2}{\rule{0.400pt}{0.312pt}}
\multiput(1131.58,283.00)(0.492,0.712){21}{\rule{0.119pt}{0.667pt}}
\multiput(1130.17,283.00)(12.000,15.616){2}{\rule{0.400pt}{0.333pt}}
\multiput(1143.58,300.00)(0.493,0.695){23}{\rule{0.119pt}{0.654pt}}
\multiput(1142.17,300.00)(13.000,16.643){2}{\rule{0.400pt}{0.327pt}}
\multiput(1156.58,318.00)(0.493,0.695){23}{\rule{0.119pt}{0.654pt}}
\multiput(1155.17,318.00)(13.000,16.643){2}{\rule{0.400pt}{0.327pt}}
\multiput(1169.58,336.00)(0.492,0.798){21}{\rule{0.119pt}{0.733pt}}
\multiput(1168.17,336.00)(12.000,17.478){2}{\rule{0.400pt}{0.367pt}}
\multiput(1181.58,355.00)(0.493,0.774){23}{\rule{0.119pt}{0.715pt}}
\multiput(1180.17,355.00)(13.000,18.515){2}{\rule{0.400pt}{0.358pt}}
\multiput(1194.58,375.00)(0.493,0.814){23}{\rule{0.119pt}{0.746pt}}
\multiput(1193.17,375.00)(13.000,19.451){2}{\rule{0.400pt}{0.373pt}}
\multiput(1207.58,396.00)(0.493,0.814){23}{\rule{0.119pt}{0.746pt}}
\multiput(1206.17,396.00)(13.000,19.451){2}{\rule{0.400pt}{0.373pt}}
\multiput(1220.58,417.00)(0.492,0.884){21}{\rule{0.119pt}{0.800pt}}
\multiput(1219.17,417.00)(12.000,19.340){2}{\rule{0.400pt}{0.400pt}}
\multiput(1232.58,438.00)(0.493,0.893){23}{\rule{0.119pt}{0.808pt}}
\multiput(1231.17,438.00)(13.000,21.324){2}{\rule{0.400pt}{0.404pt}}
\multiput(1245.58,461.00)(0.493,0.893){23}{\rule{0.119pt}{0.808pt}}
\multiput(1244.17,461.00)(13.000,21.324){2}{\rule{0.400pt}{0.404pt}}
\multiput(1258.58,484.00)(0.493,0.933){23}{\rule{0.119pt}{0.838pt}}
\multiput(1257.17,484.00)(13.000,22.260){2}{\rule{0.400pt}{0.419pt}}
\multiput(1271.58,508.00)(0.492,1.013){21}{\rule{0.119pt}{0.900pt}}
\multiput(1270.17,508.00)(12.000,22.132){2}{\rule{0.400pt}{0.450pt}}
\multiput(1283.58,532.00)(0.493,0.972){23}{\rule{0.119pt}{0.869pt}}
\multiput(1282.17,532.00)(13.000,23.196){2}{\rule{0.400pt}{0.435pt}}
\multiput(1296.58,557.00)(0.493,1.012){23}{\rule{0.119pt}{0.900pt}}
\multiput(1295.17,557.00)(13.000,24.132){2}{\rule{0.400pt}{0.450pt}}
\multiput(1309.58,583.00)(0.492,1.142){21}{\rule{0.119pt}{1.000pt}}
\multiput(1308.17,583.00)(12.000,24.924){2}{\rule{0.400pt}{0.500pt}}
\multiput(1321.58,610.00)(0.493,1.052){23}{\rule{0.119pt}{0.931pt}}
\multiput(1320.17,610.00)(13.000,25.068){2}{\rule{0.400pt}{0.465pt}}
\multiput(1334.58,637.00)(0.493,1.052){23}{\rule{0.119pt}{0.931pt}}
\multiput(1333.17,637.00)(13.000,25.068){2}{\rule{0.400pt}{0.465pt}}
\multiput(1347.58,664.00)(0.493,1.131){23}{\rule{0.119pt}{0.992pt}}
\multiput(1346.17,664.00)(13.000,26.940){2}{\rule{0.400pt}{0.496pt}}
\put(1310,673){$B$}
\multiput(1360.58,693.00)(0.492,1.229){21}{\rule{0.119pt}{1.067pt}}
\multiput(1359.17,693.00)(12.000,26.786){2}{\rule{0.400pt}{0.533pt}}
\multiput(1372.58,722.00)(0.493,1.171){23}{\rule{0.119pt}{1.023pt}}
\multiput(1371.17,722.00)(13.000,27.877){2}{\rule{0.400pt}{0.512pt}}
\multiput(1385.58,752.00)(0.493,1.171){23}{\rule{0.119pt}{1.023pt}}
\multiput(1384.17,752.00)(13.000,27.877){2}{\rule{0.400pt}{0.512pt}}
\multiput(1398.58,782.00)(0.493,1.210){23}{\rule{0.119pt}{1.054pt}}
\multiput(1397.17,782.00)(13.000,28.813){2}{\rule{0.400pt}{0.527pt}}
\multiput(1411.58,813.00)(0.492,1.358){21}{\rule{0.119pt}{1.167pt}}
\multiput(1410.17,813.00)(12.000,29.579){2}{\rule{0.400pt}{0.583pt}}
\multiput(1423.58,845.00)(0.493,1.250){23}{\rule{0.119pt}{1.085pt}}
\multiput(1422.17,845.00)(13.000,29.749){2}{\rule{0.400pt}{0.542pt}}
\put(800.0,68.0){\rule[-0.200pt]{2.891pt}{0.400pt}}
\put(176,661){\usebox{\plotpoint}}
\put(176.00,661.00){\usebox{\plotpoint}}
\put(196.76,661.00){\usebox{\plotpoint}}
\multiput(201,661)(20.756,0.000){0}{\usebox{\plotpoint}}
\put(217.51,661.00){\usebox{\plotpoint}}
\put(238.27,661.00){\usebox{\plotpoint}}
\multiput(240,661)(20.756,0.000){0}{\usebox{\plotpoint}}
\put(259.02,661.00){\usebox{\plotpoint}}
\multiput(265,661)(20.756,0.000){0}{\usebox{\plotpoint}}
\put(279.78,661.00){\usebox{\plotpoint}}
\put(300.53,661.00){\usebox{\plotpoint}}
\multiput(303,661)(20.756,0.000){0}{\usebox{\plotpoint}}
\put(321.29,661.00){\usebox{\plotpoint}}
\multiput(329,661)(20.756,0.000){0}{\usebox{\plotpoint}}
\put(342.04,661.00){\usebox{\plotpoint}}
\put(362.80,661.00){\usebox{\plotpoint}}
\multiput(367,661)(20.756,0.000){0}{\usebox{\plotpoint}}
\put(383.55,661.00){\usebox{\plotpoint}}
\put(404.31,661.00){\usebox{\plotpoint}}
\multiput(405,661)(20.756,0.000){0}{\usebox{\plotpoint}}
\put(425.07,661.00){\usebox{\plotpoint}}
\multiput(431,661)(20.756,0.000){0}{\usebox{\plotpoint}}
\put(445.82,661.00){\usebox{\plotpoint}}
\put(466.58,661.00){\usebox{\plotpoint}}
\multiput(469,661)(20.756,0.000){0}{\usebox{\plotpoint}}
\put(487.33,661.00){\usebox{\plotpoint}}
\multiput(494,661)(20.756,0.000){0}{\usebox{\plotpoint}}
\put(508.09,661.00){\usebox{\plotpoint}}
\put(528.84,661.00){\usebox{\plotpoint}}
\multiput(532,661)(20.756,0.000){0}{\usebox{\plotpoint}}
\put(549.60,661.00){\usebox{\plotpoint}}
\put(570.35,661.00){\usebox{\plotpoint}}
\multiput(571,661)(20.756,0.000){0}{\usebox{\plotpoint}}
\put(591.11,661.00){\usebox{\plotpoint}}
\multiput(596,661)(20.756,0.000){0}{\usebox{\plotpoint}}
\put(611.87,661.00){\usebox{\plotpoint}}
\put(632.62,661.00){\usebox{\plotpoint}}
\multiput(634,661)(20.756,0.000){0}{\usebox{\plotpoint}}
\put(653.38,661.00){\usebox{\plotpoint}}
\multiput(660,661)(20.756,0.000){0}{\usebox{\plotpoint}}
\put(674.13,661.00){\usebox{\plotpoint}}
\put(694.89,661.00){\usebox{\plotpoint}}
\multiput(698,661)(20.756,0.000){0}{\usebox{\plotpoint}}
\put(715.64,661.00){\usebox{\plotpoint}}
\multiput(723,661)(20.756,0.000){0}{\usebox{\plotpoint}}
\put(736.40,661.00){\usebox{\plotpoint}}
\put(757.15,661.00){\usebox{\plotpoint}}
\multiput(761,661)(20.756,0.000){0}{\usebox{\plotpoint}}
\put(777.91,661.00){\usebox{\plotpoint}}
\put(798.66,661.00){\usebox{\plotpoint}}
\multiput(800,661)(20.756,0.000){0}{\usebox{\plotpoint}}
\put(819.42,661.00){\usebox{\plotpoint}}
\multiput(825,661)(20.756,0.000){0}{\usebox{\plotpoint}}
\put(840.18,661.00){\usebox{\plotpoint}}
\put(860.93,661.00){\usebox{\plotpoint}}
\multiput(863,661)(20.756,0.000){0}{\usebox{\plotpoint}}
\put(881.69,661.00){\usebox{\plotpoint}}
\multiput(889,661)(20.756,0.000){0}{\usebox{\plotpoint}}
\put(902.44,661.00){\usebox{\plotpoint}}
\put(923.20,661.00){\usebox{\plotpoint}}
\multiput(927,661)(20.756,0.000){0}{\usebox{\plotpoint}}
\put(943.95,661.00){\usebox{\plotpoint}}
\put(964.71,661.00){\usebox{\plotpoint}}
\multiput(965,661)(20.756,0.000){0}{\usebox{\plotpoint}}
\put(985.46,661.00){\usebox{\plotpoint}}
\multiput(991,661)(20.756,0.000){0}{\usebox{\plotpoint}}
\put(1006.22,661.00){\usebox{\plotpoint}}
\put(1026.98,661.00){\usebox{\plotpoint}}
\multiput(1029,661)(20.756,0.000){0}{\usebox{\plotpoint}}
\put(1047.73,661.00){\usebox{\plotpoint}}
\multiput(1054,661)(20.756,0.000){0}{\usebox{\plotpoint}}
\put(1068.49,661.00){\usebox{\plotpoint}}
\put(1089.24,661.00){\usebox{\plotpoint}}
\multiput(1092,661)(20.756,0.000){0}{\usebox{\plotpoint}}
\put(1110.00,661.00){\usebox{\plotpoint}}
\put(1130.75,661.00){\usebox{\plotpoint}}
\multiput(1131,661)(20.756,0.000){0}{\usebox{\plotpoint}}
\put(1151.51,661.00){\usebox{\plotpoint}}
\multiput(1156,661)(20.756,0.000){0}{\usebox{\plotpoint}}
\put(1172.26,661.00){\usebox{\plotpoint}}
\put(1193.02,661.00){\usebox{\plotpoint}}
\multiput(1194,661)(20.756,0.000){0}{\usebox{\plotpoint}}
\put(1213.77,661.00){\usebox{\plotpoint}}
\multiput(1220,661)(20.756,0.000){0}{\usebox{\plotpoint}}
\put(1234.53,661.00){\usebox{\plotpoint}}
\put(1255.29,661.00){\usebox{\plotpoint}}
\multiput(1258,661)(20.756,0.000){0}{\usebox{\plotpoint}}
\put(1276.04,661.00){\usebox{\plotpoint}}
\multiput(1283,661)(20.756,0.000){0}{\usebox{\plotpoint}}
\put(1296.80,661.00){\usebox{\plotpoint}}
\put(1317.55,661.00){\usebox{\plotpoint}}
\multiput(1321,661)(20.756,0.000){0}{\usebox{\plotpoint}}
\put(1338.31,661.00){\usebox{\plotpoint}}
\put(1359.06,661.00){\usebox{\plotpoint}}
\multiput(1360,661)(20.756,0.000){0}{\usebox{\plotpoint}}
\put(1379.82,661.00){\usebox{\plotpoint}}
\multiput(1385,661)(20.756,0.000){0}{\usebox{\plotpoint}}
\put(1400.57,661.00){\usebox{\plotpoint}}
\put(1421.33,661.00){\usebox{\plotpoint}}
\multiput(1423,661)(20.756,0.000){0}{\usebox{\plotpoint}}
\put(1436,661){\usebox{\plotpoint}}
\put(1480,650){$\Delta=1/2$}
\end{picture}

\bigskip

Dans la figure ci-dessus, les points $A$ et $B$ sont les intersections de
cette conique avec la droite d'\'equation \ \m{\Delta=1/2}. On a
$$\lim_{n\rightarrow -\infty}G_n \ = \ A \ \ \ \ {\rm et} \ \ \ \
\lim_{n\rightarrow\infty}G_n \ =\ B.$$
Remarquons que \ \m{\mu(B)-\mu(A)<3}.

Si \m{F=\ko}, il existe une unique paire \m{(G_n,G_{n+1})} telle que \
\m{\mu(G_{n+1})-\mu(G_n)\geq 1}, c'est \m{(\ko(-2),\ko(-1))}. Supposons que
\ \m{-1<\mu(F)<0}. Il existe alors une unique triade de la forme
\m{(E,F,G)}, avec \ \m{-1\leq\mu(E)<\mu(G)\leq 0}. On en d\'eduit que
\m{(G(-3),E)} est une des paires \m{(G_n,G_{n+1})}. On peut supposer que
\ \m{(G(-3),E)=(G_0,G_1)}. On a \ \m{\mu(G_1)-\mu(G_0)\geq 2}, et \m{(G_0,G_1)}
est l'unique paire  \m{(G_n,G_{n+1})} telle que

\noindent\m{\mu(G_{n+1})-\mu(G_n)\geq 1}. On l'appelle la paire {\it initiale}
de la s\'erie \m{(G_n)}.

\bigskip

\begin{xlemm}
Le fibr\'e vectoriel \m{G_0^*\ot G_1} est engendr\'e par ses sections globales.
\end{xlemm}

\begin{proof} D'apr\`es la construction de \m{(G_0,G_1)}, il suffit de prouver
le r\'esultat \hbox{suivant :} si \m{(A,B,C)} est une triade de fibr\'es
exceptionnels telle que \ \m{\mu(C)-\mu(A)\leq 1}, les fibr\'es
\m{B^*\ot A(3)}, \m{C^*\ot B(3)} et \m{C^*\ot A(3)} sont engendr\'es par leurs
sections globales. On d\'emontre cela par r\'ecurrence : il faut montrer que
si c'est vrai pour une triade, c'est vrai pour les deux triades adjacentes.
Supposons que ce soit vrai pour \m{(A,B,C)}. Soient $H$ le noyau du morphisme
canonique surjectif
$$B\ot\hom(B,C)\lra C$$
et K le conoyau du morphisme canonique injectif
$$A\lra B\ot\hom(A,B)^*.$$
Il faut montrer que le r\'esultat est vrai pour les triades \m{(A,H,B)} et
\m{(B,K,C)}. En consid\'erant la triade {\it duale}
\m{(C^*(-1),B^*(-1),A^*(-1))}, on voit qu'il suffit de consid\'erer
\m{(A,H,B)}. On a une suite exacte
$$0\lra H\lra B\ot\hom(B,C)\lra C\lra 0.$$
On en d\'eduit un morphisme surjectif
$$B^*(3)\ot A\ot\hom(B,C)^*\lra H^*(3)\ot A.$$
Puisque \m{B^*(3)\ot A} est engendr\'e par ses sections globales (hypoth\`ese
de r\'ecurrence), il en est de m\^eme de \m{H^*(3)\ot A}. On a d'autre
part une suite exacte
$$0\lra C(-3)\lra A\ot\hom(C(-3),A)^*\lra H\lra 0,$$
d'o\`u on d\'eduit un morphisme surjectif
$$B^*(3)\ot A\ot\hom(C(-3),A)^*\lra B^*(3)\ot H,$$
d'o\`u on d\'eduit que \m{B^*(3)\ot H} est engendr\'e par ses sections
globales. \end{proof}

\bigskip

\begin{xlemm}
Pour tout entier $n$, on a \ \m{n\geq 1} \ si et seulement si pour tous entiers
$a$, $b$, $c$ positifs ou nuls, le fibr\'e vectoriel
$$(G_n\ot\C^{a})\oplus(G_{n+1}\ot\C^{b})\oplus(F\ot\C^{c})$$
est prioritaire.
\end{xlemm}

\begin{proof} Imm\'ediat. \end{proof}

\medskip

On d\'efinit de m\^eme la {\em s\'erie exceptionnelle \`a droite} \m{(H_n)}
associ\'ee \`a $F$. On a \m{H_n=G_n(3)} pour tout $n$.

\newpage

\begin{subsub}{\'Etude de {\bf T}}\end{subsub}

L'ensemble {\bf T} est construit comme une union croissante de sous-ensembles
$$T_0=\lbrace(\ko(-1),Q^*,\ko)\rbrace\subset T_1\subset\ldots T_n\subset
T_{n+1}\subset\ldots, \ \ \ \
{\bf T}=\bigcup_{n\geq 0}T_n,$$
o\`u $T_n$ est l'ensemble des triades \m{(E_\alpha,E_{\alpha\times\beta},
E_\beta)}, $\alpha$, $\beta$ \'etant de la forme
\ $\alpha=\epsilon(\frac{p}{2^n}), \ \ \beta=\epsilon(\frac{p+1}{2^n})$,
avec $p$ entier. Si $n>0$, les triades de \m{T_n\backslash T_{n-1}} forment
une suite \m{t_0^{(n)}}, \ldots, \m{t_{2^n-1}^{(n)}},
$$t_i^{(n)}\ = \ (E_{\alpha(\frac{i}{2^n})},E_{\alpha(\frac{2i+1}{2^{n+1}})},
E_{\alpha(\frac{i+1}{2^n})}).$$
On a
$$\mu(E_{\alpha(\frac{i}{2^n})}) \ < \ \mu(E_{\alpha(\frac{2i+1}{2^{n+1}})})
\ < \ \mu(E_{\alpha(\frac{i+1}{2^n})}),$$
et dans le plan de coordonn\'ees \m{(\mu,\Delta)},
\m{E_{\alpha(\frac{2i+1}{2^{n+1}})}} est situ\'e au dessus de la droite
\m{E_{\alpha(\frac{i}{2^n})}E_{\alpha(\frac{i+1}{2^n})}}.


\null

\vspace{-2.8cm}

\setlength{\unitlength}{0.240900pt}
\ifx\plotpoint\undefined\newsavebox{\plotpoint}\fi
\begin{picture}(1500,900)(0,0)
\font\gnuplot=cmr10 at 10pt
\gnuplot
\sbox{\plotpoint}{\rule[-0.200pt]{0.400pt}{0.400pt}}%
\put(233,68){\usebox{\plotpoint}}
\multiput(233.00,68.58)(3.543,0.500){321}{\rule{2.930pt}{0.120pt}}
\multiput(233.00,67.17)(1139.919,162.000){2}{\rule{1.465pt}{0.400pt}}
\put(233,68){\usebox{\plotpoint}}
\multiput(233.00,68.58)(0.591,0.500){967}{\rule{0.573pt}{0.120pt}}
\multiput(233.00,67.17)(571.812,485.000){2}{\rule{0.286pt}{0.400pt}}
\put(806,553){\usebox{\plotpoint}}
\multiput(806.00,551.92)(0.887,-0.500){643}{\rule{0.810pt}{0.120pt}}
\multiput(806.00,552.17)(571.320,-323.000){2}{\rule{0.405pt}{0.400pt}}
\put(233,68){\usebox{\plotpoint}}
\multiput(233,68)(9.029,18.689){20}{\usebox{\plotpoint}}
\put(405,424){\usebox{\plotpoint}}
\put(405,424){\usebox{\plotpoint}}
\multiput(405,424)(19.758,6.356){21}{\usebox{\plotpoint}}
\put(806,553){\usebox{\plotpoint}}
\put(806,553){\usebox{\plotpoint}}
\multiput(806,553)(20.608,-2.467){20}{\usebox{\plotpoint}}
\put(1207,505){\usebox{\plotpoint}}
\put(1207,505){\usebox{\plotpoint}}
\multiput(1207,505)(11.006,-17.597){16}{\usebox{\plotpoint}}
\put(1379,230){\usebox{\plotpoint}}
\put(800,280){$t_i^{(n)}$}
\put(430,330){$t_{2i}^{(n+1)}$}
\put(1110,410){$t_{2i+1}^{(n+1)}$}
\end{picture}

Le segment de conique \m{E_{-1}E_0} de \m{\kt_{(E_{-1},E_{\frac{1}{2}},E_0)}}
n'est autre que la courbe \ \m{\Delta=-\frac{\mu(\mu+1)}{2}}. On en d\'eduit
imm\'ediatement le

\bigskip

\begin{xlemm}
Soit \ \m{Z = \bigcup_{(E,F,G)\in{\bf T}}\kt_{(E,F,G)}}. Alors, si
\ \m{(\mu,\Delta)\in Z}, on a \ \m{(\mu,\Delta')\in Z} \ si
$$\frac{\mu(\mu+1)}{2} \ \leq \ \Delta' \ \leq \ \Delta.$$
\end{xlemm}

\vskip 1.2cm

\begin{sub}{\bf Fibr\'es prioritaires g\'en\'eriques}\end{sub}

\begin{subsub}{Cohomologie naturelle}\end{subsub}

\begin{xprop}
Soient $F$ un fibr\'e exceptionnel, $r$, \m{c_1}, \m{c_2} des entiers tels que
\m{r\geq 2}, \m{\mu(F)-x_F<\mu\leq\mu(F)} \ et \ \m{\Delta=\delta(\mu)}.
Alors il existe un fibr\'e vectoriel stable $\ke$ de rang $r$ et de classes de
Chern $c_1$, $c_2$, tel que \ \m{\ext^1(\ke,F)=\lbrace 0\rbrace}.
\end{xprop}

\begin{proof} On consid\`ere la suite \m{(G_n)} de fibr\'es exceptionnels du
\para 3.2.3. Soient $n$ un entier et $\ke$ un faisceau semi-stable de rang $r$
et de classes de Chern \m{c_1}, \m{c_2}. On pose
$$k = \chi(\ke,F), \ \ \ m_n \ = \ -\chi(\ke\ot G_n^*(-3)),$$
qui sont ind\'ependants de $\ke$. Ces entiers sont positifs : pour le premier,
cela d\'ecoule du fait que le point correspondant \`a $\ke$ est situ\'e sous
la conique donnant l'\'equation de \m{\delta(\mu)} sur
\m{\rbrack\mu(F),\mu(F)+x_F\lbrack}. Pour le second on utilise le fait que
\m{H^0(\ke\ot G_n^*(-3))} et \m{H^2(\ke\ot G_n^*(-3))} sont nuls.
On consid\`ere les triades \m{(F,G_{p-1}(3),G_p(3))}. Ceci sugg\`ere de trouver
$\ke$ comme noyau d'un morphisme surjectif ad\'equat
$$\theta :
(F\ot\C^{k})\oplus(G_{p-1}(3)\ot\C^{m_{p+1}})\lra G_p(3)\ot\C^{m_p}.$$
Un tel fibr\'e a en effet les bons rang et classes de Chern, et de plus on
a \ \m{\ext^1(\ke,F)=\lbrace 0\rbrace}. Pour montrer que $\ke$ se d\'eforme en
fibr\'e stable, il suffit qu'il soit prioritaire, car le champ des
faisceaux prioritaires est irr\'eductible (cf. \cite{hi_la}). On prend
\ \m{p=1}, c'est-\`a-dire qu'on consid\`ere des morphismes
$$(F\ot\C^{k})\oplus(G_0(3)\ot\C^{m_2})\lra G_1(3)\ot\C^{m_1}.$$
Alors on a \ \m{\mu(G_1(3))-\mu(G_0(3))\geq 1}, donc \m{\mu(G_1(3))-\mu(F)
> 1}, et la paire \m{(F,G_1(3))} est initiale dans la s\'erie qui la contient.
Ceci entra\^ine que le faisceau des morphismes pr\'ec\'edents est engendr\'e par
ses sections globales. Comme \ \m{r\geq 2}, il existe un morphisme
$$\theta : (F\ot\C^{k})\oplus(G_0(3)\ot\C^{m_2})\lra G_1(3)\ot\C^{m_1}$$
qui est surjectif. Soit
$$\ke \ = \ \ker(\theta).$$
Il reste \`a montrer que $\ke$ est prioritaire, c'est-\`a-dire que
\ \m{\hom(\ke,\ke(-2))=\lbrace 0\rbrace}. On a une suite exacte
$$0\lra\ke\lra (F\ot\C^{k})\oplus(G_0(3)\ot\C^{m_2})\lra G_1(3)\ot\C^{m_1}
\lra 0,$$
d'o\`u on d\'eduit que
$$\hom(\ke,\ke(-2))\ \subset \ (\hom(\ke,F(-2))\ot\C^{k})\oplus
(\hom(\ke,G_0(1))\ot\C^{m_2}).$$
Il faut montrer que
$$\hom((\ke,F(-2))=\hom(\ke,G_0(1))=\lbrace 0\rbrace.$$

Montrons d'abord que \ \m{\hom((\ke,F(-2))=\lbrace 0\rbrace}.
D'apr\`es la suite exacte pr\'ec\'edente, on a une suite exacte
$$(\hom(F,F(-2))\ot\C^{k})\oplus(\hom(G_0(3),F(-2))\ot\C^{m_2})\lra
\hom((\ke,F(-2)) \ \ \ \ \ \ \ \ \ \ \ \ \ \ \ \ $$
$$ \ \ \ \ \ \ \ \ \ \ \ \ \ \ \ \ \ \lra\ext^1(G_1(3),F(-2))\ot\C^{m_1}.$$
On a \ \m{\hom(F,F(-2))=\hom(G_0(3),F(-2))=\lbrace 0\rbrace}, car
\ \m{\mu(G_0(3))>\mu(F)>\mu(F(-2))}.

\noindent D'autre part,
$$\ext^1(G_1(3),F(-2))\ \simeq\ \ext^1(F(-2),G_1)^*$$
par dualit\'e de Serre. Pour montrer que \ \m{\ext^1(F(-2),G_1)=
\lbrace 0\rbrace}, il suffit d'apr\`es \cite{dr1}
de prouver que \ \m{\mu(F(-2))\leq\mu(G_1)}.
Si \ \m{F=\ko} \ c'est \'evident car \ \m{G_1=\ko(-1)}. Sinon, on a
\ \m{\mu(G_1)-\mu(G_0)\geq 2}, et si \ \m{\mu(F(-2))>\mu(G_1)}, on a \
\m{\mu(F)-\mu(G_0)>4}, ce qui est faux car \ \m{\mu(F)-\mu(G_0)<3}.

Montrons maintenant que \ \m{\hom(\ke,G_0(1))=\lbrace 0\rbrace}. On a une
suite exacte
$$(\hom(F,G_0(1))\ot\C^{k})\oplus(\hom(G_0(3),G_0(1))\ot\C^{m_2})
\ \ \ \ \ \ \ \ \ \ \ \ \ \ \ $$
$$\ \ \ \ \ \ \ \ \ \ \ \ \ \ \ \lra
\hom((\ke,G_0(1))\lra\ext^1(G_1(3),G_0(1))\ot\C^{m_1}.$$
On a \ \m{\hom(F,G_0(1))=\lbrace 0\rbrace} \ car \
\m{\mu(F)>\mu(G_1)\geq\mu(G_0(1))}, et
\m{\hom(G_0(3),G_0(1))=\lbrace 0\rbrace}. Il reste \`a prouver que \
\m{\ext^1(G_1(3),G_0(1))=\lbrace 0\rbrace}. On a
$$\ext^1(G_1(3),G_0(1)) \ \simeq \ \ext^1(G_0(1),G_1)^* \ = \
\lbrace 0\rbrace$$
d'apr\`es \cite{dr1} et le fait que \ \m{\mu(G_0(1))\leq\mu(G_1)}. \end{proof}

\vskip 0.8cm

\begin{subsub}{Exemples}\end{subsub}

On montre ici que sous les hypoth\`eses de la proposition 3.6, il existe en
g\'en\'eral des fibr\'es $\ke$  tels que \
\m{\ext^1(\ke,F)\not =\lbrace 0\rbrace}.
Soient \m{n\geq 1} un entier, \m{X} un ensemble de $n$ points
distincts de $\P_2$. On a \ \m{\dim(\ext^1(\ki_X(1),\ko))=n}. Soit
$$0\lra\ko\ot\C^{n}\lra E\lra\ki_X(1)\lra 0$$
une extension d\'efinie par un isomorphisme
\ \m{\ext^1(\ki_X(1),\ko)\simeq\C^{n}}.
Alors $E$ est un fibr\'e vectoriel de rang \m{n+1} et de classes de Chern
$1$ et $n$. On a \ \m{\chi(E,\ko)=\lbrace 0\rbrace}, et si \ \m{n\geq 2},
on a \ \m{0<\mu(E)<x_\ko}. De plus, $E$ est de type de d\'ecomposition
g\'en\'erique \m{(1,0,\ldots,0)}. Les droites de saut de $E$ sont les
droites joignant deux points de $X$.
Le fibr\'e $E$ est stable, \`a cause du r\'esultat
suivant, qui semble connu mais dont je donne une d\'emonstration, n'ayant pu
trouver de r\'ef\'erence :

\bigskip

\begin{xlemm}
Soit $U$ un fibr\'e vectoriel sur $\P_2$, de type de d\'ecomposition
g\'en\'erique
\m{(1,0,\ldots,0)}, ayant un nombre fini de droites de saut, et
tel que
$$h^0(U(-1))=h^0(U^*)=\lbrace 0\rbrace. $$
Alors $U$ est stable.
\end{xlemm}

\begin{proof} Il faut montrer qu'\'etant donn\'e un morphisme g\'en\'eriquement
injectif de fibr\'es vectoriels
$$f : G\lra U,$$
avec \ \m{2\leq rg(G)\leq rg(U)-2}, on a \ \m{c_1(G)\leq 0}. On a
\ \m{c_1(G)\leq 1}, car $U$ est de type de d\'ecomposition g\'en\'erique
\m{(1,0,\ldots,0)}. Supposons que \ \m{c_1(G)=1}. Alors $G$ est de type de
d\'ecomposition g\'en\'erique \m{(1,0,\ldots,0)}. Sur une droite g\'en\'erique
$\ell$ on a un morphisme g\'en\'eriquement injectif \ \m{G_{\mid\ell}\lra
U_{\mid\ell}}, et comme \m{G_{\mid\ell}} et \m{U_{\mid\ell}} sont de type de
d\'ecomposition \m{(1,0,\ldots,0)}, ce morphisme est injectif. Il en d\'ecoule
qu'il n'y a qu'un nombre fini de points de $\P_2$ o\`u $f$ n'est pas
injectif. On en d\'eduit que si \m{U_{\mid\ell}} est de type
\m{(1,0,\ldots,0)}, il en est de m\^eme de \m{G_{\mid\ell}}, et donc
\m{\coker(f)} modulo torsion est de type de d\'ecomposition trivial sur $\ell$
et n'a comme $U$ qu'un nombre fini de droites de saut. Donc \m{\coker(f)}
modulo torsion est trivial, ce qui contredit le fait que \
\m{h^0(U^*)=\lbrace 0\rbrace}. \end{proof}

\bigskip

On a \ \m{\Delta(E^*) = \delta(\mu(E^*))}, et \ \m{\chi(E^*,\ko)=3}. Donc si \
\m{n\geq 2}, le point \m{(\mu(E^*),\Delta(E^*))} est situ\'e sur \m{G(\ko)}.
De plus on a, si \m{n\geq 4},
$$\dim(\hom(E^*,\ko)) \ \geq n \ > \ 3,$$
donc
$$\ext^1(E^*,\ko) \ \not= \ \lbrace 0\rbrace.$$

\vskip 0.8cm

\begin{subsub}{D\'emonstration du th\'eor\`eme 3.1}\end{subsub}

Soient $F$ un fibr\'e exceptionnel, $r$, \m{c_1}, \m{c_2} des entiers tels que
\m{\mu(F)-x_F<\mu<\mu(F)+x_F}, \m{\Delta<\delta(\mu)} \ et \
\m{(\mu,\Delta)\not=(\mu(F),\Delta(F))}. On peut se limiter au cas o\`u
 \m{\mu(F)-x_F<\mu\leq\mu(F)}, l'autre cas s'en d\'eduisant par dualit\'e.
On a alors
$$p \ = \ r.rg(F)(P(\mu-\mu(F))-\Delta-\Delta(F)) \ \ > \ \ 0.$$

On peut supposer que \ \m{\mu > \delta'(\mu)} : dans le cas contraire on
montrera au \para 3.3.4 qu'il existe un \'el\'ement $t$ de {\bf T} tel que
\ \m{(\mu,\Delta)\in\kt_t} \ , et que la partie 1- du th\'eor\`eme 3.1 est
vraie dans ce cas. Montrons qu'on a \ \m{p. rg(F) < r} : ceci \'equivaut \`a
$$\delta(\mu)-\Delta \ < \ \frac{1}{rg(F)^2}.$$
L'in\'egalit\'e est vraie car le terme de droite est la distance des points
not\'es $P$ et $F$ sur la premi\`ere
figure du \para 3.1, et le point \m{(\mu,\Delta)} est
situ\'e \`a l'int\'erieur
du "losange" limit\'e par \m{G(F)}, \m{G'(F)}, \m{D(F)} et \m{D'(F)}.
Il existe donc des entiers \m{r'}, \m{c'_1},
\m{c'_2}, tels que $r$, \m{c_1} et \m{c_2} soient le rang et le classes de
Chern d'une somme directe d'un fibr\'e vectoriel $\ku$ de rang \m{r'} et de
classes de Chern \m{c'_1},\m{c'_2} et de \m{F\ot\C^{p}}. Le point correspondant
\`a $\ku$ est situ\'e sur la conique d'\'equation
$$\Delta = P(\mu-\mu(F))-\Delta(F)$$
(ceci parce que \ \m{\chi(F,\ku)=0}), qui est la conique contenant \m{G(F)}.
On v\'erifie imm\'ediatement qu'on a on a \ \m{\Delta \geq \delta'(\mu)} \ si
et seulement si ce point est situ\'e sur le segment \m{G(F)} de la conique.

Supposons que \ \m{\Delta \geq \delta'(\mu)} \ et \ \m{r'\geq 2}. Dans ce
cas il existe d'apr\'es la proposition 3.6 un fibr\'e stable $\ku$ de
rang \m{r'} et de
classes de Chern \m{c'_1},\m{c'_2} tel que \ \m{\ext^1(\ku,F)=\lbrace 0
\rbrace}. Le fibr\'e
$$\ke \ = \ (F\ot\C^{p})\oplus\ku$$
est prioritaire, de rang \m{r} et de classes de Chern \m{c_1},\m{c_2}. Les
fibr\'es prioritaires g\'en\'eriques sont de ce type, car les fibr\'es tels que
$\ke$ sont d\'efinis par la suite de conditions ouvertes suivante :

\medskip

\noindent (i) on a \ \m{\ext^2(F,\ke)=\lbrace 0\rbrace}.

\noindent (ii) Le morphisme canonique d'\'evaluation
$$ev : F\ot\C^{p}=F\ot\hom(F,\ke)\lra\ke$$
est injectif.

\noindent (iii) Si \ \m{\ku=\coker(ev)}, $\ku$ est un fibr\'e stable tel que
\ \m{\ext^1(\ku,F)=\lbrace 0\rbrace}.

\medskip

Supposons maintenant que \ \m{r'=1}. Dans ce cas on doit avoir \ \m{F=\ko} \ et
\m{c_2=1}. Les faisceaux de \m{M(r',c'_1,c'_2)} sont de la forme \m{\ki_x}
(id\'eal d'un point $x$ de $\P_2$). On a \ \m{\ext^1(\ki_x,\ko)=\C}, d'o\`u
le th\'eor\`eme 3.1 dans ce cas.

Il reste \`a traiter le cas o\`u \ \m{\Delta < \delta'(\mu)}.
C'est une cons\'equence du th\'eor\`eme 3.2,
dont la d\'emonstration suit.

\vskip 1.2cm

\begin{subsub}{D\'emonstration du th\'eor\`eme 3.2}\end{subsub}

Soit \m{(E,F,G)\in{\bf T}}. En consid\'erant la suite spectrale de Beilinson
g\'en\'eralis\'ee associ\'ee \`a \m{(E,F,G)}, on voit imm\'ediatement que les
points \m{(\mu,\Delta)} de \m{\kt_{(E,F,G)}} (\`a coordonn\'ees rationnelles)
sont les paires \m{(\mu(\ke),\Delta(\ke))}, o\`u $\ke$ est de la forme
$$\ke\ = (E\ot\C^{a})\oplus(F\ot\C^{b})\oplus(G\ot\C^{c}),$$
avec \ $a,b,c\geq 0$ \ non tous nuls. Le fibr\'e pr\'ec\'edent est
prioritaire et rigide, c'est donc un fibr\'e prioritaire g\'en\'erique.

On pose comme dans le lemme 3.5,
$$Z \ = \ \bigcup_{(E,F,G)\in{\bf T}}\kt_{(E,F,G)}.$$
La partie 1- du th\'eor\`eme 3.2 est une cons\'equence imm\'ediate du \para
3.2.4. Il reste donc \`a prouver que
$$Z \ = \ \ks.$$
Soit \ \m{(\mu,\Delta)\in Z}. Alors on a \ \m{\Delta\leq\delta'(\mu)}, car
les fibr\'es prioritaires g\'en\'eriques ayant les invariants \m{\mu} et
\m{\Delta} sont rigides, comme on vient de le voir. On a donc \
\m{Z\subset\ks}.

Soit $F$ un fibr\'e exceptionnel tel que \ \m{-1<\mu(F)\leq 0}, \m{(G_n)} la
s\'erie exceptionnelle \`a gauche associ\'ee \`a $F$. On va montrer que lorsque
$n$ tend vers l'infini, le segment de conique \m{G_nF} de
\m{T_{(G_{n-1},G_n,F)}} tend vers le segment de conique
$$\lbrace(\mu,\delta'(\mu)), \mu(F)-x_F<\mu\leq\mu(F)\rbrace.$$
On montrerait de m\^eme que si  \ \m{-1\leq\mu(F)<0}, et si \m{(H_n)} est la
s\'erie exceptionnelle \`a droite associ\'ee \`a $F$, alors lorsque
$n$ tend vers moins l'infini, le segment de conique \m{FH_n} de
\m{T_{(F,H_n,H_{n+1})}} tend vers le segment de conique
$$\lbrace(\mu,\delta'(\mu)), \mu(F)\leq\mu<\mu(F)+x_F\rbrace.$$
D'apr\`es le lemme 3.5, ceci entra\^ine que \ \m{\ks\subset Z}.

L'\'equation du segment de conique \m{G_nF} de \m{T_{(G_{n-1},G_n,F)}} est
$$\Delta\ = \ P(\mu-\mu(G_{n-1})-3)-\Delta(G_{n-1}).$$
On a
$$\lim_{n\rightarrow\infty}(\mu(G_{n-1})) \ = \ \mu(F)-x_F, \ \ \
\lim_{n\rightarrow\infty}(\Delta(G_{n-1})) \ = \ \frac{1}{2}.$$
Donc le segment \m{G_nF} tend vers la courbe
$$\lbrace(\mu,\phi(\mu)), \mu(F)-x_F<\mu\leq\mu(F)\rbrace.$$
avec
$$\phi(\mu) \ = \ P(\mu-\mu(F)+x_F-3)-\frac{1}{2}.$$
On v\'erifie imm\'ediatement que \ \m{\phi(\mu)=\delta'(\mu)}, ce qui ach\`eve
la d\'emonstration du th\'eor\`eme 3.2.

\vskip 1.5cm

\begin{sub}{\bf
Vari\'et\'es de modules fins de faisceaux prioritaires instables}
\end{sub}

Soient $r$, \m{c_1}, \m{c_2} des entiers tels que \m{r\geq 1} et \
\m{M(r,c_1,c_2)=\emptyset}. Soient
$$\mu \ = \ \frac{c_1}{r}, \ \ \ \ \Delta \ = \ \frac{1}{r}(c_2 -
\frac{r-1}{2r}c_1^2), \ \ \ \ \chi = \frac{c_1(c_1+3)}{2}+r-c_2.$$
On suppose que
$$-1 \ < \ \mu \leq \ 0 \ \ \ \ {\rm et} \ \ \ \ \Delta \ \geq
\frac{\mu(\mu+1)}{2}.$$
Il existe donc des fibr\'es prioritaires de rang $r$ et de classes de Chern
\m{c_1}, \m{c_2}.

\vskip 1.2cm

\begin{subsub}{Le cas \ \m{\Delta < \delta'(\mu)}}\end{subsub}

Dans ce cas il existe d'apr\`es les th\'eor\`emes 3.1 et 3.2 une triade
\m{(E,F,G)} de fibr\'es exceptionels de pentes comprises entre \m{-1} et $0$
et des entiers $m$, $n$, $p$ tels que le faisceau prioritaire g\'en\'erique de
rang $r$ et de classes de Chern \m{c_1}, \m{c_2} soit isomorphe \`a
$$\ke \ = \ (E\ot\C^{m})\oplus(F\ot\C^{n})\oplus(G\ot\C^{p}).$$
Le fibr\'e $\ke$ est rigide. On en d\'eduit imm\'ediatement la

\medskip

\begin{xprop}
Tout ensemble ouvert $\kx$ de faisceaux de rang $r$ et de classes de Chern
\m{c_1}, \m{c_2} sur $\P_2$, constitu\'e de classes d'isomorphisme de
faisceaux prioritaires, contient celle de $\ke$, et $\kx$ admet une vari\'et\'e
de modules fins d\'efinie localement si et seulement si $\kx$ est r\'eduit \`a
la classe d'isomorphisme de $\ke$.
\end{xprop}

\vskip 1.2cm

\begin{subsub}{Le cas \ \m{\Delta > \delta'(\mu)}}\end{subsub}

\medskip

\subsubsubsection{3.4.2.1}{Cas d'existence de vari\'et\'es de modules fins}

\medskip

Soit $F$ l'unique fibr\'e exceptionnel tel que \
\m{\mu(F)-x_F<\mu\leq\mu(F)+x_F}. On supposera que \
\m{\mu(F)-x_F<\mu\leq\mu(F)} (l'autre cas est analogue). On supposera aussi
que \ \m{(c_1,c_2)\not = (0,1)}. Soient
$$p \ = \ r.rg(F)(P(\mu-\mu(F))-\Delta-\Delta(F)) \ = \ r.rg(F)(\delta(\mu)
-\Delta)$$
et \m{r'}, \m{c'_1}, \m{c'_2} les entiers tels que le faisceau prioritaire
g\'en\'erique de rang $r$ et de classes de Chern \m{c_1}, \m{c_2} soit de la
forme
$$(F\ot\C^{p})\oplus\ku,$$
$\ku$ \'etant un fibr\'e stable de rang \m{r'}
et de classes de Chern \m{c'_1}, \m{c'_2}. Soient
$$\mu' \ = \ \frac{c'_1}{r'}, \ \ \ \ \Delta' \ = \ \frac{1}{r'}(c'_2 -
\frac{r'-1}{2r'}{c'_1}^2), \ \ \ \ \chi' = \frac{c'_1(c'_1+3)}{2}+r'-c'_2.$$
On a \ \m{\mu(F)-x_F<\mu'\leq\mu(F)} \ et \ \m{\Delta'=\delta(\mu')}.

Supposons que \m{r'}, \m{c'_1} et \m{\chi'} soient premiers entre eux. Dans ce
cas \m{M(r',c'_1,c'_2)} est lisse et il existe un {\em faisceau universel}
$\ke$ sur \ \m{M(r',c'_1,c'_2)\times\P_2}. Soit \
\m{M_0\subset M(r',c'_1,c'_2)} \ l'ouvert constitu\'e des points $x$ tels que
\ \m{\ext^1(\ke_x,F)=\lbrace 0\rbrace}. Soit \m{\kx_0} l'ensemble des classes
d'isomorphisme des faisceaux \ \m{(F\ot\C^{p})\oplus\ke_x} , $x$ parcourant
\m{M_0}. C'est un ensemble ouvert, et il ais\'e de voir que \
\m{(M_0, (F\ot\C^{p})\oplus\ke_{M_0})} \ est une vari\'et\'e de modules fins
pour \m{\kx_0}. On a une sorte de r\'eciproque :

\bigskip

\begin{xprop}
Soit $\kx$ un ensemble ouvert de classes d'isomorphisme de faisceaux sur
$\P_2$ contenant des classes d'isomorphisme de
faisceaux prioritaires de rang $r$ et de classes de
Chern \m{c_1}, \m{c_2}. Alors s'il existe une vari\'et\'e de modules fins
d\'efinie localement pour $\kx$, les entiers $r'$, \m{c'_1} et \m{\chi'} sont
premiers entre eux.
\end{xprop}

\begin{proof} Quitte \`a remplacer $\kx$ par un ensemble plus petit, on peut
supposer que $\kx$ est constitu\'e uniquement de classes d'isomorphisme de
faisceaux prioritaires de la forme \ \m{(F\ot\C^{p})\oplus\ku}, o\`u $\ku$ est
un fibr\'e stable de rang $r'$ et de classes de Chern \m{c'_1}, \m{c'_2}. On
peut aussi supposer que $\kx$ admet une vari\'et\'e de modules fins
\m{(M,\ke_0)}. Alors $M$ est lisse car pour tout faisceau prioritaire $E$ on a
\ \m{\ext^2(E,E)=\lbrace 0\rbrace}. Il en d\'ecoule que le faisceau coh\'erent
\ \m{p_{M*}(F^*\ot\ke_0)} \ sur $M$ est localement libre de rang $p$, et le
morphisme canonique de fibr\'es vectoriels sur \ \m{M\ot\P_2}
$$\Phi : F\ot p_{M*}(F^*\ot\ke_0)\lra\ke_0$$
est injectif. Son conoyau est une famille de fibr\'es stables de rang \m{r'} et
de classes de Chern \m{c'_1}, \m{c'_2} param\'etr\'ee par $M$. Le morphisme
induit par \m{\coker(\Phi)}
$$f : M\lra M(r',c'_1,c'_2)$$
est une immersion ouverte, comme on peut le voir en examinant le morphisme
tangent \m{Tf}. Il en d\'ecoule que $M$, vu comme ouvert de
\m{M(r',c'_1,c'_2)}, poss\`ede un fibr\'e universel. D'apr\`es \cite{dr3}, les
entiers \m{r'}, \m{c'_1} et \m{\chi'} sont premiers entre eux. \end{proof}

\vskip 0.8cm

\subsubsubsection{3.4.2.2}{Application du th\'eor\`eme 2.13}

\medskip

On suppose ici que \m{r'}, \m{c'_1} et \m{\chi'} sont premiers entre eux.
On note $\ke$ un faisceau universel sur \ \m{M(r',c'_1,c'_2)\times\P_2}.
Soit $\ky$ l'ensemble des classes d'isomorphisme de faisceaux coh\'erents $\ku$
qui peuvent s'\'ecrire comme extensions
$$0\lra F\ot\C^{p}\lra\ku\lra E\lra 0,$$
$E$ \'etant un faisceau stable de rang \m{r'} et de classes de Chern \m{c'_1},
\m{c'_2}, de telle sorte que l'application lin\'eaire induite
$${\C^{p}}^*\lra\ext^1(E,F)$$
soit surjective. De tels faisceaux sont prioritaires, de rang $r$ et de classes
de Chern \m{c_1}, \m{c_2}. On d\'eduit du th\'eor\`eme 2.13 le

\medskip

\begin{xtheo}
Soit $M$ l'ouvert de \m{M(r',c'_1,c'_2)} constitu\'e des points $x$ tels que \

\noindent
\m{\dim(\ext^1(\ke_x,F))\leq p} . Alors $\ky$ est un ensemble ouvert et il
existe une vari\'et\'e de modules fins d\'efinie localement \m{(M,(\ke_i))}
pour $\ky$.
\end{xtheo}

\noindent{\bf Remarques : } {\bf 1 -} Si $p$ est suffisamment grand, on a \
\m{M \ = \ M(r',c'_1,c'_2)} , et $M$ est projective. C'est donc une vari\'et\'e
de modules fins d\'efinie localement {\em maximale} (cf. \para 2.6). J'ignore
si en g\'en\'eral $M$ est maximale.

\medskip

\noindent {\bf 2 -} Il existe bien s\^ur d'autres faisceaux prioritaires de
rang $r$ et de classes de Chern \m{c_1}, \m{c_2} que ceux de $M$. Mais j'ignore
s'il existe une vari\'et\'e de modules fins d\'efinie localement ne contenant
que des faisceaux prioritaires de rang $r$ et de classes de Chern \m{c_1},
\m{c_2}, mais contenant au moins un faisceau qui n'est pas dans $M$.

\medskip

\noindent {\bf 3 -} La dimension du groupe d'automorphismes des faisceaux de
$M$ est
$$p^2+1+3pr'.rg(F)(\mu(F)-\mu').$$

\newpage

\subsubsubsection{3.4.2.3}{Le probl\`eme de l'existence d'un faisceau universel
global}

\medskip

On consid\`ere la situation du th\'eor\`eme 3.10. Supposons qu'il existe une
vari\'et\'e de modules fins \m{(M,\kf)} pour $\ky$ (et donc que les \m{\ke_i} se
recollent). Comme dans la d\'emonstration de la proposition 3.9 on trouve une
suite exacte de faisceaux sur \ \m{M\times\P_2}
$$0\lra F\ot p_M^*(\gamma)\lra\kf\lra\ke_M\ot p_M^*(L)\lra 0,$$
$\gamma$ \'etant un fibr\'e vectoriel de rang $p$ sur $M$ et $L$ un fibr\'e en
droites sur $M$. Cette suite exacte est associ\'ee \`a
$$\sigma \ \in \ \ext^1(\ke_M\ot p_M^*(L), F\ot p_M^*(\gamma)).$$
Pour tout \ \m{x\in M}, soient
$$\theta_x : \ext^1(\ke_M\ot p_M^*(L), F\ot p_M^*(\gamma))\lra
\ext^1(\ke_x\ot L_x,F\ot\gamma_x)\simeq L({\C^{p}}^*,\ext^1(\ke_x,F))$$
le morphisme canonique et \ \m{\sigma_x=\theta_x(\sigma)}. Alors \m{\sigma_x}
est surjective.

R\'eciproquement, soient $\gamma$ un fibr\'e vectoriel de rang $p$ sur $M$, $L$
un fibr\'e en droites sur $M$ et \ \m{\sigma\in\ext^1(\ke_M\ot p_M^*(L),
F\ot p_M^*(\gamma))} \ tels que pour tout $x$ dans $M$, \m{\sigma_x} soit
surjective. Alors l'extension de \m{\ke_M\ot p_M^*(L)} par \m{F\ot
p_M^*(\gamma)} associ\'ee \`a $\sigma$ fournit un faisceau universel pour
$\ky$, qui admet donc une vari\'et\'e de modules fins d\'efinie globalement.

D'apr\`es la d\'emonstration du th\'eor\`eme 2.13 il existe un faisceau
coh\'erent $\A$ sur $M$ tel que pour tout \m{x\in M(r',c'_1,c'_2)} on ait un
isomorphisme canonique
$$\A_x \ \simeq \ \ext^1(\ke_x,F),$$
et que pour toute sous-vari\'et\'e localement ferm\'ee $Z$ et tout \
\m{s\in H^0(Z,\C^{p}\ot\A)} \ il existe une suite exacte de faisceaux
coh\'erents sur \ \m{Z\times\P_2}
$$0\lra F\ot\C^{p}\lra\ku\lra\ke_Z\lra 0$$
telle qu'en tout point \m{x\in Z} l'extension
$$0\lra F\lra\ku_x\lra\ke_x\lra 0$$
soit associ\'ee \`a \ \m{s(x)\in\A_x=\ext^1(\ke_x,F)}. Lorsqu'on remplace
\m{F\ot\C^{p}} par \m{F\ot\gamma} et $\ke$ par \m{\ke\ot p_m^*(L)} on a un
faisceau \m{\A(\gamma,L)} analogue \`a \m{\A\ot\C^{p}} tel qu'en tout
\m{x\in M} on ait un isomorphisme canonique
$$\A(\gamma,L)_x \ \simeq \ \ext^1(\ke_x\ot L_x,F\ot\gamma_x).$$
On a
$$\A(\gamma,L) \ \simeq \ \A(\gamma,\ko)\ot L^*.$$

Soit $n$ la dimension maximale de \m{\dim(\ext^1(\ke_x,F))}, $x$ parcourant
\m{M(r',c'_1,c'_2)}, et pour tout entier $i$ tel que \m{1\leq i\leq n}, \m{d_i}
la dimension de la sous-vari\'et\'e localement ferm\'ee de \m{M(r',c'_1,c'_2)}
constitu\'ee des points $x$ tels que \ \m{\dim(\ext^1(\ke_x,F))=i}. Soit
\m{\ko_M(1)} un fibr\'e en droites tr\`es ample sur $M$. On prend $\gamma$
trivial et \
\m{L=\ko_M(-m)}, avec $m$ assez grand pour que le faisceau \
\m{\A(\ko\ot\C^{p},L)=\A(\ko\ot\C^{p},\ko)\ot L^*} \ soit
engendr\'e par ses sections globales. Cela implique que pour tout \m{x\in M},
\m{\theta_x} est surjective. On cherche $p$ sections de ce faisceau
engendrant \m{\ext^1(\ke_x,F)} pour tout \m{x\in M}. Un calcul simple permet
de d\'emontrer la

\medskip

\begin{xprop}
Si \ \m{p\geq{\rm Max}(d_i+i, 1\leq i\leq n)}, on a \
\m{M=M(r',c'_1,c'_2)} \ et il existe une vari\'et\'e de
modules fins globale \m{(M,\kf)} pour $\ky$.
\end{xprop}

\vskip 0.8cm

\subsubsubsection{3.4.2.4}{Exemple : le cas \ \m{(r',c'_1,c'_2)=(4,-2,4)}}

\medskip

Dans ce cas $F$ est le fibr\'e vectoriel habituellement not\'e $Q^*$, qui est
le noyau du morphisme d\'evaluation \ \m{\ko\ot H^0(\ko(1))\lra\ko(1)}. On a
ici \ \m{\chi'=1}, donc $r'$, \m{c'_1} et \m{\chi'} sont bien premiers entre
eux. Les faisceaux des vari\'et\'es de modules fins \'etudi\'ees ici seront des
extensions
$$0\lra Q^*\ot\C^{p}\lra \ku\lra E\lra 0,$$
$E$ \'etant dans \m{M(4,-2,4)}.
On a alors
$$(r,c_1,c_2)=(2p+4, -p-2,\frac{p(p+5)}{2}+4), \ \ \ \ \mu=-\frac{1}{2},
\ \ \ \ \Delta=\frac{3}{8}+\frac{1}{2(p+2)}.$$
On montre dans la proposition 3.12 ci-dessous qu'avec les notations du \para
3.4.2.3, on a \ \m{n=2}, \m{d_1=4} \ et \ \m{d_2=2}.
On a un isomorphisme canonique
$$M(r',c'_1,c'_2) \ \simeq \ \P_{5}=\P({S^2V}).$$
Si $D$ est une droite de \m{S^2V}, le faisceau semi-stable correspondant \`a
$D$ est le conoyau \m{\ke_D} du morphisme canonique injectif
$$\ko(-3)\lra\ko(-1)\ot(S^2V/D).$$
On a donc un faisceau universel $\ke$ sur \m{\P({S^2V})} et une suite exacte
de faisceaux sur \ \m{\P({S^2V})\times\P_2}
$$0\lra p_1^*(\ko(-3))\lra p_1^*(\ko(-1))\ot p_2^*(Q)\lra\ke\lra 0,$$
\m{p_1}, \m{p_2} d\'esignant les projections \ \m{\P({S^2V})\times\P_2
\lra\P({S^2V})} \ et \ \m{\P({S^2V})\times\P_2\lra\P_2} \
respectivement. Soit \ \m{Z\subset\P({S^2V})} \ la sous-vari\'et\'e
localement ferm\'ee correspondant aux \'el\'ements d\'ecomposables, $P$
la sous-vari\'et\'e ferm\'ee isomorphe \`a \m{\P_2} compos\'ee des points
\m{\C.u^2}, avec \m{u\in V\backslash\lbrace 0\rbrace}.

\medskip

\begin{xprop}
Si \ \m{p=1} \ on a \ \m{M=\P(S^2V)\backslash P}. Si \ \m{p\geq 2}
\ on a \ \m{M=\P({S^2V})}. Il existe une vari\'et\'e de modules fins
globale \m{(M,\kf)} pour $\ky$ si et seulement si \ \m{p\not = 2}.
\end{xprop}

\begin{proof} On note \m{M(V)} le noyau du morphisme canonique
$$S^2V\ot V^*\lra V$$
et \m{\phi_D} l'application canonique
$$M(V)\lra (S^2V/D)\ot V^*.$$
Alors de la suite exacte
$$0\lra\ko(-3)\lra\ko(-1)\ot(S^2V/D)\lra\ke_D\lra 0$$
on d\'eduit un isomorphisme canonique
$$\ext^1(\ke_D,Q^*) \ \simeq \ \coker({}^t\phi_D)\simeq\ker(\phi_D)^*.$$
Il en d\'ecoule que 
$$\dim(\ext^1(\ke_D,Q^*)) \ = \ 0 \ \ \ \ \ \ {\rm si} \ \ D\in
\P({S^2V})\backslash (Z\cup P),$$
$$\dim(\ext^1(\ke_D,Q^*)) \ = \ 1 \ \ \ \ \ \ {\rm si} \ \ D\in Z,$$
$$\dim(\ext^1(\ke_D,Q^*)) \ = \ 2 \ \ \ \ \ \ {\rm si} \ \ D\in P.$$
On en d\'eduit d\'ej\`a que \ \m{M=\P({S^2V})\backslash P} \ si \m{p=1},
\ \m{M=\P({S^2V})} si \ \m{p\geq 2}, et qu'il existe une vari\'et\'e de
modules fins \m{(M,\kf)} pour $\ky$ si \ \m{p\geq 5} (d'apr\`es la proposition
3.11).

Etant donn\'e que la restriction
$$Pic(\P(S^2V))\lra Pic(Z)$$
est surjective, il d\'ecoule de la discussion du \para 3.4.2.3 qu'on peut
trouver une vari\'et\'e de modules fins \m{(M,\kf)} pour $\ky$ si \ \m{p=1}.
Plus pr\'ecis\'ement, le fibr\'e en droites \ \m{\A} \ sur $Z$ peut \^etre
rendu trivial si on remplace $\ke$ par \m{\ke\ot p_M^*(\L)}, pour un $\L$
convenable dans \m{Pic(M)}.

Il reste \`a traiter les cas \ \m{2\leq p\leq 4}. On a une suite exacte de
faisceaux sur \m{M=\P(S^2V)}
$$0\lra Q^*\ot V\lra \ko\ot M(V)^*\lra \A\lra 0,$$
d'o\`u on d\'eduit ais\'ement que
$$\A_{\mid Z} \ \simeq \ \ko_M(1)_{\mid Z}, \ \ \ \ \A_{\mid P} \ \simeq \
Q_P$$
o\`u \m{Q_P} est le fibr\'e $Q$ sur \m{\P_2}. Puisque \m{Q_P} n'est
pas la restriction \`a $P$ d'un fibr\'e vectoriel sur \m{\P(S^2V)} (son
d\'eterminant n'est pas la restriction \`a $P$ d'un \'el\'ement de
\m{Pic(\P(S^2V))}), il n'est pas possible
de choisir $L$ et $\gamma$ (cf. \para 3.4.2.3) de telle sorte que
\m{\kg_{\mid P}} devienne
engendr\'e par deux sections globales (c'est-\`a-dire trivial). On n'a donc
pas de vari\'et\'e de modules fins globale dans ce cas.

Supposons que \ \m{p=3} ou $4$. Alors trois sections globales de \m{\ko_Z(1)}
suffisent \`a engendrer ce fibr\'e, et trois sections globales de \m{Q_P}
suffisent aussi \`a l'engendrer. Il en d\'ecoule que trois sections
g\'en\'eriques de $\A$ engendrent ce faisceau. Donc dans ce cas il existe une
vari\'et\'e de modules fins pour $\ky$.
\end{proof}

\vskip 0.8cm

\subsubsubsection{3.4.2.5}{Les cas o\`u \ \m{c_1=0} \ et \ \m{c_2=1}}

\medskip

On a alors \ \m{(r',c'_1,c'_2)=(1,0,1)}, et les faisceaux de
\m{M(r',c'_1,c'_2)} sont les id\'eaux de points de \m{\P_2}. On a donc
$$M(r',c'_1,c'_2) \ \simeq \ {\rm Hilb}^1(\P_2) \ \simeq \P_2.$$
Pour tout point $x$ de \m{\P_2} on a une suite exacte canonique
$$0\lra\ko(-2)\ot\Lambda^2x^\bot\lra\ko(-1)\ot x^\bot\lra\ki_x\lra 0,$$
o\`u \ \m{x^\bot\subset V^*} \ est l'orthogonal de $x$. On en d\'eduit la
suite exacte
$$0\lra\C\lra V^*\ot Q_x\lra S^2V^*\ot\Lambda^2Q_x\lra\ext^1(\ki_x,\ko)
\lra 0.$$
Soient \m{p_M}, \m{p_2} les projections \ \m{{\rm Hilb}^1(\P_2)\times
\P_2\lra{\rm Hilb}^1(\P_2)} \ et \m{{\rm Hilb}^1(\P_2)\times\P_2
\lra\P_2} \ respectivement. On d\'eduit de ce qui pr\'ec\`ede un faisceau
universel \m{\ki'} sur \ \m{{\rm Hilb}^1(\P_2)\times\P_2} \ et une suite
exacte de fibr\'es vectoriels sur \m{{\rm Hilb}^1(\P_2)\times\P_2}
$$0\lra p_2^*(\ko(-2))\ot p_M^*(\Lambda^2Q^*)\lra p_2^*(\ko(-1))\ot p_M^*(Q^*)
\lra\ki'\lra 0,$$
et une suite exacte de fibr\'es vectoriels sur \m{{\rm Hilb}^1(\P_2)}
$$0\lra\ko\lra p_M^*(Q)\ot V^*\lra p_M^*(\Lambda^2Q)\ot S^2V^*\lra\A\lra 0,$$
o\`u \m{\A} est le faisceau de la d\'emonstration du th\'eor\`eme 2.13 (qui est
ici localement libre de rang 1) tel qu'en tout point $x$ de \m{\P_2} on
ait \ \m{\A_x\simeq\ext^1(\ki_x,\ko)}. On a
$$\A \ \simeq \ko(3).$$
Soit
$$\ki \ = \ \A\ot p_M^*(\ko(3)).$$
Pour ce nouveau faisceau universel, le fibr\'e $\A$ devient trivial, et la
section $1$ d\'efinit une extension
$$0\lra\ko\lra\kv\lra\ki\lra 0$$
telle qu'en tout point $x$ de \ \m{\P_2={\rm Hilb}^1(\P_2)}, sa
restriction \`a \ \m{\lbrace x\rbrace\times\P_2} \ soit non triviale.

Soit $\ky$ l'ensemble des classes d'isomorphisme des faisceaux \m{\kv_x}, $x$
parcourant \m{\P_2}. Alors $\ky$ est ouvert et \m{(\P_2,\kv)} est une
vari\'et\'e de modules fins pour $\ky$.

Plus g\'en\'eralement, soit $\ky_p$ l'ensemble des classes d'isomorphisme des
faisceaux  \m{\kv_x\oplus(\ko\ot\C^{p})}, $x$ parcourant \m{\P_2}. Alors
$\ky_p$ est ouvert et \m{(\P_2,\kv\oplus(\ko\ot\C^{p}))} est une vari\'et\'e de
modules fins pour $\ky_p$.

\sepsec

\section{Vari\'et\'es de modules fins de faisceaux de rang 1}

\begin{sub}{\bf Faisceaux simples de rang 1}\end{sub}

On suppose dans toute cette section que \ \m{X=\P_2} \ ou \m{X=\P_1\times\P_1}.
On note \m{\omega_X} le fibr\'e canonique sur $X$. 

\medskip

\begin{xlemm}
Soit $E$ un faisceau de rang non nul simple sur $X$.
Alors on a
$$\ext^2(E,E)=\lbrace 0\rbrace.$$
\end{xlemm}

\begin{proof} On a par dualit\'e de Serre
$$\ext^2(E,E) \ \simeq \ \hom(E,E\ot\omega_X)^*.$$
Supposons que \ \m{\ext^2(E,E)\not=\lbrace 0\rbrace}. On en d\'eduit un
morphisme non nul
$$f : E\lra E\ot\omega_X.$$
Soit \ \m{\phi:\omega_X\lra\ko_X} \ un morphisme non nul. Il s'annule donc sur
une hypersurface de $X$.
Alors \m{(I_E\ot\phi)\circ f} est un endomorphisme de $E$ qui n'est pas une
homoth\'etie. On a donc \ \m{\ext^2(E,E)=\lbrace 0\rbrace}. \end{proof}

\medskip

On s'int\'eresse aux vari\'et\'es de modules fins constitu\'ees de faisceaux
simples de rang 1 et de d\'eterminant trivial. D'apr\`es le lemme
pr\'ec\'edent, de telles vari\'etes sont lisses.

\medskip

\begin{xtheo}
On suppose que \ \m{X=\P_2} \ ou \ \m{X=\P_1\times\P_1}.
Soient $Y$ une vari\'et\'e alg\'ebrique int\`egre, $\ke$ une famille plate
compl\`ete de faisceaux de rang $1$ simples sur $X$ et de d\'etermiant trivial.
Soit $U$ l'ouvert de $Y$ constitu\'e des
points $y$ tels que \m{\ke_y} soit sans torsion, et  \
\m{Z=Y\backslash U}.
Alors si \ \m{X=\P_2}, $Z$ est de codimension
au moins $4$ dans $Y$, et si \ \m{X=\P_1\times\P_1}, $Z$ est de
codimension au moins $3$ dans $Y$.
\end{xtheo}

\medskip

Il en d\'ecoule que dans toute vari\'et\'e de modules fins de faisceaux simples
de rang 1 et de d\'eterminant trivial
sur $X$, l'ouvert correspondant aux faisceaux qui sont des
id\'eaux de sous-sch\'emas de $X$ de dimension 0 est non vide et dense.

Le reste du \para 4.1 est consacr\'e \`a la d\'emonstration du
th\'eor\`eme 4.2. Soit $m$ un entier tel
que pour tout point $y$ de $Y$, $\ke_y(m)$ soit engendr\'e par ses sections
globales, et que \ \m{h^1(\ke_y(m)) = h^2(\ke_y(m)) = \lbrace 0\rbrace}.
Soient \ \m{p=h^0(\ke_y(m))}, qui est ind\'ependant de $y$, et $P$ le
polyn\^ome de Hilbert des faisceaux \m{\ke_y}. Soient enfin \
$${\bf Q} \ = \ \hilb^P(\ko_X(-m)\ot\C^{p}),$$
et
$$\pi : \ko_X(-m)\ot\C^{p}\lra\kf$$
le morphisme surjectif universel sur \ \m{{\bf Q}\ot X}. On se ram\`ene
ais\'ement, $\ke$ \'etant compl\`ete, au cas o\`u $Y$ est l'ouvert de $Q$
constitu\'e des points $q$ tels que \m{\kf_q} soit simple, et o\`u \
\m{\ke=\kf_Y}. Soit \ \m{\N=\ker(\pi)}.

Commen\c cons par donner une description des faisceaux simples de rang 1 sur
$X$ et une caract\'erisation de ceux qui sont sans torsion.

\medskip

\begin{xlemm}
Soit $E$ un faisceau simple de rang $1$ sur $X$ et de d\'eterminant trivial.
Alors on a \ \m{\ext^2(E,\ko_X)=\lbrace 0\rbrace}, et $E$ est sans torsion si
et seulement si \ \m{\dim(\hom(E,\ko_X))=1}.
\end{xlemm}

\begin{proof} On a une suite exacte
$$0\lra T\lra E\lra \ki_S\ot L\lra 0,$$
$T$ d\'esignant le sous-faisceau de torsion de $E$, $S$ un sous-sch\'ema de
dimension 0 de $X$, \m{\ki_S} le faisceau d'id\'eaux de $S$ et $L$ un
fibr\'e en droites sur $X$ (pour un \'etude pr\'ecise des faisceaux de
torsion, voir \cite{lp2}, \cite{lp4} ou \cite{si}). Supposons $T$ non nul.
Alors $T$ est pur de dimension 1, car
s'il avait un sous-faisceau de dimension 0, il existerait des morphismes non
nuls \ \m{\ki_S\ot L\lra T}, et $E$ ne serait pas simple. Soit $D$ le support
sch\'ematique de $T$. On a alors
$$c_1(T) \ = \ \ko(D), \ \ \ \ L \ = \ \ko(-D).$$
On a donc
$$\hom(E,\ko_X) \ \simeq \ \hom(\ki_S\ot L,\ko_X) \ \simeq \ H^0(\ko(D)).$$
et \ \m{\dim(\hom(E,\ko_X))>1}. Si maintenant $T$ est nul, on a \
\m{E\simeq \ki_S}, et donc

\noindent \m{\dim(\hom(E,\ko_X))=1}.

Il reste \`a montrer que \ \m{\ext^2(E,\ko_X)=\lbrace 0\rbrace}. Par dualit\'e
de Serre, on a
$$\ext^2(E,\ko_X) \ \simeq \ \hom(\omega_X^{-1},E)^*.$$
Supposons qu'il existe un \ \m{\phi\in\hom(\omega_X^{-1},E)} \ non nul.
Puisque \ \m{\ki_S\ot L\subset\ko_X}, il existe des morphismes
\ \m{E\lra\omega_X^{-1}} \ dont la restriction au support de $T$ est non
nulle. On en d\'eduit par composition avec $\phi$ un endomorphisme de $E$ qui
n'est pas une homoth\'etie. Donc \ \m{\ext^2(E,\ko_X)=\lbrace 0\rbrace}.
\end{proof}

\medskip

Soit \m{q\in Y}. On a une suite exacte
$$0\lra\hom(\kf_q,\ko_X)\lra H^0(\ko_X(m))\ot\C^{p*}\lra\hom(\N_q,\ko_X)\lra
\ext^1(\kf_q,\ko_X)\lra 0.$$
Puisque \ \m{\ext^2(\kf_q,\ko_X)=\lbrace 0\rbrace}, on en d\'eduit que la
dimension de \m{\hom(\N_q,\ko_X)} ne d\'epend pas de $q$. On obtient donc un
morphisme de fibr\'es vectoriels sur $Y$
$$\Phi : F_0=\ko_Y\ot H^0(\ko_X(m))\ot\C^{p*}\lra
F_1=p_{Y*}(\underline{\hom}(\N,\ko_{Y\times X}))$$
tel que le lieu \m{Y_i}  des points $q$ de $Y$ tels que \
\m{\dim(\hom(\kf_q,\ko_X))=i} \ est pr\'ecis\'ement celui des points o\`u le
noyau de $\Phi$ est de dimension $i$. Soient \ \m{G=F_0^*\ot F_1}, vu comme
vari\'et\'e alg\'ebrique, et \m{G_i} la sous-vari\'et\'e localement ferm\'ee
constitu\'ee des morphismes dont le noyau est de dimension $i$. Alors $G$ et
les \m{G_i} sont lisses, et en tout point $f$ de \m{G_i}, l'espace normal
de \m{G_i} en $f$ est canoniquement isomorphe \`a \m{\hom(\ker(f),\coker(f))}.

Soit $q$ un point de $Y$ tel que $\kf_q$ ait de la torsion, et \
\m{i=\dim(\hom(\kf_q,\ko_X))\geq 2}. Alors \m{q\in G_i}. Le morphisme $\Phi$
peut \^etre vu comme un morphisme
$$Y\lra G.$$
L'application
$$TY_q\lra \hom(\hom(\kf_q,\ko_X),\ext^1(\kf_q,\ko_X))$$
d\'eduite de \m{T\phi} est nulle sur les \m{PGL(p)}-orbites, donc se factorise
de la fa\c con suivante :
$$TY_q \ \hfl{\alpha}{} \ \ext^1(\kf_q,\kf_q) \ \hfl{\beta}{} \
\hom(\hom(\kf_q,\ko_X),\ext^1(\kf_q,\ko_X)),$$
o\`u $\alpha$ est le morphisme de d\'eformation infinit\'esimale de
Koda\"\i ra-Spencer de $\kf$ en $q$ et $\beta$ le morphisme canonique (pour des
d\'emonstrations analogues, voir par exemple \cite{hi}). Le th\'eor\`eme 4.2
est donc une cons\'equence du

\medskip

\begin{xlemm}
Soit $E$ un faisceau simple de rang $1$ sur $X$ de d\'eterminant trivial et
ayant de la torsion. Alors l'application canonique
$$\beta : \ext^1(E,E)\lra\hom(\hom(E,\ko_X),\ext^1(E,\ko_X))$$
est de rang au moins $4$ si \ \m{X=\P_2}, et de rang au moins $3$ si \
\m{X=\P_1\times\P_1}.
\end{xlemm}

\begin{proof} On consid\`ere la suite exacte dej\`a vue dans le lemme 4.3
$$0\lra T\lra E\lra \ki_S\ot L\lra 0,$$
$T$ d\'esignant le sous-faisceau de torsion de $E$ (qui est pur de dimension
1), $S$ un sous-sch\'ema de dimension 0 de $X$, \m{\ki_S} le faisceau
d'id\'eaux de $S$ et $L$ un fibr\'e en droites sur $X$. On a \ \m{L=\ko(-D)},
$D$ d\'esignant le support sch\'ematique de $T$. On a alors
$$c_1(T) \ = \ -c_1(L), \ \ \ \ c_2(T) \ = c_1(L)^2-s+c_2,$$
$s$ d\'esignant la longueur de $S$, et \ \m{c_2=c_2(E)}.

\vskip 0.8cm

\noindent {\bf Etape 1 : }{\em quelques calculs}

\medskip

Le th\'eor\`eme de Riemann-Roch fournit les r\'esultats suivants, dont on se
servira par la suite :
$$\chi(T,\ki_S\ot L) \ = \ \frac{1}{2}c_1(L)^2-\frac{\omega_X c_1(L)}{2}-c_2+s,
\ \ \ \ \
\chi(\ki_S\ot L,T) \ = \ \frac{1}{2}c_1(L)^2+\frac{\omega_X c_1(L)}{2}-c_2+s,$$
$$\chi(T,\ko_X) \ = \ -\frac{1}{2}c_1(L)^2-\frac{\omega_X c_1(L)}{2}-c_2+s,
\ \ \ \ \chi(E,\ko_X) \ = \ 1 - c_2, \ \ \ \ \chi(E,\ko_S) \ = \ s.$$
On en d\'eduit
$$\dim(\hom(E,\ko_X)) \ = \ h^0(L^{-1}), \ \ \ \ \
\dim(\ext^1(E,\ko_X)) \ = \ h^0(L^{-1})+c_2-1,$$
(on a \ \m{\ext^2(E,\ko_X)=\lbrace 0\rbrace} d'apr\`es le lemme 4.3).

Puisque \ \m{\hom(\ki_S\ot L,T)=\lbrace 0\rbrace} (car $E$ est simple) et
\ \m{\ext^2(\ki_S\ot L,T)=\lbrace 0\rbrace} (par dualit\'e de Serre) on a
$$\dim(\ext^1(\ki_S\ot L,T)) \ = \ -\chi(\ki_S\ot L,T),$$
et c'est la dimension g\'en\'erique. Il en d\'ecoule qu'en bougeant $S$ il
est suffisant de d\'emontrer le lemme 4.4 dans le cas o\`u $S$ ne rencontre
pas le support de $T$ (autrement dit, on peut toujours d\'eformer $E$ en
faisceau simple tel que $S$ ne rencontre pas le support de $T$). On a
$$\ext^2(E,\ko_S) \ = \ \lbrace 0\rbrace.$$
En effet, par dualit\'e de Serre, on a
$$\ext^2(E,\ko_S) \ \simeq \ \hom(\ko_S,E\ot\omega_X)$$
qui est nul car $T$ est pur de dimension 1. On en d\'eduit que
$$\dim(\ext^1(E,\ko_S)) \ = \ s .$$
D'autre part il d\'ecoule de la d\'emonstration du lemme 3.2 de \cite{lp4} que
$T$ peut se d\'eformer en faisceau lisse (un fibr\'e en droites sur une courbe
lisse). On peut donc supposer que \
\m{\ext^2(T,T)=\lbrace 0\rbrace}. D'apr\`es la suite exacte
$$0\lra T\lra E\lra\ki_S\ot L\lra 0$$
on en d\'eduit que
$$\ext^2(E,T) \ = \ \lbrace 0\rbrace.$$

\vskip 0.8cm

\noindent{\bf Etape 2 : }{\em interpr\'etation de \m{\coker(\beta)}}

\medskip

On note \m{V(L)} le fibr\'e vectoriel sur $X$ conoyau du morphisme canonique
injectif
$$L\lra \ko_X\ot H^0(L^{-1})^*,$$
et \m{V_S} celui du morphisme canonique injectif
$$\ki_S\ot L\lra \ko_X\ot H^0(L^{-1})^*.$$
On a une suite exacte
$$0\lra L\ot\ko_S\lra V_S\lra V(L)\lra 0.$$
Puisque \ \m{\ext^2(E,T) = \lbrace 0\rbrace}, le morphisme canonique
$$f : \ext^1(E,E)\lra\ext^1(E,\ki_S\ot L)$$
est surjectif. On consid\`ere maintenant la suite exacte
$$0\lra \ki_S\ot L\lra \ko_X\ot H^0(L^{-1})^*\lra V_S\lra 0.$$
On en d\'eduit la suite exacte

$$\ext^1(E,\ki_S\ot L) \ \hfl{g}{} \
\hom(\hom(E,\ko_X),\ext^1(E,\ko_X))\lra\ext^1(E,V_S)\lra$$
$$\ \ \lra\ext^2(E,\ki_S\ot L).$$
On a \ \m{\ext^2(E,\ki_S\ot L)=\lbrace 0\rbrace}. En effet, par dualit\'e de
Serre, il suffit de v\'erifier que \
\m{\hom(\ki_S\ot L,E\ot\omega_X)=\lbrace 0\rbrace}. Cela d\'ecoule du fait que
\m{\hom(\ki_S\ot L,E)=\lbrace 0\rbrace} (car $E$ est simple), et de ce qu'il
existe des morphismes non nuls \ \m{\omega_X\lra\ko_X}. Comme \
\m{\beta=g\circ f}, on en d\'eduit un isomorphisme
$$\coker(\beta) \ \simeq \ \ext^1(E,V_S).$$
Pour d\'emontrer le lemme 4.4, il suffit de montrer que
$$\dim(\hom(\hom(E,\ko_X),\ext^1(E,\ko_X)) \ - \ \dim(\ext^1(E,V_S))$$
est sup\'erieur ou \'egal \`a 4 si \ \m{X=\P_2}, et sup\'erieur ou \'egal
\`a 3 si \ \m{X=\P_1\times\P_1}.

\vskip 0.8cm

\noindent{\bf Etape 3 : }{\em estimation de \ \
\m{\dim(\hom(\hom(E,\ko_X),\ext^1(E,\ko_X))) - \dim(\ext^1(E,V_S))}}

\medskip

On consid\`ere la suite exacte
$$0\lra\ko_S\ot L\lra V_S\lra V(L)\lra 0.$$
On en d\'eduit la suite exacte
$$\hom(E,V(L))\lra\ext^1(E,\ko_S\ot L)\lra\ext^1(E,V_S)\lra\ext^1(E,V(L)),$$
d'o\`u
$$\dim(\ext^1(E,V_S)) \ \leq \dim(\ext^1(E,\ko_S\ot L))+\dim(\ext^1(E,V(L))).$$
On commence par calculer \ \m{\dim(\ext^1(E,\ko_S\ot L))}. On a \
\m{\ext^2(E,\ko_S\ot L)=\lbrace 0\rbrace} \ par dualit\'e de Serre. Cette
derni\`ere raison implique que
$$\hom(E,\ko_S\ot L) \ \simeq \ \hom(\ki_S\ot L,\ko_S\ot L),$$
qui est de dimension \m{2s}. Comme \ \m{\chi(E,\ko_S\ot L)=\chi(E,\ko_S)=s},
on en d\'eduit que
$$ \dim(\ext^1(E,\ko_S\ot L)) \ = \ s.$$
On calcule maintenant \ \m{\dim(\ext^1(E,V(L)))}. On a une suite exacte
$$\ext^1(\ki_S\ot L,V(L))\lra\ext^1(E,V(L))\lra\ext^1(T,V(L)).$$
Donc
$$\dim(\ext^1(E,V(L))) \ \leq \ \dim(\ext^1(T,V(L)))+
\dim(\ext^1(\ki_S\ot L,V(L))).$$
On a \ \m{\ext^2(T,V(L))=\lbrace 0\rbrace}. En effet, par dualit\'e de Serre il
suffit de montrer que \Nligne \m{\hom(V(L),T\ot\omega_X)=\lbrace 0\rbrace}. Si 
ce n'\'etait pas le cas, en utilisant une section de \m{V(L)} et un
morphisme non nul \ \m{\omega_X\lra\ko_X} \ on obtiendrait une section non
nulle de $T$, donc un morphisme non nul \ \m{\ki_S\ot L\lra T} \ qui
d\'efinirait un endomorphisme de $E$ qui n'est pas une homoth\'etie. On a
aussi \ \m{\hom(T,V(L))=\lbrace 0\rbrace} \ car \m{V(L)} est localement libre.
Donc, en utilisant le th\'eor\`eme de Riemann-Roch on trouve
$$\dim(\ext^1(T,V(L))) \ = \ -\chi(T,V(L)) \ \ \ \ \ \ \
 \ \ \ \ \ \ \ \ \ \ \ \ \ $$
$$= \ c_1(L)^2+
(h^0(L^{-1})-1)(\frac{\omega_Xc_1(L)}{2}+\frac{c_1(L)^2}{2}+c_2-s).$$
Pour calculer \ \m{\dim(\ext^1(\ki_S\ot L,V(L)))} \ on consid\`ere les suites
exactes
$$0\lra\ki_S\ot L\lra L\lra\ko_S\ot L\lra 0,$$
$$0\lra L\lra\ko_X\ot H^0(L^{-1})^*\lra V(L)\lra 0.$$
De la seconde on d\'eduit que
$$\ext^1(L,V(L)) \ = \ext^2(L,V(L)) \ = \lbrace 0\rbrace,$$
d'o\`u il d\'ecoule avec la premi\`ere que
$$\ext^1(\ki_S\ot L,V(L)) \ \simeq \ \ext^2(\ko_S\ot L,V(L)) \ \simeq \
\hom(V(L),\ko_S\ot L)^*.$$
On a donc
$$\dim(\ext^1(\ki_S\ot L,V(L))) \ = \ s(h^0(L^{-1})-1).$$
On a donc
$$\dim(\ext^1(E,V_S)) \ \leq \ s+c_1(L)^2+(h^0(L^{-1})-1)
(\frac{\omega_Xc_1(L)}{2}+\frac{c_1(L)^2}{2}+c_2-s).$$
On obtient finalement
$$\dim(\hom(\hom(E,\ko_X),\ext^1(E,\ko_X)) - \dim(\ext^1(E,V_S)) \ \geq \
-\chi(T,\ki_S\ot L).$$

\vskip 0.8cm

\noindent{\bf Etape 4 : }{\em fin de la d\'emonstration du lemme 4.4}

\medskip

On a
$$\chi(T,\ki_S\ot L) \ = \chi(\ki_S\ot L,T) - \omega_Xc_1(L).$$
On a \ \m{\hom(\ki_S\ot L,T)=\ext^2(\ki_S\ot L,T)=\lbrace0\rbrace}, donc
\ \m{\chi(\ki_S\ot L,T) < 0}. On a d'autre part \ \m{\omega_Xc_1(L)\geq 3} \ si
\ \m{X=\P_2} \ et \ \m{\omega_Xc_1(L)\geq 2} \ si \
\m{X=\P_1\times\P_1}. On en d\'eduit \ \m{-\chi(T,\ki_S\ot L)\geq 4} \ si
\ \m{X=\P_2} \ et \ \m{-\chi(T,\ki_S\ot L)\geq 3} \ si
\m{X=\P_1\times\P_1}. Ceci ach\`eve la d\'emonstration du lemme 4.4.
\end{proof}

\bigskip

\noindent{\bf Remarque : }
On consid\`ere la suite spectrale \m{E_r^{p,q}} convergeant vers
\m{\ext^{p+q}(E,E)} d\'efinie par la filtration \ \m{T\subset E}. Les termes
\m{E_1^{p,q}} \'eventuellement non nuls sont repr\'esent\'es ci-dessous :
\[
\begin{array}{ccc}
\ext^2(T,\ki_S\ot L) & 0 & 0 \\
\ext^1(T,\ki_S\ot L) & \ext^2(T,T)=\lbrace 0\rbrace & 0 \\
0 & \ext^1(T,T)\oplus\ext^1(\ki_S,\ki_S) & 0 \\
0 & \hom(T,T)\oplus\C & \ext^1(\ki_S\ot L,T)\\
\end{array}
\]
On en d\'eduit un isomorphisme non canonique
$$\ext^1(E,E) \ \simeq \ \ext^1(T,T)\oplus\ext^1(\ki_S,\ki_S)\oplus
\ext^1(\ki_S\ot L,T)\oplus\ext^1(T,\ki_S\ot L).$$
Les trois premiers termes correspondent aux d\'eformations de $E$ obtenues en
d\'eformant $T$, $S$ et l'extension. Le troisi\`eme doit repr\'esenter des
d\'eformations en faisceaux sans torsion ou avec un \m{\dim(\hom(E,\ko_X))}
plus petit. Je n'ai pas trouv\'e de d\'emonstration r\'eellement
convaincante du th\'eor\`eme 4.2 utilisant cette suite spectrale.

\vskip 1.2cm

\begin{sub}{\bf Un exemple de vari\'et\'e de modules fins maximale
non projective}\end{sub}

On donne ici des exemples de vari\'et\'es de modules fins qui sont des cas
particuliers de ce qui est d\'ecrit au \para 5.1. Mais montre ici que ces
vari\'et\'es de modules fins sont maximales.

\vskip 0.8cm

\begin{subsub}
{Sous-sch\'emas finis de $\P_2$ de longueur \m{\frac{n(n+1)}{2}}}\end{subsub}

\begin{xprop}
Soient \m{n} un entier positif, $Z$ un sous-sch\'ema ferm\'e de \m{\P_2} de
longueur \m{\frac{n(n+1)}{2}}. Alors il existe une suite exacte
$$0\lra\ko(-n-1)\ot\C^{n}\lra\ko(-n)\ot\C^{n+1}\lra\ki_Z\lra 0$$
si et seulement si $Z$ n'est pas contenu dans une courbe de degr\'e \m{n-1}.
\end{xprop}

\begin{proof} Il est clair que si il existe une telle suite exacte on a \
\m{H^0(\ki_Z(n-1))=\lbrace 0\rbrace}, c'est-\`a-dire que $Z$ n'est pas contenu
dans une courbe de degr\'e \m{n-1}. R\'eciproquement, supposons que $Z$ ne soit
pas contenu dans une courbe de degr\'e \m{n-1}. On a alors
$$H^0(\ki_Z(n-1)) \ = \ H^0(\ki_Z\ot Q^*(n-1)) = H^0(\ki_Z(n-2)) = \lbrace 0
\rbrace,$$
et aussi
$$H^2(\ki_Z(n-1)) \ = \ H^2(\ki_Z\ot Q^*(n-1)) = H^2(\ki_Z(n-2)) = \lbrace 0
\rbrace$$
par dualit\'e de Serre. On a donc \ \m{H^1(\ki_Z(n-1))=\lbrace 0\rbrace}, car
\ \m{\chi(\ki_Z(n-1))=0}. La suite spectrale de Beilinson donne donc une suite
exacte
$$0\lra\ko(-n-1)\ot H^1(\ki_Z(n-2))\lra\ko(-n)\ot H^1(\ki_Z\ot Q^*(n-1))\lra
\ki_Z\lra 0.$$
Un calcul simple montre que \ \m{h^1(\ki_Z(n-2))=n}, \m{h^1(\ki_Z\ot Q^*(n-1))
=n+1}. \end{proof}

\bigskip

Soit \ \m{W=\hom(\ko(-n-1)\ot\C^{n},\ko(-n)\ot\C^{n+1})}. Sur l'espace
projectif

\noindent\m{\P=\P(W)} agit le groupe alg\'ebrique r\'eductif \
\m{SL(n)\times SL(n+1)}, et cette action se prolonge de mani\`ere \'evidente
\`a une action sur \m{\ko_{\P}(1)}. Un point de \m{\P}, correspondant
\`a un morphisme
$$f : \ko(-n-1)\ot\C^{n}\lra\ko(-n)\ot\C^{n+1}$$
est semi-stable pour cette action si et seulement si il est stable, et ceci
est vrai si et seulement si pour tout entier $p$ tel que \ \m{1\leq p\leq n},
et tous sous-espaces vectoriels de dimension $p$, \m{M_p\subset\C^{n}},
\m{N_p\subset\C^{n+1}}, \m{f(\ko(-n-1)\ot M_p)} n'est pas contenu dans
\m{\ko(-n)\ot N_p} (cf. \cite{dr2}). Si \m{\P^s} d\'esigne l'ouvert des
points stables de \m{\P}, il existe un quotient g\'eom\'etrique
$$N(3,n,n+1) \ = \ \P^s/(SL(n)\times SL(n+1)),$$
et c'est une vari\'et\'e projective lisse de dimension \m{n(n+1)}. Soient
\ \m{p_N}, \m{p_2} les projections \ \m{N(3,n,n+1)\times\P_2\lra
N(3,n,n+1)} \ et \m{N(3,n,n+1)\times\P_2\lra\P_2} \ respectivement.
Alors il existe sur \ \m{N(3,n,n+1)\times\P_2} \ un {\em morphisme
universel}
$$\Phi : p_2^*(\ko(-n-1))\ot p_N^*(M_n)\lra p_2^*(\ko(-n))\ot p_N^*(N_n),$$
\m{M_n} (resp. \m{N_n}) \'etant un fibr\'e vectoriel de rang $n$ (resp.
\m{n+1}) sur \m{N(3,n,n+1)}. On pose \ \m{\ke=\coker(\Phi)}.

Soit $\kx$ l'ensemble ouvert des classes d'isomorphisme de faisceaux
coh\'erents $E$ sur \m{\P_2}, de rang 1 et de classes de Chern 0 et
\m{\frac{n(n+1)}{2}}, v\'erifiant
$$H^0(E(n-1)) \ = \ H^1(E(n-1)) \ = \ H^2(E(n-2)) \ = \ \lbrace 0\rbrace.$$
D'apr\`es la suite spectrale de Beilinson, ce sont les faisceaux qui sont
isomorphes \`a des conoyaux de morphismes injectifs de $W$. Soit \m{\kx^s}
l'ensemble ouvert de faisceaux constitu\'e des classes d'isomorphisme de
conoyaux de morphismes injectifs stables de $W$. Alors, si
$U$ d\'esigne l'ouvert de \m{N(3,n,n+1)} constitu\'e des morphismes injectifs
(comme morphismes de faisceaux), \m{(U,\ke_U)} est une vari\'et\'e de
modules fins pour \m{\kx^s}.

\medskip

\begin{xprop}
1 - Soient $Z$ un sous-sch\'ema fini de longueur \m{\frac{n(n+1)}{2}} de
\m{\P_2}, et
$$f : \ko(-n-1)\ot\C^{n}\lra\ko(-n)\ot\C^{n+1}$$
un morphisme dont le conoyau est isomorphe \`a \m{\ki_Z}. Alors $f$ est stable.

\noindent 2 - Il existe des points $x$ de $U$ tels que \m{\ke_x} ait de la
torsion.

\noindent 3 - Si \ \m{n\geq 5} \ on a \ \m{U \ \not= \ N(3,n,n+1)}, et donc $U$
n'est pas projective.
\end{xprop}

\begin{proof} Soit
$$g : \ko(-n-1)\ot\C^{n}\lra\ko(-n)\ot\C^{n+1}$$
un morphisme injectif non stable. Alors il existe un entier $p$ tel que
\ \m{1\leq p\leq n} \ et des sous-espaces vectoriels \ \m{M_p\subset\C^{n}},
\m{N_p\subset\C^{n+1}} de dimension $p$ tels que
$$g(\ko(-n-1)\ot M_p) \ \subset \ \ko(-n)\ot N_p.$$
Soient
$$g' : \ko(-n-1)\ot M_p\lra\ko(-n)\ot N_p, \ \ \ \
g'' : \ko(-n-1)\ot (\C^{n}/M_p)\lra\ko(-n)\ot (\C^{n+1}/N_p)$$
les morphismes d\'eduits de $g$. On a alors une suite exacte
$$0\lra\ker(g'')\lra\coker(g')\lra E\lra\coker(g'')\lra 0.$$
Mais \m{\ker(g'')} est localement libre et \m{\coker(g')} de torsion (car de
rang nul). Donc \m{\ker(g'')} est nul, et \m{\coker(g')} est un sous-faisceau
de torsion non trivial de $E$. Ceci prouve 1- (c'est aussi une cons\'equence
de la proposition 5.1).

Pour d\'emontrer 2-, on consid\`ere un sous-sch\'ema fini $S$ de \m{\P_2}
de longueur
$$s \ = \ h^0(\ko(n-2))-1 \ = \ \frac{(n+1)(n-2)}{2}$$
contenu dans une unique courbe de degr\'e \m{n-2}. On a alors
$$H^1(\ki_S(n-2)) \ = \ H^2(\ki_S(n-2)) \ = \ \lbrace 0\rbrace.$$
Soit \ \m{\ell=\P(D)} \ une droite de \m{\P_2} ne rencontrant pas
le support de $S$. On consid\`ere une extension
$$0\lra\ko_\ell(-n-1)\lra E\lra\ki_S(-1)\lra 0$$
donn\'ee par un \'el\'ement $\sigma$ de \m{\ext^1(\ki_S(-1),\ko_\ell(-n-1))}.
On a
$$\ext^1(\ki_S(-1),\ko_\ell(-n-1)) \ \simeq \ \ext^1(\ko(-1),\ko_\ell(-n-1)) \
\simeq \ S^{n-2}D.$$
Le morphisme canonique
$$\ext^1(\ki_S(-1),\ko_\ell(-n-1))\ot H^0(\ko(n-2))\lra H^1(\ko_\ell(-2))$$
est l'accouplement canonique
$$S^{n-2}D\ot S^{n-2}V^*\lra\C.$$
Il en d\'ecoule qu'il est possible de choisir $\sigma$ tel que le morphisme
induit
$$H^0(\ki_S(n-2))\lra H^1(\ko_\ell(-2))$$
soit un isomorphisme. On a alors
$$H^0(E(n-1)) \ = \ H^1(E(n-1)) \ = \ H^2(E(n-2)) \ = \ \lbrace 0\rbrace.$$
Il en d\'ecoule que $E$ est isomorphe au conoyau d'un morphisme injectif
$$g : \ko(-n-1)\ot\C^{n}\lra\ko(-n)\ot\C^{n+1}.$$
Il reste \`a voir que $g$ est stable. On a vu dans la d\'emonstration de -1
que si $g$ n'est pas stable, $E$ contient un sous-faisceau de torsion non nul
qui est un quotient d'un fibr\'e du type \ \m{\ko(-n)\ot\C^{p}}. Ceci est
impossible car
$$\hom(\ko(-n),\ko_\ell(-n-1)) \ = \ \lbrace 0\rbrace.$$
Donc $g$ est stable.

La d\'emonstration de 3- est donn\'ee au \para 4.2.2. \end{proof}

\bigskip

\begin{xtheo}
La vari\'et\'e de modules fins $U$ est maximale.
\end{xtheo}

\begin{proof} Soit \m{U'} une vari\'et\'e de modules fins d\'efinie localement
contenant strictement $U$. Soient \ \m{x\in U'\backslash U} \ et $\kf$ un
faisceau universel d\'efini sur un voisinage $W$ de $x$.
Supposons que \ \m{h^i(\kf_x(n-1))=0} \ pour \ \m{i\geq 0}. Alors d'apr\`es la
proposition 5.1 (avec \m{E=\ko(-n-1)}, \m{G=\ko(-n)}, \m{F=\ko(-n+1)}),
on a \ \m{x\in U}. Ceci \'etant faux, $x$ est un point de
l'hypersurface $\kh$ de \m{W} constitu\'e des points $y$ tels que les
\m{h^i(\kf_y(n-1))} ne soient pas tous nuls. D'apr\'es le th\'eor\`eme 4.2,
l'ouvert de $\kh$ des points $y$ tels que \m{\kf_y} soit un faisceau d'id\'eaux
est dense, et $\kh$ est irr\'eductible (car dans
\m{\hilb^{n(n+1)/2}(\P_2)},
l'hypersurface constitu\'ee des sous-sch\'emas finis contenus dans au moins une
courbe de degr\'e \m{n-1} est irr\'eductible).
Il existe aussi un point $y$ de $\kh$
tel que $\kf_y$ soit le faisceau d'id\'eaux d'un sous-sch\'ema $S$ contenu dans
une unique courbe de degr\'e \m{n-1}. On a un diagramme commutatif avec
colonnes exactes
\[
\begin{array}{ccc}
0 & & 0 \cr
\vfl{}{} & & \vfl{}{} \\
H^0(Q^*(n-1)) & \hfl{i'}{} & H^0(Q^*(n-1)_{\mid S}) \\
\vfl{}{} & & \vfl{}{} \\
H^0(\ko(n-1))\ot V^* & \hfl{i}{} & H^0(\ko(n-1)_{\mid S})\ot V^* \\
\vfl{}{} & & \vfl{}{} \\
H^0(\ko(n)) & \hfl{i''}{} & H^0(\ko(n)_{\mid S}) \\
\vfl{}{} & & \vfl{}{} \\
0 & & 0 \\
\end{array}
\]
Soit $\sigma$ une \'equation de l'unique courbe $C$ de degr\'e \m{n-1}
contenant $S$. Alors
$$\ker(i)=V^*\ot\C\sigma.$$
Donc l'image dans \m{H^0(\ko(n))}
d'un \'el\'ement non nul de \m{\ker(i)} n'est pas nulle. Ceci implique que
\m{i'} est injective. On en d\'eduit que
$$H^0(\ki_S\ot Q^*(n-1)) \ = \ \lbrace 0\rbrace.$$
Il en d\'ecoule que
$$h^1(\ki_S(n-2)) \ = \ n, \ \ \ \ h^1(\ki_S\ot Q^*(n-1)) \ = \ n+1.$$

Soit
$$g : \ko(-n-1)\ot H^1(\ki_S(n-2))\lra \ko(-n)\ot H^1(\ki_S\ot Q^*(n-1))$$
le morphisme canonique. On va montrer qu'on peut choisir $S$ de telle sorte
que $g$ soit stable.
On a \ \m{\chi(\ko_C(-S)(n-1))=-1}. On peut donc choisir $S$ suffisamment
g\'en\'eral sur $C$ pour que \ \m{h^0(\ko_C(-S)(n-1))=0}. On a une suite exacte
$$0\lra\ko(1-n) \ \hfl{\sigma}{} \ \ki_S\lra\ko_C(-S)\lra 0.$$
On en d\'eduit des isomorphismes canoniques
$$H^1(\ki_S(n-2)) \ \simeq \ H^1(\ko_C(-S)(n-2)), \ \ \ \
H^1(\ki_S\ot Q^*(n-1)) \ \simeq \ H^1(\ko_C(-S)\ot Q^*(n-1)).$$
Soit $E$ l'unique extension non triviale
$$0\lra\ko_C(-S)\lra E\lra\ko(1-n)\lra 0.$$
Alors $E$ est un faisceau de rang 1 et de classes de Chern 0,
\m{\frac{n(n+1)}{2}}. De plus on a \ \m{h^i(E(n-1))=0} \ pour tout $i$ et \
\m{h^2(E(n-2))=0}. D'autre part le morphisme canonique
$$\ko(-n-1)\ot H^1(E(n-2))\lra \ko(-n)\ot H^1(E\ot Q^*(n-1))$$
est canoniquement isomorphe \`a $g$. Il en d\'ecoule que $g$ est injectif et \
\m{\coker(g)\simeq E}.

Montrons maintenant que $g$ est stable. Dans le cas contraire, il existe un
entier $p$ tel que \ \m{1\leq p\leq n-1} \ et des sous-espaces vectoriels de
dimension $p$, \m{M_p\subset H^1(E(n-2))}, \m{N_p\subset H^1(E\ot Q^*(n-1))},
tels que
$$g(\ko(-n-1)\ot M_p) \ \subset \ \ko(-n)\ot N_p.$$
Soient
$$g' : \ko(-n-1)\ot M_p\lra\ko(-n)\ot N_p, \ \ \ \
g'' : \ko(-n-1)\ot (\C^{n}/M_p)\lra\ko(-n)\ot (\C^{n+1}/N_p)$$
les morphismes d\'eduits de $g$. On a alors une suite exacte
$$0\lra\ker(g'')\lra\coker(g')\lra E\lra\coker(g'')\lra 0.$$
Comme dans la d\'emonstration de la proposition 4.6 on en d\'eduit que \m{g''}
est injectif. Puisque le support de \m{\coker(g')} est une courbe de degr\'e
$p$ et que \ \m{\coker(g')\subset\ko_C(-S)}, on a \ \m{p=n-1}. On a aussi
un morphisme surjectif
$$\coker(g'')\lra\ko(1-n).$$
Mais on a une suite exacte
$$0\lra\ko(-n-1)\lra\ko(-n)\ot\C^{2}\lra\coker(g'')\lra 0,$$
donc \m{\coker(g'')} est isomorphe \`a \ \m{\ko_\ell(-n)\oplus\ko(-n)} \ ou
\m{\ki_x(-n)} ($\ell$ \'etant une droite et $x$ un point de \m{\P_2}). Dans
aucun des deux cas il ne peut exister de morphisme surjectif 
\m{\coker(g'')\lra\ko(1-n)}. Donc $g$ est stable.

Soient $x$ le point de \m{U'} correspondant \`a \m{\ki_S}, et $y$ celui qui
correspond \`a $E$. On consid\`ere un germe de courbe lisse $S$ sur \m{U'}
d'origine $x$ et tel que \ \m{S\backslash\lbrace x\rbrace\subset U}. On note
\m{p_S} la projection \ \m{S\times\P_2\lra S}. En consid\'erant
le conoyau du morphisme injectif de fibr\'es sur \ $S\times\P_2$
$$\ko(-n-1)\ot R^1p_{S*}(\kf_S\ot\ko(n-2))\lra \ko(-n-1))\ot
R^1p_{S*}(\kf_S\ot Q^*(n-1))$$
on obtient un morphisme \ \m{\phi : S\lra U'} \ tel que \ \m{\phi(s)=s} \
si \ \m{s\not =x}, et \ \m{\phi(x)=y}. Ceci est absurde, donc \ \m{U'\subset U}.
\end{proof}

\vskip 0.8cm

\begin{subsub}{Morphismes stables non injectifs}\end{subsub}

On d\'emontre ici le 3- de la proposition 4.6. Rappelons d'abord quelques
r\'esultats de \cite{dr2}. La vari\'et\'e \m{N(3,n,n+1)} est isomorphe \`a la
vari\'et\'e de modules \m{M(n+2,-1,n+1)} des faisceaux semi-stables de rang
\m{n+2} et de classes de Chern \m{-1}, \m{n+1} sur l'espace projectif dual
\m{\P_2^*}. De tels faisceaux sont d'ailleurs stables et localement libres.
L'isomorphisme est d\'efini de la fa\c con suivante : \`a un morphisme
stable sur \m{\P_2}
$$g : \ko(-n-1)\ot\C^{n}\lra\ko(-n)\ot\C^{n+1}$$
on associe d'abord le morphisme sur \m{\P_2^*} :
$$\kov{g} : \ko(-1)\ot\C^{n}\lra Q^*\ot\C^{n+1}$$
(en utilisant le fait que \ \m{H^0(\P_2,\ko(1)) = H^0(\P_2^*,Q^*(1))}).
Le fibr\'e stable associ\'e \`a $g$ est \m{\coker(\kov{g})}. La partie 3- de la
proposition 4.6 \'equivaut alors \`a la

\medskip

\begin{xprop}
Si \ \m{n\geq 5} \ il existe des fibr\'es stables de rang \m{n+2} et de classes
de Chern \m{-1}, \m{n+1} sur \m{\P_2} qui ne sont pas de type de
d\'ecomposition g\'en\'erique rigide.
\end{xprop}

\begin{xlemm}
Soit $E$ un fibr\'e stable de rang $n$ et de classes de Chern $0$, $n$ sur
$\P_2$, qui n'est pas de type de d\'ecomposition g\'en\'erique
\m{(0,\ldots,0)} ou \m{(1,0,\ldots,0,-1)}. Alors
il existe des extensions non triviales
$$0\lra Q^*\lra\ke\lra E\lra 0.$$
Le fibr\'e $\ke$ obtenu est stable de rang \m{n+2} et de classes de Chern
\m{-1}, \m{n+1} et n'est pas de type de d\'ecomposition g\'en\'erique rigide.
\end{xlemm}

\begin{proof} Puisque $E$ est stable, on a \ \m{\hom(E,Q^*)=\ext^2(E,Q^*)=
\lbrace 0\rbrace}. On a donc
$$\dim(\ext^1(E,Q^*)) \ = \ -\chi(E,Q^*) \ = \ 2n \ > 0,$$
donc il existe bien des extensions non triviales.

Montrons maintenant que $\ke$ est stable. Il faut montrer que pour tout
sous-fibr\'e propre $\kf$ de $\ke$, on a \ \m{\mu(\kf)<0}. Soient \m{\kf''}
l'image de $\kf$ dans $E$, et \m{\kf'} le noyau de \ \m{\kf\lra\kf''}. On a
alors une suite exacte
$$0\lra\kf'\lra\kf\lra\kf''\lra 0.$$
Supposons d'abord que \m{\kf'} est non nul. On a alors \ \m{\mu(\kf')<0} \ et \
\m{\mu(\kf'')\leq 0}. Donc \ \m{\mu(\kf)<0}. On peut donc supposer que \m{\kf'}
est nul, et \ \m{\kf=\kf''}. Donc, $E$ \'etant stable, \m{\mu(\kf)} ne peut
\^etre positif ou nul que si \ \m{\kf=E}. Mais ceci est impossible car alors
l'extension serait triviale.

Il reste \`a montrer que $\ke$ est de type de d\'ecomposition g\'en\'erique non
rigide, c'est-\`a-dire que pour toute droite $\ell$ de \m{\P_2} on a \
\m{h^0(\ke_\ell(-1))>0}. On a une suite exacte
$$0\lra H^0(\ke_\ell(-1))\lra H^0(E_\ell(-1))\lra H^1(Q^*_\ell(-1))=\C,$$
et le r\'esultat d\'ecoule du fait que \ \m{h^0(E_\ell(-1))\geq 2}, \`a cause
du type de d\'ecomposition g\'en\'erique de $E$. \end{proof}

\medskip

La proposition 4.8 est donc une cons\'equence de la

\medskip

\begin{xprop}
Si \ \m{n\geq 5}, il existe un fibr\'e stable de rang $n$ et de classes de
Chern $0$, $n$ sur $\P_2$, qui n'est pas de type de d\'ecomposition
g\'en\'erique \m{(0,\ldots,0)} ou \m{(1,0,\ldots,0,-1)}.
\end{xprop}

\begin{proof} On montre d'abord qu'il suffit de trouver un fibr\'e semi-stable
$E$ de rang $n$ et de classes de Chern $0$, $n$ sur $\P_2$, qui n'est pas de
type de d\'ecomposition g\'en\'erique \m{(0,\ldots,0)} ou
\m{(1,0,\ldots,0,-1)}. En effet, si un tel $E$ existe, dans toute d\'eformation
compl\`ete de $E$ param\'etr\'ee par un germe de vari\'et\'e lisse $S$, les
fibr\'es semi-stables non-stables constituent une sous-vari\'et\'e de
codimension au moins \m{\frac{n(n-1)}{2}} de $S$
(cf. la d\'emonstration du th\'eor\`eme
4.10 de \cite{dr_lp}). D'autre part, pour toute droite $\ell$ de \m{\P_2},
le morphisme de restriction
$$\ext^1(E,E)\lra\ext^1(E_\ell,E_\ell)$$
est surjectif. Il en d\'ecoule (cf. \cite{br}) que la sous-vari\'et\'e
ferm\'ee de $S$ correspondant aux fibr\'es dont le type de d\'ecomposition
g\'en\'erique est diff\'erent de \m{(0,\ldots,0)} et \m{(1,0,\ldots,0,-1)} est
de codimension 6. Comme \ \m{\frac{n(n-1)}{2}>6} \ on peut d\'eformer $E$ en
fibr\'e stable de type de d\'ecomposition g\'en\'erique diff\'erent de
\m{(0,\ldots,0)} et \m{(1,0,\ldots,0,-1)}.

Pour toute famille compl\`ete de faisceaux
semi-stables de rang $n$ et de classes de Chern $0$, $n$ param\'etr\'ee par
une vari\'et\'e lisse $S$, les points de $S$ correspondant aux faisceaux
non localement libres constituent une sous-vari\'et\'e ferm\'ee de
codimension au moins \m{n-1}. Il en d\'ecoule que si \ \m{n>7} \ il suffit
de trouver un faisceau semi-stable $E$ de
rang $n$ et de classes de Chern $0$, $n$ sur $\P_2$, qui n'est pas de type
de d\'ecomposition g\'en\'erique \m{(0,\ldots,0)} ou \m{(1,0,\ldots,0,-1)}.

On va maintenant construire par r\'ecurrence sur $n$ un faisceau semi-stable
$E$ de rang $n$ et de classes de Chern $0$, $n$ sur $\P_2$, qui n'est pas de
type de d\'ecomposition g\'en\'erique \m{(0,\ldots,0)} ou
\m{(1,0,\ldots,0,-1)}, et qui est localement libre si \ \m{n\leq 7}. Soit $x$
un point de \m{\P_2}. On prend
$$E \ = \ (S^4Q^*)(-2)\oplus (n-5)\ki_x$$
si \ \m{n\not =6,7}. Pour \ \m{n=6}, on prend une extension non triviale
$$0\lra (S^4Q^*)(-2)\lra E\lra\ki_x\lra 0,$$
et pour \ \m{n=7} \ on prend \ \m{E=(S^6Q^*)(-3)}.
\end{proof}

\vskip 0.8cm

\begin{subsub}{Exemples de faisceaux simples de rang 1 \'equivalents}
\end{subsub}

On reprend les notations de la d\'emonstration du th\'eor\`eme 4.7. Soit $S$ un
sous-sch\'ema fini
de \m{\P_2} de longueur \m{\frac{n(n+1)}{2}} contenu dans une
unique courbe lisse $C$ de degr\'e \m{n-1}, tel que \
\m{h^0(\ko_C(-S)(n-1))=0}. On a une unique extension non triviale
$$0\lra\ko_C(-S)\lra E\lra\ko(1-n)\lra 0.$$
Ce faisceau d\'epend en fait uniquement de \ \m{L=\ko_C(S)} (et pas de $S$).
On notera donc \ \m{E=E_L}. Le faisceau \m{E_L} fait partie de la vari\'et\'e
de modules fins $U$ d\'efinie pr\'ec\'edemment et d'apr\`es la d\'emonstration
du th\'eor\`eme 4.7 on a
$$E_L \ \equiv \ \ki_S.$$
Donc si $S'$ est un autre diviseur de \m{C} ayant les m\^emes propri\'et\'es
que $S$ tel que  \m{\ko_C(S)\simeq\ko_C(S')} \ on a
$$\ki_S \ \equiv \ \ki_{S'}.$$
Il existe une sous-vari\'et\'e localement ferm\'ee de dimension au moins
\m{3n-3} de

\noindent\m{{\rm Hilb}^{n(n+1)/2}(\P_2)} constitu\'ee de tels $S'$.

Soit
$$\phi : \ext^1(E,E)\lra\ext^1(\ko_C(-S),\ko(1-n))$$
le morphisme canonique. Alors $\phi$ est surjectif et son noyau est
l'espace tangent en $E$ de \
\m{N=U\backslash{\rm Hilb}^{n(n+1)/2}(\P_2)}. En utilisant la dualit\'e de
Serre sur $C$ et \m{\P_2} on obtient un isomorphisme canonique
$$\ext^1(\ko_C(-S),\ko(1-n)) \ \simeq \ H^0(\ko(S)).$$
Soient \ \m{\sigma\in\P({H^0(\ko(S))})} \ et \m{(D,E)} un germe de courbe
lisse sur $U$ passant par $E$ tel que l'image par $\phi$ de la tangente en
$E$ \`a $D$ soit $\sigma$. Alors on montre ais\'ement que le point limite de
correspondant dans \m{{\rm Hilb}^{n(n+1)/2}(\P_2)} est le sous-sch\'ema
fini de $C$ d\'efini par $\sigma$.

\vskip 2.5cm

\section{Vari\'et\'es de modules fins et vari\'et\'es de modules de
morphismes}

Dans ce chapitre on \'etudie des vari\'et\'es de modules fins de faisceaux
coh\'erents sur \m{\P_2}. On utilise la th\'eorie des fibr\'es
exceptionnels, d\'ej\`a rappel\'ee au \para 3. Les r\'esultats peuvent
s'\'etendre sans difficult\'es aux autres surfaces pour lesquelles une
th\'eorie similaire existe, c'est-\`a-dire les {\em surfaces de Del Pezzo} (cf.
\cite{go}, \cite{go_ru}, \cite{ka1}, \cite{ka2}).

\vskip 1.2cm

\begin{sub}{\bf Vari\'et\'es de modules de morphismes}\end{sub}

\begin{subsub}{Modules de Kronecker}\end{subsub}

Les r\'esultats de ce chapitre sont d\'emontr\'es dans \cite{dr2}. Soient
\m{L} un espace vectoriel complexe de dimension finie avec \
\m{q=\dim(L)\geq 3}, $m$ et $n$ des entiers positifs. Les applications
lin\'eaires
$$L\ot \C^{m}\lra \C^{n}$$
sont appel\'ees des \m{L}-{\em modules de Kronecker}. Soit
$$W \ = \ \hom(L\ot \C^{m},\C^{n}).$$
Sur \m{W} op\`ere de mani\`ere \'evidente le groupe alg\'ebrique r\'eductif
$$G=(GL(m)\times GL(n))/\C^{*}.$$
L'action de \ \m{SL(m)\times SL(n)} \ sur \m{\P(W)} se lin\'earisant de
fa\c con \'evidente, on a une notion de point {\em (semi-)stable} de
\m{\P(W)} (au sens de la g\'eom\'etrie invariante). On montre que si
\m{f\in W}, \m{f} est semi-stable (resp.
stable) si et seulement si pour tous sous-espaces vectoriels \m{M'} de
\m{\C^{m}} et \m{N'} de \m{\C^{n}}, tels que \m{M'\not = \lbrace 0\rbrace},
\m{N'\not =\C^{n}}, et \ \m{f(L\ot M')\subset N'}, on a
$$\frac{\dim(N')}{\dim(M')}\geq\frac{n}{m}\ \ {\rm (resp.}
\ > \ {\rm)}.$$
Soit \m{W^{ss}} (resp. \m{W^s}) l'ouvert des points semi-stables (resp.
stables) de \m{W}. Alors il existe un bon quotient (resp. un quotient
g\'eom\'etrique)
$$N(q,m,n) = W^{ss}//G \ \ \ {\rm (resp. } \ \ N_s(q,m,n)=W^s/G \ {\rm)},$$
\m{N(q,m,n)} est projective, et \m{N_s(q,m,n)} est un ouvert lisse de
\m{N(q,m,n)}. De plus, si $m$ et $n$ sont premiers entre eux, on a \
\m{N(q,m,n)=N_s(q,m,n)}.

\vskip 0.8cm

\begin{subsub}{Applications aux vari\'et\'es de modules fins de faisceaux
coh\'erents}\end{subsub}

Soient \m{(E,G,F)} une triade de fibr\'es exceptionnels sur \m{\P_2}.
Soient $m$, $n$ des entiers positifs. On consid\`ere des morphismes
$$E\ot\C^{m}\lra G\ot\C^{n}.$$
On suppose que \ \m{n.rg(G)>m.rg(E)} (le cas \ \m{n.rg(G)<m.rg(E)} \ est
analogue). Soient \ \m{W=\hom(\hom(E,G)^*\ot\C^{m},\C^{n})}, \m{W_0} l'ouvert
de $W$ constitu\'e des morphismes injectifs (comme morphismes de faisceaux).
Le morphisme canonique universel
$$\phi : E\ot\C^{m}\lra G\ot\C^{n}$$
sur \ \m{W_0\times\P_2} \ est injectif. Soit $\kf$ son conoyau.
Alors pour tout point $x$ de \m{W_0}, $\kf$ est une d\'eformation compl\`ete de
\m{\kf_x}. C'est ce qui fait l'int\'er\`et de ces morphismes. On les utilise
dans \cite{dr2} pour d\'ecrire certaines vari\'et\'es de modules de faisceaux
semi-stables sur \m{\P_2} (pour une g\'en\'eralisation aux surfaces de Del
Pezzo, voir \cite{ka1}, \cite{ka2}).

\medskip

\begin{xprop}
Soient \m{(E,G,F)} une triade de fibr\'es exceptionnels sur \m{\P_2}, $H$
le fibr\'e exceptionnel noyau du morphisme d'\'evaluation \
\m{E\ot\hom(E,G)\lra G}. Soit $\ke$ un faisceau coh\'erent simple sur
\m{\P_2} tel que \ \m{h^i(\ke\ot F^*)=0} \ pour \ \m{i\geq 0}. Supposons
que \m{\chi(\ke,E)\leq 0} (resp. \m{\chi(\ke,E)>0}). Alors on a
$$H^2(\ke\ot E^*(-3)) \ = \ H^2(\ke\ot H^*(-3)) \ = \lbrace 0\rbrace,$$
et une suite exacte
$$0\lra E\ot H^1(\ke\ot E^*(-3)) \ \hfl{f}{} \ G\ot H^1(\ke\ot H^*(-3))\lra
\ke\lra 0$$
(resp.
$$H^1(\ke\ot E^*(-3)) \ = \ H^1(\ke\ot H^*(-3)) \ = \lbrace 0\rbrace,$$
et une suite exacte
$$0\lra \ke\lra E\ot H^2(\ke\ot E^*(-3)) \ \hfl{f}{} \ G\ot H^2(\ke\ot H^*(-3))
\lra 0).$$
Si de plus il existe un ensemble ouvert de classes d'isomorphisme de faisceaux
coh\'erents contenant celle de $\ke$ et admettant une vari\'et\'e de modules
fins d\'efinie localement, le morphisme $f$ est stable, et
\m{h^1(\ke\ot E^*(-3))} et \m{h^1(\ke\ot H^*(-3))} (resp. \m{h^1(\ke\ot
E^*(-3))} et \m{h^1(\ke\ot H^*(-3))}) sont premiers entre eux.
\end{xprop}

\begin{proof} Le diagramme de Beilinson de $\ke$ correspondant \`a la triade
\m{(E,G,F)} a l'allure suivante :
\[
\begin{array}{cccc}
E\ot H^2(\ke\ot E^*(-3)) & \hfl{d_1^{-2,2}}{} &
G\ot H^2(\ke\ot H^*(-3)) & 0 \\
 & & &  \\
E\ot H^1(\ke\ot E^*(-3)) & \hfl{d_1^{-2,1}}{} & G\ot H^1(\ke\ot H^*(-3)) & 0 \\
 & & & \cr 0 & & 0 & 0\\
 \end{array}
 \]
Il en d\'ecoule que \m{d_1^{-2,2}} est surjectif, \m{d_1^{-2,1}} injectif, et
qu'on a une suite exacte
$$0\lra\coker(d_1^{-2,1})\lra\ke\lra\ker(d_1^{-2,2})\lra 0.$$
Mais un calcul simple montre que
$$\ext^1(\ker(d_1^{-2,2}), \coker(d_1^{-2,1})) \ = \ \lbrace 0\rbrace.$$
On a donc un isomorphisme
$$\ke \ \simeq \ \ker(d_1^{-2,2})\oplus\coker(d_1^{-2,1}).$$
Puisque $\ke$ est simple, on a \ \m{\ke=\ker(d_1^{-2,2})} \ ou \
\m{\ke=\coker(d_1^{-2,1})}.

Supposons maintenant qu'il existe un ensemble ouvert $\kx$ de classes
d'isomorphisme de faisceaux coh\'erents contenant celle de $\ke$ et admettant
une vari\'et\'e de modules fins d\'efinie localement. En prenant un ensemble
plus petit que $\kx$, on se ram\`ene au cas o\`u $\kx$ admet une vari\'et\'e de
modules fins \m{(M,\kf)}. On peut supposer que \ \m{\chi(\ke,E)\leq 0} (l'autre
cas est analogue). Soit \m{M_0} l'ouvert de $M$ constitu\'e des points $x$
tels que
$$h^i(\kf_x\ot F^*) \ = \ h^j(\kf_x\ot E^*(-3)) \ = \
h^j(\kf_x\ot H^*(-3)) \ = \ 0$$
pour \ \m{i\geq 0} et \m{j=0,2}. Il existe alors un point \m{x_0} de \m{M_0}
tel que \ \m{\kf_{x_0}\simeq\ke}. Soient \ \m{p_M : M\times\P_2\lra M} \ la
projection et $U$ un voisinage ouvert de \m{x_0} dans \m{M_0} sur lequel les
fibr\'es vectoriels \m{R^1p_{M*}(\kf\ot E^*(-3))} et \m{R^1p_{M*}(\kf\ot
H^*(-3))} sont triviaux. Soient enfin
$$m \ = \ -\chi(\ke\ot E^*(-3)), \ \ \ \ n \ = \ -\chi(\ke\ot H^*(-3)),$$
\m{W=\hom(\hom(E,G)^*\ot\C^{m},\C^{n})}, \m{W_0} l'ouvert de $W$ constitu\'e
des morphismes injectifs de faisceaux dont le conoyau est isomorphe \`a un
\m{\kf_x}, avec \m{x\in U}. Le morphisme canonique
$$\pi : \P(W_0)\lra U$$
(associant \`a $f$ le point $x$ tel que \ \m{\kf_x=\coker(f)}) a une section
(d\'efinie par les trivialisations pr\'ec\'edentes). Il en d\'ecoule ais\'ement
que c'est un quotient g\'eom\'etrique par \ \m{SL(m)\times SL(n)}. D'apr\`es
\cite{mu_fo} (converse 1.13) il existe un \
\m{L \ \in \ Pic^{SL(m)\times SL(n)}(\P(W))} \ tel que \ \m{\P(W_0)
\subset\P(W)^{ss}(L)}. Mais \ \m{SL(m)\times SL(n)} \ n'ayant pas de
caract\`ere non
trivial, la seule lin\'earisation possible \'equivaut \`a la
lin\'earisation canonique (autrement dit il n'y a qu'une seule notion de
stabilit\'e pour les points de \m{\P(W)}). Donc $f$ est stable. La
d\'emonstration du reste de la proposition 5.1 est analogue \`a celle du
th\'eor\`eme G de \cite{dr_na}. \end{proof}

\medskip

On en d\'eduit la

\bigskip

\begin{xprop}
Soient \m{(E,G,F)} une triade de fibr\'es exceptionnels sur \m{\P_2}, $H$
le fibr\'e exceptionnel noyau du morphisme d'\'evaluation \
\m{E\ot\hom(E,G)\lra G}, $r$, \m{c_1}, \m{c_2} des entiers, tels que \
\m{r\geq 0}. Soit $\ke$ un faisceau coh\'erent de rang $r$ et de classes de
Chern \m{c_1}, \m{c_2}. On suppose que \m{\chi(F,\ke)=0}. Soient
$$m \ = \ \mid\chi(E(3),\ke)\mid, \ \ \ \ n \ = \ \mid\chi(H(3),\ke)\mid.$$
On suppose que $m$ et $n$ sont non nuls et premiers entre eux.

\noindent 1 - Soit \m{\ky_F} l'ensemble ouvert des classes d'isomorphismes de
faisceaux simples de rang $r$ et de classes de Chern \m{c_1}, \m{c_2} qui
sont des conoyaux de morphismes stables injectifs de faisceaux \
\m{E\ot\C^{m}\lra G\ot\C^{n}} \ si \ \m{\chi(E(3),\ke)<0}, et noyaux de tels
morphismes stables et surjectifs si \ \m{\chi(E(3),\ke)>0}. Alors il existe une
vari\'et\'e de modules fins \m{(U,\kf)} pour \m{\ky_F}, $U$ \'etant un ouvert
de \m{N(3rg(F),m,n)}.

\noindent 2 - Soit \m{\kx_F} l'ensemble ouvert des classes d'isomorphismes de
faisceaux simples de rang $r$ et de classes de Chern \m{c_1}, \m{c_2} tels que
\ \m{h^i(F^*\ot\ke)=0} \ pour \ \m{i\geq 0}. Soit $\kx$ un ensemble ouvert de
faisceaux de rang $r$ et de classes de Chern \m{c_1}, \m{c_2} admettant une
vari\'et\'e de modules fins d\'efinie localement. Alors on a
$$\kx\cap\kx_F \ \subset \ \ky_F.$$
\end{xprop}

\vskip 1.2cm

\begin{sub}{\bf Autres types de vari\'et\'es de modules de morphismes}
\end{sub}

Soient \m{(E,G,F)} une triade de fibr\'es exceptionnels sur \m{\P_2}.
Soient $m$, $n$, $p$ des entiers positifs. On consid\`ere des morphismes
$$E\ot\C^{m}\lra (G\ot\C^{n})\oplus(F\ot\C^{p}).$$
On suppose que \ \m{n.rg(G)+p.rg(F)>m.rg(E)} (le cas \
\m{n.rg(G)+p.rg(F)<m.rg(E)} \ est analogue). Soient
$$W \ = \ \hom(E\ot\C^{m},(G\ot\C^{n})\oplus(F\ot\C^{p})),$$
\m{W_0} l'ouvert de $W$ constitu\'e des morphismes injectifs (comme morphismes
de faisceaux). Le morphisme canonique universel
$$\phi : E\ot\C^{m}\lra (G\ot\C^{n})\oplus(F\ot\C^{p})$$
sur \ \m{W_0\times\P_2} \ est injectif. Soit $\kf$ son conoyau. Alors pour
tout point $x$ de \m{W_0}, $\kf$ est une d\'eformation compl\`ete de \m{\kf_x}.
On peut donc envisager de d\'ecrire des vari\'et\'es de modules fins de
faisceaux en utilisant de tels morphismes. C'est ce qui est fait dans
\cite{dr5} pour le cas des vari\'et\'es de modules {\em extr\'emales} de
faisceaux stables. En g\'en\'eral, il faut pouvoir construire des quotients
d'ouverts ad\'equats de $W$ par le groupe \
$$G \ = \ \Aut(E\ot\C^{m})\times\Aut((G\ot\C^{n})\oplus(F\ot\C^{p}))$$
qui n'est pas r\'eductif. C'est ce qui est fait dans \cite{dr_tr}, \cite{dr6}
et \cite{dr7}.

Le groupe $G$ poss\`ede un sous-groupe normal unipotent maximal $\gamma$
\'evident, isomorphe au groupe additif \m{\hom(G\ot\C^{n}, F\ot\C^{p})}.
Soient \m{\lambda}, \m{\mu} des nombres rationnels positifs tels que \
\m{\lambda n+\mu p =1}, et
$$f : E\ot\C^{m}\lra (G\ot\C^{n})\oplus(F\ot\C^{p})$$
un morphisme. On dit $f$ est {\em semi-stable} (resp. {\em stable})
relativement \`a \m{(\lambda,\mu)} si pour tous sous-espaces vectoriels \
\m{M'\subset\C^{m}}, \m{N'\subset\C^{n}}, \m{P'\subset\C^{p}}, avec \m{M'} non
nul et \m{N'\not =\C^{n}} ou \m{P'\not =\C^{p}}, et tout \m{f'\in\gamma.f} tel
que \ \m{f'(E\ot M')\subset(G\ot N')\oplus(F\ot P')}, on a
$$\lambda\dim(N')+\mu\dim(P') \ \geq \ \frac{\dim(M')}{m} \ \ \ \
{\rm (resp. \ \ }>{\rm \ \ )}.$$
La paire \m{(\lambda,\mu)} s'appelle une {\em polarisation} de l'action de $G$
sur $W$. On note \m{W^{ss}} (resp. \m{W^s}) l'ouvert de $W$ constitu\'e des
morphismes semi-stables (resp. stables) relativement \`a \m{(\lambda,\mu)}.
Pour certaines valeurs de \m{(\lambda,\mu)} on sait construire un bon quotient
\m{W^{ss}//G} (resp. un quotient g\'eom\'etrique \m{W^s/G}).

Soient $r$ le rang et \m{c_1}, \m{c_2} les classes de Chern des conoyaux de
morphismes injectifs de $W$. Soit \ \m{\chi = \frac{c_1(c_1+3)}{2}+r-c_2} \ la
caract\'eristique d'Euler-Poincar\'e de ces conoyaux. On suppose que $r$,
\m{c_1} et $\chi$ sont premiers entre eux. Soient \m{W^s/G} un quotient
g\'eom\'etrique construit dans \cite{dr7} (cf. th\'eor\`eme 5.6 de cet
article), et \m{M_0} l'ouvert correspondant aux morphismes injectifs. Alors,
en utilisant la proposition 2.4 de \cite{dr7} on voit ais\'ement qu'il existe
un {\em faisceau universel} sur \ \m{M_0\times\P_2}. On obtient ainsi une
vari\'et\'e de modules fins de faisceaux de rang $r$ et de classes de Chern
\m{c_1}, \m{c_2}.

Dans le cas des vari\'et\'es de modules extr\'emales \m{M(r,c_1,c_2)}
\'etudi\'ees dans \cite{dr5} (dans le cas o\`u $r$, \m{c_1} et \m{\chi} sont
premiers entre eux), la vari\'et\'e \m{W^s/G} est isomorphe \`a
\m{M(r,c_1,c_2)} pour une polarisation donn\'ee \m{(\lambda,\mu)}. En faisant
varier la polarisation on obtient des modifications de \m{M(r,c_1,c_2)} qui
sont d'autres vari\'et\'es de modules fins de faisceaux de rang $r$ et de
classes de Chern \m{c_1}, \m{c_2}. Un tel exemple est donn\'e au \para 5.3.2.

\newpage

\begin{sub}{\bf Exemples}\end{sub}

\begin{subsub}{Fibr\'es simples non admissibles}\end{subsub}

Soit \m{n>1} un entier. On consid\`ere des morphismes de fibr\'es vectoriels
sur \m{\P_2}
$$\ko\ot\C^{2n-1}\lra Q\ot\C^{2n+1}.$$
D'apr\`es \cite{dr2}, les morphismes stables de ce type sont injectifs, et
leurs conoyaux sont les fibr\'es stables de rang \m{2n+3} et de classes de
Chern \m{2n+1}, \m{(n+1)(2n+1)}. On a en fait un isomorphisme
$$M(2n+3,2n+1,(n+1)(2n+1)) \ \simeq \ N(3,2n-1,2n+1).$$

\bigskip

\begin{xprop}
Il existe des morphismes injectifs non stables
$$\phi : \ko\ot\C^{2n-1}\lra Q\ot\C^{2n+1}$$
tels que \m{\coker(\phi)} soit sans torsion et simple.
\end{xprop}

\begin{proof} Puisque la sous-vari\'et\'e de \m{\hom(V^*,\C^{n+1})}
constitu\'ee des applications de rang inf\'erieur ou \'egal \`a 1 est de
codimension $2n$ il existe une application lin\'eaire
$$\phi'' : \C^{n}\lra\hom(V^*,\C^{n+1})$$
telle que pour tout \ \m{u\in\C^{n}} \ non nul \m{\phi''(u)} soit de rang au
moins 2. On peut m\^eme choisir \m{\phi''} de telle sorte que le morphisme
de fibr\'es associ\'e
$$\phi'' : \ko\ot\C^{n}\lra Q\ot\C^{n+1}$$
soit stable et injectif. Il en d\'ecoule que \m{\phi''} ne contient pas
de sous-module de Kronecker du type
$$V^*\ot\C\lra\C,$$
et ne poss\`ede donc pas de modules de Kronecker quotients du type
$$V^*\ot\C^{n-1}\lra\C^{n}.$$
Soit
$$\phi' : \ko\ot\C^{n-1}\lra Q\ot\C^{n}$$
un morphisme de fibr\'es stable et injectif. Soient \ \m{E'=\coker(\phi')} \ et
\ \m{E''=\coker(\phi'')}. Alors on a
$$\hom(E'',E') \ = \ \ext^2(E'',E') \ = \ \lbrace 0\rbrace.$$
La seconde \'egalit\'e d\'ecoule de la stabilit\'e de \m{E'} et \m{E''}. Pour
montrer la premi\`ere, on remarque qu'un morphisme non nul \ \m{E''\lra E'} \
provient d'un morphisme de modules de Kronecker \ \m{\phi''\lra\phi'},
c'est-\`a-dire d'applications lin\'eaires \ \m{f:\C^{n}\lra\C^{n-1}},
\m{g:\C^{n+1}\lra\C^{n}} \ telles qu'on ait un diagramme commutatif
\[
\begin{array}{ccc}
V^*\ot\C^{n} & \hfl{\phi''}{} & \C^{n+1} \\
\vfl{I_{V^*}\ot f}{} & & \vfl{g}{} \\
V^*\ot\C^{n-1} & \hfl{\phi'}{} & \C^{n} \\
\end{array}
\]
Les stabilit\'es de \m{\phi''} et \m{\phi'} entra\^inent que $f$ et $g$ sont
surjectives, et \m{\phi'} est un module de Kronecker quotient de \m{\phi''},
ce qui est impossible. Donc
$$\dim(\ext^1(E'',E')) \ = \ -\chi(E'',E') \ = n^2 > 0.$$
Il existe donc des extensions non triviales
$$0\lra E'\lra E\lra E''\lra 0.$$
Le fibr\'e $E$ est isomorphe au conoyau d'un morphisme injectif
$$\phi : \ko\ot\C^{2n-1}\lra Q\ot\C^{2n+1}$$
qui est une extension de \m{\phi''} par \m{\phi'}, et n'est donc pas stable.
On v\'erifie ais\'ement que $E$ est simple. \end{proof}

\bigskip

On obtient donc des fibr\'es simples qui se d\'eforment en fibr\'es
stables mais qui ne peuvent pas faire partie d'une vari\'et\'e de modules fins
de faisceaux simples (d'apr\`es la proposition 5.2 avec \ \m{E=\ko} \ et \
\m{G=Q}). 

\vskip 0.8cm

\begin{subsub}{Exemples non triviaux de vari\'et\'es de modules fins
projectives de faisceaux simples}\end{subsub}

Soit \ \m{n>3} \ un entier. On consid\`ere des morphismes de fibr\'es
sur $\P_2$ du type
$$\phi : \ko(-2)\ot\C^{n-1}\lra\ko(-1)\ot\C^{2n-3}.$$
Si un tel $\phi$ est injectif (comme morphisme de faisceaux) on a, avec \
\m{\ke=\coker(\phi)},
$$rg(\ke)=n-2, \ \ \ c_1(\ke)=1, \ \ \ c_2(\ke)=n.$$
D'apr\`es la suite spectrale de Beilinson, un faisceau coh\'erent $E$ est
isomorphe au conoyau d'un morphisme injectif $\phi$ si et seulement si on a
$$\chi(E) \ = \ h^0(E) \ = \ h^2(E(-1)) \ = \ 0.$$
Dans ce cas on a \ \m{h^1(E(-1))=n-1}, \m{h^1(E\ot Q^*)=2n-3} \ et $\phi$ est
isomorphe au morphisme canonique
$$\ko(-2)\ot H^1(E(-1))\lra\ko(-1)\ot H^1(E\ot Q^*).$$
Soit
$$W_n \ = \ \hom(V\ot\C^{n-1},\C^{2n-3}).$$

\bigskip

\begin{xprop}
Soit \ \m{\phi\in W^{ss}_n}. Alors $\phi$ est injectif comme morphisme de
faisceaux.
\end{xprop}

\begin{proof} On voit $\phi$ comme une application lin\'eaire \
\m{V\ot\C^{n-1}\lra\C^{2n-3}}. Soit \ \m{x\in\P_2} \ tel que \m{\phi_x} ne
soit pas injectif. Il existe donc un \'el\'ement non nul $u$ de \m{\C^{n-1}}
tel que \ \m{\phi(x\ot u)=0}. Il en d\'ecoule que
$$\dim(\phi(V\ot u)) \ = \ 2$$
(la dimension ne peut pas \^etre $1$ car $\phi$ est semi-stable). On a alors un
diagramme commutatif
\[
\begin{array}{ccc}
V & \lra & V/x \\
\vfl{\alpha}{} & & \vfl{}{} \\
V\ot\C^{n-1} & \hfl{\phi}{} & \C^{2n-3} \\
\end{array}
\]
les fl\`eches verticales \'etant injectives, et $\alpha$ associant
\m{v\ot u} \`a
$v$. Soit \ \m{H\subset\C^{n-1}} \ le sous-espace vectoriel engendr\'e par les
$u$ comme pr\'ec\'edemment. Il existe donc une base \m{(u_1,\ldots,u_m)} de
$H$ et des points \m{x_1,\ldots,x_m} de \m{\P_2} tels qu'on ait un
diagramme commutatif
\[
\begin{array}{ccc}
V\ot H & \hfl{\phi_0}{} & \som_{1\leq i\leq m}(V/x_i) \cr
\vfl{\alpha_0}{} & & \vfl{\beta}{} \cr
V\ot\C^{n-1} & \hfl{\phi}{} & \C^{2n-3} \\
\end{array}
\]
o\`u \m{\alpha_0} est l'inclusion, et \m{\phi_0} est d\'efini par
$$\phi_0(\sigg_{1\leq i\leq m}v_i\ot u_i) \ = \ \sigg_{1\leq i\leq m}
\kov{v_i}$$
pour tous \m{v_i\in V}, \m{\kov{v_i}} d\'esignant l'image de \m{v_i} dans
\m{V/x_i}. La semi-stabilit\'e de $\phi$ entra\^ine que deux cas seulement
peuvent se produire :

\medskip

\noindent (i) $\beta$ est injective,

\noindent (ii) \m{\ker(\beta)} est de dimension $1$ et \ \m{H=\C^{n-1}}.

\medskip

\noindent Dans le cas (i), l'ensemble des points de \m{\P_2} o\`u $\phi$
n'est pas injectif est le m\^eme que celui o\`u \m{\phi_0} ne l'est pas, et
c'est exactement \m{\lbrace x_1,\ldots,x_m\rbrace}, qui est fini.

On peut donc supposer qu'on est dans le cas (ii), et que \m{(u_1,\ldots,u_m)}
est la base canonique de \m{\C^{n-1}}. Le module de Kronecker $\phi$
est donc \`a isomorphisme pr\`es de la forme
$$V\ot\C^{n-1}\lra\bigg(\som_{1\leq i\leq n-1}V/x_i\bigg)/D$$
$$(v_1,\ldots,v_{n-1})\longmapsto q(\kov{v_1},\ldots,\kov{v_{n-1}}),$$
$D$ \'etant une droite de \m{\som_{1\leq i\leq {n-1}}V/x_i}, \m{\kov{v_i}} la
classe de \m{v_i} dans \m{V/x_i} et
$$q :  \som_{1\leq i\leq n-1}V/x_i\lra
\bigg(\som_{1\leq i\leq n-1}V/x_i\bigg)/D$$
la projection. Posons
$$D \ = \ \C.\sigg_{1\leq i\leq n-1}\kov{w_i},$$
avec \m{w_i\in V}. Alors on a \ \m{\kov{w_i}\not = 0}, car sinon la restriction
de $\phi$
$$V\ot(\som_{j\not =i}\C u_i)\lra\bigg(\som_{i\not =j}V/x_i\bigg)/D$$
est un sous-module de Kronecker de $\phi$ contredisant sa semi-stabilit\'e.
Soit \m{z_i} une \'equation de la droite de \m{\P_2} contenant \m{w_i} et
\m{x_i}. On voit sans peine que le conoyau $E$ du morphisme injectif
$$\ko(-1) \ \hfl{(z_1,\ldots,z_{n-1})}{} \ \som_{1\leq i\leq n-1}\ki_{x_i}$$
v\'erifie
$$rg(E) \ = \ n-2, \ \ \ \ c_1(E) \ = \ 1, \ \ \ \ c_2(E) \ = \ n,$$
$$h^0(E) \ = \ h^2(E(-1)) \ = \ 0,$$
et que son module de Kronecker est isomorphe \`a $\phi$, qui est donc injectif
comme morphisme de faisceaux. \end{proof}

\bigskip

Avec les notations de la proposition 5.2 on obtient donc une vari\'et\'e de
modules fins \m{(N(3,n-1,2n-3),\kf)} pour \m{\ky_\ko}, et \m{N(3,n-1,2n-3)}
\'etant projective, cette vari\'et\'e de modules fins est maximale. Pr\'ecisons
que \m{\ky_\ko} contient des faisceaux ayant de la torsion, ainsi que des
faisceaux stables. Le groupe de Picard de \m{N(3,n-1,2n-3)} est isomorphe \`a
$\Z$, tandis que celui de la vari\'et\'e de modules correspondante
\m{M(n-2,1,n)} est isomorphe \`a $\Z^2$ (cf. \cite{dr4}).

\medskip

\noindent{\bf Remarque : } Soient $p$ un entier tel que
$$2n-2 \ \leq \ p \ < \ \frac{3+\sqrt{5}}{2}(n-1),$$
et
$$r \ = \ p-n, \ \ \ \ c_1 \ = \ 2n-2-p, \ \ \ \ c_2 \ = \
\frac{p(p-1)}{2}+\frac{n(n+1)}{2}-2(n-1)p.$$
Alors d'apr\`es \cite{dr2}, les faisceaux (semi-)stables de rang $r$ et de
classes de Chern \m{c_1}, \m{c_2} sont pr\'ecis\'ement les conoyaux des
morphismes (semi-)stables
$$\ko(-2)\ot\C^{n-1}\lra\ko(-1)\ot\C^{p}$$
et on a un isomorphisme \ \m{M(r,c_1,c_2) \ \simeq \ N(3,n-1,p)}.

\vskip 0.8cm

\begin{subsub}{Exemple de vari\'et\'e de modules fins de faisceaux simples sans
torsion, projective et contenant des faisceaux instables}\end{subsub}

On consid\`ere ici des faisceaux de rang 6 et de classes de Chern -3,8 sur
\m{\P_2}.
Dans ce cas on a \ \m{\chi=-2}, donc $r$, \m{c_1} et $\chi$ sont premiers entre
eux. On va d\'ecrire une vari\'et\'e de modules fins \m{M_1} de faisceaux de
rang 6 et de classes de Chern -3,8, qui est une modification de \m{M(6,-3,8)}.
On utilise les vari\'et\'es de modules de morphismes d\'ecrites au \para 5.2.
On consid\`ere des morphismes de fibr\'es
$$\ko(-3)\ot\C^{2}\lra\ko(-2)\oplus(\ko(-1)\ot\C^{7}).$$
Soient $W$ l'espace vectoriel de ces morphismes, $\lambda$, $\mu$ des nombres
rationnels positifs tels que \ \m{\lambda+7\mu=1}, d\'efinissant une
polarisation de l'action du groupe
$$G \ = \ GL(2)\times\Aut(\ko(-2)\oplus(\ko(-1)\ot\C^{7}))$$
sur $W$. On sait (d'apr\`es \cite{dr_tr}, \cite{dr6}, \cite{dr7}) construire
un quotient g\'eom\'etrique
$$M(\lambda,\mu)=W^s/G$$
d\`es que
$$\rho \ = \ \frac{\lambda}{\mu} > 3.$$
D'apr\`es \cite{dr5}, si \ \m{\rho=3+\epsilon}, avec $\epsilon$ positif
suffisamment petit, les morphismes stables du type pr\'ec\'edent sont injectifs
et leurs conoyaux sont les faisceaux stables de rang 6 et de classes de Chern
-3,8. On a en fait un isomorphisme \ \m{M(3,-6,8)\simeq M(\lambda,\mu)}.

On montre ais\'ement qu'il n'existe qu'une seule valeur de \m{\rho>3} pour
laquelle il existe des morphismes semi-stables non stables. C'est \ \m{\rho=7}.
Pour \ \m{\rho<7}, on a

\noindent\m{M(\lambda,\mu)=M(3,-6,8)}, et pour \ \m{\rho>7} \
on obtient une nouvelle vari\'et\'e de modules fins, not\'ee \m{M_1},
qu'on va d\'ecrire, munie d'un faisceau universel $\ke$.

Tout d'abord on v\'erifie ais\'ement qu'un morphisme stable pour une des
polarisations pr\'ec\'edentes est n\'ecessairement injectif. Soit
$$\phi=(\phi_1,\phi_2) : \ko(-3)\ot\C^{2}\lra\ko(-2)\oplus(\ko(-1)\ot\C^{7})$$
un morphisme stable. Alors \m{\phi_1} est non nul. Les morphismes \m{\phi}
peuvent donc \^etre de deux types, selon le rang de l'application induite par
\m{\phi_1}, \m{f_1 : \C^{2}\lra H^0(\ko(1))}.

\medskip

\noindent{\em Type 1 : } \m{f_1} est de rang 2. Il existe dans ce cas un
unique point $x$ de \m{\P_2} tel que \ \m{\imm(\phi_1)=\ki_x(-2)}, et on a
une suite exacte
$$0\lra\ko(-4)\lra\ko(-3)\ot\C^{2} \ \hfl{\phi_1}{} \ \ko(-2)\lra\C_x(-2)
\lra 0.$$
On a aussi une suite exacte
$$0\lra\ko(-4) \ \hfl{\beta}{} \ \ko(-1)\ot\C^{7}\lra\coker(\phi)\lra\C_x(-2)
\lra 0$$
(le morphisme $\beta$ \'etant la restriction de \m{\phi_2} \`a
\m{\ker(\phi_1)}). La stabilit\'e de $\phi$ entra\^ine que $\beta$ induit une
injection \ \m{{\C^{7}}^*\lra H^0(\ki_x(3))}. La $G$-orbite de $\phi$ est
enti\`erement d\'etermin\'ee par $x$ et par l'image de cette injection.

\medskip

\noindent{\em Type 2 : } \m{f_1} est de rang 1. Soit $\ell$ la droite de
\m{\P_2} d\'efinie par l'image de \m{f_1}. On a alors une suite exacte
$$0\lra\ko(-3)\lra\ko(-3)\ot\C^{2} \ \hfl{\phi_1}{} \ \ko(-2)\lra\ko_\ell(-2)
\lra 0.$$
On montre dans ce cas que la semi-stabilit\'e de $\phi$ implique que le
morphisme
$$\C^{2}\lra H^0(\ko_\ell(2))\ot\C^{7}$$
d\'eduit de \m{\phi_2} est injectif. La $G$-orbite de $\phi$ est alors
enti\`erement d\'etermin\'ee par la \m{GL(7)}-orbite de l'image de ce
morphisme.

\vskip 2cm

\noindent{\bf Description de \m{M_1} et des faisceaux correspondants}

Ces faisceaux sont les \m{\coker(\phi)}, $\phi$ \'etant de type 1.
On a donc un isomorphisme
$$M_1 \ \simeq \ {\bf G},$$
{\bf G} \'etant le fibr\'e en grassmanniennes sur \m{\P_2} dont la fibre
en $x$ est \m{{\bf Gr}(7, H^0(\ki_x(3)))}. L'intersection \
\m{M_1\cap M(6,-3,8)} \
au dessus du point $x$ de \m{\P_2}, est l'ouvert de la grassmannienne
\m{{\bf Gr}(7, H^0(\ki_x(3)))}
constitu\'e des \ \m{K\subset H^0(\ki_x(3))} \ ne contenant aucun sous-espace
vectoriel de la forme \m{s.H^0(\ko(2))}, $s$ \'etant une section non nulle de
\m{\ki_x(1)}.

On donne maintenant une description plus pr\'ecise des faisceaux param\'etr\'es
par $M_1$. Soient \ \m{x\in\P_2} \ et \ \m{K\in {\bf Gr}(7, H^0(\ki_x(3)))}.
On en d\'eduit un morphisme injectif
$$\beta_K : \ko(-4)\lra\ko(-1)\ot\C^{7}$$
nul en $x$. Soit \ \m{\kf_K=\coker(\beta_K)}. Du fait que \m{\beta_K}
s'annule en $x$ on a
$$\dim(\ext^1(\C_x,\kf_K)) \ = \ 1.$$
On obtient donc le faisceau \m{\ke_K}, unique extension non triviale
$$0\lra\kf_K\lra\ke_K\lra\C_x\lra 0.$$
Ce faisceau est sans torsion et a au plus deux points singuliers. C'est le
faisceau correspondant au point $K$ de $M_1$.

On obtient ainsi une vari\'et\'e de modules fins de faisceaux simples sans
torsion, projective, et qui n'est pas une vari\'et\'e de modules de faisceaux
stables.

\end{document}